\documentclass{article}
\usepackage{amsthm,thmtools}
\usepackage{amsmath}
\usepackage{setspace}
\usepackage{amsfonts}
\usepackage{mathtools}
\usepackage{wasysym}
\DeclareMathAlphabet{\mathbbb}{U}{bbold}{m}{n}
\usepackage{mathalfa} % just load the package, no options
\DeclareMathAlphabet{\mathdutchcal}{U}{dutchcal}{m}{n}
\usepackage{enumitem}% http://ctan.org/pkg/enumitem
\usepackage{cite}

\usepackage{stmaryrd}%\sslash
\usepackage[depth=subsection]{bookmark}

\theoremstyle{plain}
\newtheorem{prop}{Proposition}[section]
\newtheorem{theorem}[prop]{Theorem}
\newtheorem{maintheorem}[prop]{Main Theorem}
\newtheorem{corollary}[prop]{Corollary}

\theoremstyle{definition}
\newtheorem{df}[prop]{Definition}
\newtheorem{exa}[prop]{Example}
\newtheorem{lemma}[prop]{Lemma}
\newtheorem{remark}[prop]{Remark}

\newtheorem{construction}[prop]{Construction}
\newtheorem{variant}[prop]{Variant}
\newtheorem{recoll}[prop]{Recollection}
\newtheorem{notation}[prop]{Notation}
\newtheorem{alert}[prop]{Caveat}

\newtheorem{pb}[prop]{Problem}
\newtheorem{inputdata}[prop]{Input}

\theoremstyle{plain}
\newtheorem*{theorem*}{Theorem}

\usepackage{url}
\usepackage{scalerel}
\usepackage{latexsym, amssymb,amscd,mathrsfs,stackrel}
\usepackage{tikz-cd,tikz}
\usetikzlibrary{backgrounds}
\usepackage{graphicx} %插入图片的宏包
\usepackage{float} %设置图片浮动位置的宏包
\usepackage{subfigure} %插入多图时用子图显示的宏包
\usetikzlibrary{arrows}
\usepackage{geometry}

\newcommand{\infcat}{$\infty$-category}
\newcommand{\infcats}{$\infty$-categories}
\DeclareMathOperator{\map}{Map}
\DeclareMathOperator{\dalg}{DAlg}
\DeclareMathOperator{\calg}{CAlg}
\DeclareMathOperator{\daff}{dAff}

\DeclareMathOperator{\dArt}{dArtin}
\DeclareMathOperator{\dSt}{dSt}

\DeclareMathOperator{\Sp}{Sp}
\DeclareMathOperator{\vect}{Vect}
\DeclareMathOperator{\scr}{SCR}
\DeclareMathOperator{\alg}{Alg}
\DeclareMathOperator{\m}{Mod}
\DeclareMathOperator{\fun}{Fun}
\DeclareMathOperator{\QC}{QC}
\DeclareMathOperator{\aperf}{APerf}
\DeclareMathOperator{\perf}{Perf}
\DeclareMathOperator{\fil}{Fil}
\DeclareMathOperator{\mon}{CMon}
\DeclareMathOperator{\grp}{CGrp}

\DeclareMathOperator{\gr}{Gr}
\DeclareMathOperator{\tot}{Tot}
\DeclareMathOperator{\adic}{adic}

\DeclareMathOperator{\frob}{Frob}
\DeclareMathOperator{\colim}{colim}

\DeclareMathOperator{\Catoo}{\widehat{\mathcal{C}\mathrm{at}}_{\infty}}

\newcommand{\infcohfil}[2]{F\mathbbb{\Pi}\mathstrut_{#1/#2}}
\newcommand{\infcohcplfil}[2]{\widehat{F\mathbbb{\Pi}}\mathstrut_{#1/#2}}

\newcommand{\infeqtrcohcpl}[3]{\widehat{F_{ #3}\mathbbb{\Pi}}\mathstrut_{#1/#2}}
\newcommand{\infoneexcohcpl}[2]{F_{[0,1]}\mathbbb{\Pi}\mathstrut_{#1/#2}}
\newcommand{\inftwoexcohcpl}[2]{F_{[0,2]}\mathbbb{\Pi}\mathstrut_{#1/#2}}
\newcommand{\infcohcpl}[2]{\widehat{\mathbbb{\Pi}}\mathstrut_{#1/#2}}
\newcommand{\cpl}[1]{#1\textit{-}\mathrm{cpl}}

\DeclareMathOperator{\globalsection}{\mathbb{R}\Gamma}
\newcommand{\linsp}[1]{\langle#1\rangle}

\DeclareMathOperator{\Inf}{Inf}
\DeclareMathOperator{\deRham}{dR}
\newcommand{\infform}[4]{\mathcal{A}^{#1,\Inf}_{#3}(#4,#2)}
\newcommand{\dRform}[4]{\mathcal{A}^{#1,\deRham}_{#3}(#4,#2)}
\newcommand{\ordform}[4]{\mathcal{A}^{#1}_{#3}(#4,#2)}

\newcommand{\infformtwo}[3]{\mathcal{A}^{2,\Inf}_{#2}(#3,#1)}
\newcommand{\ordformtwo}[3]{\mathcal{A}^{2}_{#2}(#3,#1)}
\newcommand{\ndformtwo}[3]{\mathcal{A}^{2,nd}_{#2}(#3,#1)}
\newcommand{\exformtwo}[3]{\mathcal{A}^{2,ex}_{#2}(#3,#1)}

\DeclareMathOperator{\symp}{\mathcal{S}ymp}
\newcommand{\infsymp}[3]{\symp_{#2}(#3,#1)}

\newcommand{\infsympex}[3]{\symp^{ex}_{#2}(#3,#1)}

\DeclareMathOperator{\isotop}{Isot}
\newcommand{\isot}[2]{\isotop(#1,#2)}
\DeclareMathOperator{\spec}{Spec}
\DeclareMathOperator{\spf}{Spf}
\DeclareMathOperator{\lsym}{LSym}

\DeclareMathOperator{\sslashop}{\sslash}
\DeclareMathOperator{\ctr}{cotr}
\DeclareMathOperator{\stperf}{\underline{Perf}}

\DeclareMathOperator{\stvect}{\underline{Vect}}
\newcommand{\dT}{\mathbb{T}}
\newcommand{\dL}{\mathbb{L}}
\newcommand{\bF}{\mathbb{F}}
\newcommand{\bG}{\mathbb{G}}
\newcommand{\bA}{\mathbb{A}}
\newcommand{\ccA}{\mathdutchcal{A}}
\newcommand{\ccC}{\mathdutchcal{C}}

\newcommand{\cO}{\mathcal{O}}
\newcommand{\cA}{\mathcal{A}}

\newcommand{\cV}{\mathcal{V}}
\newcommand{\cE}{\mathcal{E}}
\newcommand{\cL}{\mathcal{L}}
\newcommand{\cJ}{\mathcal{J}}

\newcommand{\structuresheaf}[1]{\mathcal{O}_{#1}}

\DeclareMathOperator{\dcrit}{dCrit}
\DeclareMathOperator{\incl}{inc}
\DeclareMathOperator{\Liouv}{Liouv}
\DeclareMathOperator{\res}{Res}

\newcommand{\what}[1]{\widehat{#1}\mathstrut}
\DeclareMathOperator{\HH}{HH}
\DeclareMathOperator{\HC}{HC}
\DeclareMathOperator{\GL}{GL}

\newcommand{\formalpower}[1]{\llbracket#1\rrbracket}
\DeclareMathOperator{\ch}{Ch}

\DeclareMathOperator{\Fix}{\underline{Fix}}
\DeclareMathOperator{\Ker}{\underline{Ker}}
\DeclareMathOperator{\witt}{\mathbb{W}}
\DeclareMathOperator{\stwitt}{\underline{\mathbb{W}}}
\newcommand{\pinfty}{p^{\infty}}
\DeclareMathOperator{\atiyah}{At}
\DeclareMathOperator{\trace}{Tr}
\tikzcdset{scale cd/.style={every label/.append style={scale=#1},
		cells={nodes={scale=#1}}}}

\title{$(-1)$-Shifted Darboux theorem of derived schemes in characteristic $p>2$}
\author{Jiaqi Fu}
\begin{document}
	\maketitle
	\begin{abstract}
	The derived geometry approach to Donaldson--Thomas theory (over $\mathbb{C}$) is built on Pantev--To\"en--Vezzosi--Vaqui\'e's existence theorem of $(-1)$-shifted symplectic forms \cite{pantev2013shifted} and Brav--Bussi--Joyce's shifted Darboux theorem \cite{brav2019darboux}.	
		
	In this paper, we prove a Darboux theorem in characteristic $p>2$ for the $(-1)$-shifted symplectic forms endowed with an \textit{infinitesimal structure}. A key ingredient is Antieau's derived infinitesimal cohomology \cite{antieau2025filtrations}, which enjoys a Poincar\'e-type lemma. Our argument is in fact characteristic-free and provides a conceptual understanding of the Brav--Bussi--Joyce theorem.
	
	Moreover, we extend the existence theorem of Pantev--To\"en--Vaqui\'e--Vezzosi by constructing a de Rham $(-1)$-shifted symplectic form on $\map_{k}(X,\stperf)$, where $X$ is a Calabi--Yau $3$-fold over a field $k$ in characteristic $p>2$. We conjecture that this $(-1)$-shifted symplectic form admits an infinitesimal structure.\end{abstract}

	\section{Introduction}
	
	Donaldson--Thomas type (DT-type) invariants are virtual counts of (stable) sheaves on complex Calabi--Yau $3$-folds. Following work of Behrend \cite{behrend2009donaldson} and Pantev--To\"en--Vaqui\'e--Vezzosi \cite{pantev2013shifted}, classical DT theory can be regarded as a shadow of $(-1)$-shifted symplectic structures. Although introduced only recently, $(-1)$-shifted symplectic forms already appear implicitly in Thomas’s foundational work \cite[page 8]{thomas2000holomorphic}, in the guise of intersections of Lagrangians in a (holomorphic) symplectic manifold, each of which naturally carries a $(-1)$-shifted symplectic form \cite[\S2.2]{pantev2013shifted}. The simplest example can be given by considering the derived critical locus $\dcrit(f)$ of a function $f:U\to \mathbb{C}$, i.e. the derived pullback of $U\xrightarrow{0}T^*U\xleftarrow{df}U$, where $T^*U$ carries the tautological symplectic form. Invariants of DT-type then emanate from a shifted Darboux theorem of Brav--Bussi--Joyce stating that $(-1)$-shifted symplectic derived schemes are Zariski locally equivalent to the simplest Lagrangian intersections $\mathrm{Crit}(f)$ (see \cite{brav2019darboux}). For instance, classical DT-numbers are essentially Milnor numbers associated to the singularities of $f$.
	
	The introduction of derived foliations \cite{toen2023algebraic} by To\"en--Vezzosi helps to revisit Brav--Bussi--Joyce's result in a more conceptual way. Pantev--To\"en pointed out that, on a $(-1)$-symplectic derived scheme $(X,\omega)$, the following data are equivalent (see \cite{PT2015youtube}):
	\begin{enumerate}[label=(\arabic*)]
		\item a Lagrangian derived foliation $\mathscr{F}$ on $X$, and the exactness of $\omega$;
		\item a symplectic identification of $X$ to the derived critical locus $\dcrit(f)$ for some function.
	\end{enumerate}Hence, the shifted Darboux theorem can be interpreted as the local existence of Lagrangian derived foliations, which places derived foliation theory at centre stage.

	This work extends Brav–Bussi–Joyce's Darboux theorem to positive characteristics ($p>2$) via a foliation-theoretic approach, see Theorem \ref{main theorem}. We hope this work will initiate the investigation of Donaldson–Thomas theory in positive characteristics.

	\paragraph{Poincar\'e-type lemma.} One of the main challenge in characteristic $p>0$ is that de Rham cohomology admits too many cycles via the Cartier isomorphism, even in the derived setting. To remedy this, we instead consider a derived analogue of Grothendieck's infinitesimal cohomology, namely the \textit{Hodge-completed derived infinitesimal cohomology} $\widehat{\mathbbb{\Pi}}$ (Construction \ref{c: inf coh, ordinary}), following Antieau \cite{antieau2025filtrations}. The universality of derived infinitesimal cohomology guarantees a comparison map
	\[c:\widehat{F\mathbbb{\Pi}}(X/k)\to \widehat{F\mathrm{DR}}(X/k),\]which respects the Hodge filtrations. In contrast to the derived de Rham cohomology, derived infinitesimal cohomology satisfies a Poincar\'e-type lemma (Theorem \ref{theorem: computing inf coh of B/k}):
	\begin{theorem}
		Assume that $B$ is a simplicial commutative ring which is {\normalfont almost of finite presentation} over $k$. Then there is a natural equivalence of derived commutative rings
		\[\infcohcpl{B}{k}\simeq \lim_n B^{(n)}\simeq \lim_n \pi_0 B^{(n)},\]where $B^{(n)}$ denotes the Frobenius twist relative to $k$. In particular, $\pi_{i}\infcohcpl{B}{k}=0$ for $i\ne0,-1$.
	\end{theorem}
	\begin{remark}
		\begin{enumerate}[label=(\arabic*)]
			\item Antieau--Mathew independently discovered this Poincar\'e-type lemma in the context of quasisyntomic descent \cite{antieau2025filtrations}.
			\item The graded pieces of $c:\widehat{F\mathbbb{\Pi}}(X/k)\to \widehat{F\mathrm{DR}}(X/k)$ are given by
			\[\lsym_B(\dL_{B/k}[-1])\xrightarrow{\mathrm{Norm}} \wedge_B(\dL_{B/k})[-1].\]This difference is related to the higher $p$-power on the Lie side, see \cite[Example 1.12]{fu2024duality}.
		\end{enumerate}
	\end{remark}

	The Hodge filtration induces a fibre sequence for every $q\ge 1$,
	\[F_{[0,q-1]}\mathbbb{\Pi}(X/k)[-1]\xrightarrow{d_{\Inf}} \what{F_{q}\mathbbb{\Pi}}(X/k)\xrightarrow{\incl} \what{\mathbbb{\Pi}}(X/k).\]In this context, an \textit{infinitesimal closed $q$-form of degree $n$} is defined as a cycle
	\[\omega\in \pi_{-q-n}\what{F_q\mathbbb{\Pi}}(X/k),\]and an \textit{exact structure for $\omega$} amounts to a lift $\lambda\in \pi_{1-q-n}F_{[0,q-1]}\mathbbb{\Pi}(X/k)$ such that $d_{\Inf}\lambda=\omega$, or in other words, a null homotopy $\incl\omega\sim 0$ in $\what{\mathbbb{\Pi}}(X/k)$.

	\paragraph{Jet-Frobenius residue.}When $X=\spec B$ is affine and of finite presentation, the Poincar\'e-type lemma shows that the obstructions to exactness live in
	\[\pi_{-1}\what{F_q\mathbbb{\Pi}}_{B/k}[1]\cong\lim_n\mathstrut^1\pi_0 B^{(n)}.\]Furthermore, for each $N\ge \log_{p}(q)$, there is a \textit{jet-Frobenius power map} (Construction \ref{c: jF power map})
	\begin{equation*}
		\cJ:\prod_{n\ge N}B^{(n)}\to \what{F_q\mathbbb{\Pi}}_{B/k}[1],
	\end{equation*}where $\incl\circ \pi_0\cJ: \prod_{n\ge N}\pi_0B^{(n)}\twoheadrightarrow\lim_n^1\pi_0B^{(n)}$ agrees with the surjection in Milnor's sequence. A \textit{jet-Frobenius residue} of a $(-1)$-shifted $q$-form $\omega$ is a cycle \[\res_{JF}:=\cJ(g)\]which satisfies $\incl \res_{JF}\sim \incl \omega$ in $\infcohcpl{B}{k}$. For an arbitrarily large $N$, there always exists a jet-Frobenius residue that lifts to Hodge filtration degree $N$.
	\begin{exa}
		Consider $\infcohcpl{\mathbb{A}^1_{\bF_p}}{\bF_p}$ as $\bF_p$-module, which can be modeled by a chain complex
		\[\bF_p[x]\xrightarrow{\partial_0:f\mapsto f\otimes1-1\otimes f}\bF_p[x]\what{\otimes}_{\bF_p}\bF_p[x]\to\cdots,\]where the Hodge filtration is given by the degreewise adic filtration. Here, $\bF_p[x]\what{\otimes}_{\bF_p}\bF_p[x]$ is the $\infty$-jet of the affine line over $\bF_p$, and a jet-Frobenius residue is an infinite sum
		\[f_q(x^{p^q})\otimes1-1\otimes f_q(x^{p^q})+\cdots+f_n(x^{p^n})\otimes 1-1\otimes f_n(x^{p^n})+\cdots.\]The jet-Frobenius residue is not uniquely defined, since the degree-zero copy of $\bF_p[x]$ lies in Hodge filtration degree $0<q$.
	\end{exa}
	
	\paragraph{Infinitesimal structures and Darboux charts.} In this article, we focus on $(-1)$-shifted symplectic forms that admit a lift to derived infinitesimal cohomology,% https://q.uiver.app/#q=WzAsMyxbMCwxLCJrWy0xXSJdLFsxLDEsIlxcd2lkZWhhdHtGXzJcXG1hdGhybXtEUn19KFgvaykiXSxbMSwwLCJcXHdpZGVoYXR7Rl8yXFxtYXRoYmJie1xcUGl9fShYL2spIl0sWzAsMSwiXFxvbWVnYV97XFxkZVJoYW19IiwyXSxbMiwxLCJjIl0sWzAsMiwiXFxvbWVnYV97XFxJbmZ9IiwwLHsic3R5bGUiOnsiYm9keSI6eyJuYW1lIjoiZGFzaGVkIn19fV1d
	\[\begin{tikzcd}[ampersand replacement=\&]
		\& {\widehat{F_2\mathbbb{\Pi}}(X/k)} \\
		{k[-1]} \& {\widehat{F_2\mathrm{DR}}(X/k).}
		\arrow["c", from=1-2, to=2-2]
		\arrow["{\omega_{\Inf}}", dashed, from=2-1, to=1-2]
		\arrow["{\omega_{\deRham}}"', from=2-1, to=2-2]
	\end{tikzcd}\]There are examples from shifted cotangent stacks (Proposition \ref{prop: Liouv}), the AKSZ formalism and exact Lagrangian intersections (\S\ref{sec: 5 exa}).
	
	Meanwhile, there is a higher notion of algebraic foliations, \textit{infinitesimal derived foliations}, that has kinship to the derived infinitesimal cohomology, introduced in \cite{toen2023infinitesimal} and studied in \cite{fu2024duality}. Following Pantev--To\"en's perspective, we establish:
	\begin{maintheorem}\label{main theorem}
		Let $k$ be a field of characteristic $p>2$, and let $X$ be a derived scheme over $k$ of finite presentation. Fix an {\normalfont infinitesimal symplectic form} of degree $-1$ on $X$ relative to $k$
		\begin{equation}
			\omega\in \infsymp{-1}{k}{X}
		\end{equation} (see Definition \ref{df: infinitesimal symp}). Then, Zariski locally, $\omega$ has a Darboux chart modulo a chosen jet-Frobenius (JF) residue. More precisely, for each $x\in |X|$, there exist data $(i,j,f,\phi,\res_{JF},\eta)$:

		% https://q.uiver.app/#q=WzAsNixbMCwwLCJYIl0sWzEsMCwiVSJdLFsyLDAsIlxcbWF0aGNhbHtVfSJdLFsxLDEsIlxcbWF0aGNhbHtVfSJdLFsyLDEsIlReKlxcbWF0aGNhbHtVfSJdLFszLDAsIlxcbWF0aGJie0F9XnsxfV9rIl0sWzEsMCwiaSxcXDtcXHRleHR7XFxub3JtYWxmb250IG9wZW59IiwyLHsic3R5bGUiOnsidGFpbCI6eyJuYW1lIjoiaG9vayIsInNpZGUiOiJib3R0b20ifX19XSxbMSwyLCJqLFxcO1xcdGV4dHtcXG5vcm1hbGZvbnQgY2xvc2VkfSIsMCx7InN0eWxlIjp7InRhaWwiOnsibmFtZSI6Imhvb2siLCJzaWRlIjoidG9wIn19fV0sWzEsMywiaiIsMCx7InN0eWxlIjp7InRhaWwiOnsibmFtZSI6Imhvb2siLCJzaWRlIjoidG9wIn19fV0sWzMsNCwiMCJdLFsyLDQsImRmIl0sWzEsNCwiIiwwLHsic3R5bGUiOnsibmFtZSI6ImNvcm5lciJ9fV0sWzIsNSwiZiJdXQ==
		\[\begin{tikzcd}[ampersand replacement=\&, column sep=1.8cm]
			X \& U \& {\mathcal{U}} \& {\mathbb{A}^{1}_k} \\
			\& {\mathcal{U}} \& {T^*\mathcal{U}}
			\arrow["{i,\;\text{\normalfont open}}"', hook', from=1-2, to=1-1]
			\arrow["{j,\;\text{\normalfont closed}}", hook, from=1-2, to=1-3]
			\arrow["j", hook, from=1-2, to=2-2]
			\arrow["\lrcorner"{anchor=center, pos=0.125}, draw=none, from=1-2, to=2-3]
			\arrow["f", from=1-3, to=1-4]
			\arrow["df", from=1-3, to=2-3]
			\arrow["0", from=2-2, to=2-3]
		\end{tikzcd}\]
		\begin{enumerate}[label=(\arabic*)]
			\item $i:U=\spec B\hookrightarrow X$ an open immersion through which $x$ factors;
			\item a closed immersion $j:U\to \mathcal{U}$ into a smooth variety $\mathcal{U}$ over $k$, a function $f:\mathcal{U}\to \bA^1_k$;
			\item an equivalence $\phi:U\simeq \dcrit(\mathcal{U},f)$;
			\item a jet-Frobenius residue $\res_{JF}\in \pi_{-1}\infeqtrcohcpl{B}{k}{2}$;
			\item a homotopy $i^*\omega\stackrel{\eta}{\sim}\phi^*\omega_f+\res_{JF}$, where $\omega_f$ is the infinitesimal symplectic form of degree $-1$ induced by the derived critical locus structure of $f$ (Example \ref{exa: dcrit}).
		\end{enumerate}
		If $X$ is a derived Deligne--Mumford stack, the same statement holds except that {\normalfont Zariski} should be replaced by {\normalfont \'etale} everywhere.
	\end{maintheorem}
	
	\begin{remark}
		\begin{enumerate}[label=(\arabic*)]
			\item The choice of $\res_{JF}$ introduces ambiguity but only in higher $p$-powers. Suppose that $(i,j,f,\phi,\res_{JF}=\cJ(g),\eta)$ is a Darboux chart as above, and $\res'_{JF}=\cJ(g')$ is another jet-Frobenius residue of $i^*\omega$. If $g,g'$ are both in $\prod_{n\ge N+1}\pi_0B^{(n)}$, there is another Darboux chart $(i,j,f',\phi,\res'_{JF},\eta')$ with $f-f'\in \globalsection(\mathcal{U}^{(n)},\cO)$ (Proposition \ref{prop: perturbation of Res}). 
			\item This construction contains several choices, which could cause trouble in defining Donaldson--Thomas invariants. We plan to fix this issue in the future by studying the \textit{Darboux stacks} introduced in \cite{hennion2024gluing} in the mod $p$ context. Our theorem asserts that Darboux stacks are not empty in characteristic $p>2$.
		\end{enumerate}
	\end{remark}

	\paragraph{Sheaf moduli.} 
	Over $\mathbb{C}$, Pantev--To\"en--Vezzosi--Vaqui\'e showed that Chern character induces a $2$-shifted de Rham symplectic form on the derived moduli stack of perfect complexes. Theorem \ref{theorem: perf with symplectic structure} extends this result to characteristic $p>2$:
	\begin{theorem}
		Suppose that $R$ is a simplicial commutative ring over $\bF_p$ ($p>2$). The classifying stack $\stperf$ of perfect complexes over $R$ carries a canonical {\normalfont de Rham symplectic form of degree $2$}, which is induced by a  second Chern character $\ch^{\deRham}_2\in\what{F_2\mathrm{DR}}(\stperf/R)$.	 
	\end{theorem}
	The second Chern character $\ch^{\deRham}_2$ does not lift to derived infinitesimal cohomology in general.
	\begin{pb}
		Let $X$ be a Calabi--Yau $3$-fold over a field $k$ of characteristic $p>2$. The derived moduli stack of perfect complexes over $X$
		\[\stperf_X:=\map_{\dSt_k}(X,\stperf)\]carries a de Rham symplectic form $\int_{[X]}\ch^{\deRham}_2$ of degree $-1$ (Corollary \ref{cor: perf_X}). Does $\int_{[X]}\ch^{\deRham}_2$ admit \textit{an infinitesimal structure}?
	\end{pb}

	\paragraph{Higher Lie algebroids and Calabi--Yau $4$-folds.}The work by Brantner--Mathew and Brantner--Magidson--Nuiten generalizes Lurie--Pridham's theorem to general bases using (derived) partition Lie algebroids, see \cite{Lur,pridham2010unifying, brantner2019deformation, brantner2025formal}, which realizes Deligne and Drinfeld's principle that every formal moduli problem is controlled by a higher Lie algebra(-oid). Meanwhile, the infinitesimal derived foliations are essentially equivalent to the partition Lie algebroids under finiteness conditions \cite{fu2024duality}. In future work, we plan to construct Lagrangian foliations of $(-2)$-shifted infinitesimal symplectic stacks on the Lie side. This is of potential interest for studying the moduli $\stperf_X$ when $X$ is a Calabi--Yau $4$-fold.
	
	\paragraph{Arrangement of this paper.} We recollect the basics of derived algebraic geometry away from characteristic $0$ in \S\ref{sec: 1 DAG away from Q}, where we no longer have access to cdgas as over $\mathbb{Q}$. In \S\ref{sec: 2 Shifted symplectic structure via infinitesimal cohomology}, we introduce the notion of \textit{infinitesimal symplectic forms (of degree $n$)} over general bases. In \S\ref{sec: Formal integration of infinitesimal symplectic forms}, we develop the crucial tools including an analogue of the Poincar\'e lemma for complete derived infinitesimal cohomology and the construction of Lagrangian infinitesimal derived foliations. The main theorem is then proved in \S\ref{sec: 4 Proof} wielding the theory of square-zero extensions. Finally, \S\ref{sec: 5 exa} is devoted to the construction of examples.\vspace{-0.4em}% https://q.uiver.app/#q=WzAsNSxbMCwwLCJcXFMyIl0sWzEsMCwiXFxTMyJdLFsyLDAsIlxcUzQiXSxbMiwxLCJcXFM2Il0sWzMsMCwiXFxTNSJdLFswLDFdLFsxLDJdLFsxLDNdLFsyLDRdXQ==
	\[\begin{tikzcd}[ampersand replacement=\&,row sep=-0.1cm]
		\S2 \& \S3 \& \S4 \& \S5 \\
		\&\& \S6
		\arrow[from=1-1, to=1-2]
		\arrow[from=1-2, to=1-3]
		\arrow[from=1-2, to=2-3]
		\arrow[from=1-3, to=1-4]
	\end{tikzcd}\]
	\paragraph{Acknowledgement.} Ce projet est né d’une idée de Bertrand Toën et Marco Robalo d’utiliser ma thèse pour étendre la théorie de Donaldson–Thomas. Je les remercie chaleureusement.
	
	I am deeply grateful to Benjamin Antieau for generously sharing his insights on infinitesimal cohomology. I would like to thank Lukas Brantner and Joost Nuiten for our continuous discussions, as well as Jingbang Guo and Yupeng Wang for their helpful conversations on cohomology theories in arithmetic geometry. I am also grateful to Julian Holstein for reading an early draft.

	\section{Derived geometry away from $\mathbb{Q}$}\label{sec: 1 DAG away from Q}
	
	\subsection{Filtrations and derived rings}
	\begin{notation}
		Let $\mathbb{N}^{\ge}$ and $\mathbb{N}^{dis}$ be the $1$-categories of natural numbers, where the morphisms are given by $\ge$ and $=$ respectively. For every a stable \infcat\ $\mathcal{C}$, the \infcats\ of filtered and graded objects in $\mathcal{C}$ are defined as
		\[\fil\mathcal{C}:=\fun\big(N(\mathbb{N}^{\ge}),\mathcal{C}\big),\ \ \ \ \ \gr\mathcal{C}:=\fun\big(N(\mathbb{N}^{dis}),\mathcal{C}\big).\]The functor of taking graded pieces $\gr:\fil\mathcal{C}\to \gr\mathcal{C}$ maps every $(\ldots\to F_{2}C\to F_{1}C\to F_{0}C)$ to $\{\gr_nC\}_{n\in\mathbb{N}}$ with $\gr_nC:=\mathrm{cofib}(F_{n+1}C\to F_{n}C)$. See \cite{gwilliam2018enhancing} for a detailed discussion of filtrations in stable \infcats.
		
		For $n\in \mathbb{N}$, let $F_{\ge n}$ be the endo-functor of $\fil\mathcal{C}$ truncating $(\ldots\to F_{2}C\to F_{1}C\to F_{0}C)$ as
		\[\ldots \to F_{n+1}C\to F_{n}C\to F_{n}C\to \ldots\to F_{n}C,\] and denote $F_{[m,n]}:=\mathrm{cofib}(F_{\ge n}\to F_{\ge m})$ for all $n\ge m\ge0$. We also denote by $|-|:=F_0$ the underlying module, and abbreviate $|F_{[m,n]}-|$ as $F_{[m,n]}$ if there is no risk of confusion.
	\end{notation}

	\begin{df}
		If $\mathcal{C}$ has sequential limits, we say a filtered object $C$ is \textit{complete} if
		\[\lim_nF_nC\simeq 0.\]Let $\widehat{\fil}\mathcal{C}\subset\fil\mathcal{C}$ be the full subcategory of complete filtered objects. This inclusion has a left adjoint $(-)^\wedge$ given by the assignment $C\mapsto \lim_n F_{[0,n]}C$. 
	\end{df}We will mainly consider the case of $\mathcal{C}=\m_\mathbb{Z}$, the derived \infcat\ of $\mathbb{Z}$-modules.
	\begin{recoll}[Derived algebras]
		Let $\vect_{\mathbb{Z}}$ be the $1$-category of free abelian groups. The ordinary symmetric power \[V\mapsto \oplus_{n\ge 0} (V^{\otimes n})_{\Sigma_n}\] defines a monad on $\vect_{\mathbb{Z}}$ filtered by polynomial degrees, whose graded pieces $(V^{\otimes n})_{\Sigma_n}$ are of finite degree and preserve $\vect^{\omega}_{\mathbb{Z}}$, finite free $\mathbb{Z}$-modules. Hence it can be canonically extended to a sifted-colimit-preserving monad $\lsym$ on $\m_\mathbb{Z}$, see \cite[Theorem 2.52]{brantner2021pd} or \cite[\S4.2]{raksit2020hochschild}, called the \textit{derived symmetric power monad}.
		
		Let $\dalg:=\alg_{\lsym}(\m_{\mathbb{Z}})$ be the \infcat\ of \textit{derived commutative rings}. There is a monadic free-forgetful adjoint pair
		\[\m_{\mathbb{Z}}\rightleftarrows \dalg,\]where the forgetful functor factors through $\calg_{\mathbb{Z}}$, the \infcat\ of $\mathbb{E}_{\infty}$-algebras over $\mathbb{Z}$. The forgetful functor $\dalg\to \calg_{\mathbb{Z}}$ is small-colimit-preserving and conservative. In particular, the pushouts in $\dalg$ are computed by the tensor products over $\mathbb{Z}$, and the slice category $\dalg_{A/}$ can be identified with $\dalg_A:=\mathrm{LMod}_A(\dalg)$.
		
		Similarly, we can define the derived symmetric power monad $\lsym$ on $\fil\m_{\mathbb{Z}}$ and $\gr\m_{\mathbb{Z}}$, and
		also the \infcats\ of filtered (graded) derived commutative algebras
		\[\fil\dalg:=\alg_{\lsym}(\fil\m_\mathbb{Z}),\ \ \ \ \ \gr\dalg:=\alg_{\lsym}(\gr\m_{\mathbb{Z}}).\]
		The functors $F_0$ and $\gr$ intertwine the free-forgetful adjunctions of derived commutative algebras in these three contexts. Additionally, we may regard a derived algebra $R$ as filtered or graded derived algebra concentrated at weight $0$ by abuse of language.
	\end{recoll}
	\begin{remark}
		The full subcategory $\dalg_{\ge0}\subset\dalg$ of connective derived algebras is exactly the \infcat\ underlying the projective model structure on simplicial commutative rings, so we also write it as $\scr$. A name more \textit{styl\'e} for $\scr$ is animated rings. 
	\end{remark}
	
	\begin{notation}[Derived quotient]\label{n: derived quotient}
		Given $A\in \dalg$ and a finite collection of elements $f_i\in \pi_0 A$ ($1\le i\le r$), the derived quotient of $A$ by $f_i$ is defined as a (derived) pushout\begin{equation*}
			A\sslashop(f_1,\ldots,f_r):=A\otimes_{\mathbb{Z}[x_1,\ldots,x_r]}\mathbb{Z},
		\end{equation*}where $x_i$ is sent to $f_i$ in $A$ and $0$ in $\mathbb{Z}$. This pushout $A\sslashop(f_1,\ldots,f_r)$ is a perfect $A$-module.
	\end{notation}
	\begin{remark}
		When $A$ is an ordinary commutative ring, derived quotients are computed by Koszul complexes. It agrees with the classical quotient if $f_i$ ($1\le i\le r$) form a regular sequence.
	\end{remark}
	\begin{construction}[Cotangent complex]\label{c: non-connective cotangent complex}
		The assignment $A\mapsto \m_A$ determines a functor
		\[\m_{(-)}:\dalg\to \Catoo.\]
		Let $\mathcal{M}\to \dalg$ be its unstraightening, so that $\mathcal{M}$ is heuristically the \infcat\ of pairs $(A,M)$ where $A\in \dalg$ and $M\in \m_A$. We may identify $\mathcal{M}$ with the full subcategory $\gr_{[0,1]}\dalg\subset\gr\dalg$ spanned by graded algebras $A_*$ such that $A_i\simeq 0$ for all $i\ge 2$. Given $(A,M)$, the direct sum $A\oplus \Sigma M$ has a natural $A$-augmented algebra structure, which can be regarded as a square-zero extension. Furthermore, the assignment $(A,M)\mapsto (A\to A\oplus\Sigma M)$ has a left adjoint
		\begin{equation}
			\dL[-1]:\dalg^{\Delta^1}\to \mathcal{M}
		\end{equation}that maps $C\to B$ to $(B,\dL_{B/C}[-1])$. Here, $\dL_{B/C}$ is a non-connective version of the algebraic cotangent complex.
		
		The notion of cotangent complexes can be reproduced in every derived algebraic context, see \cite[\S4.2, \S4.4]{raksit2020hochschild}, notably the filtered or graded context. The cotangent complexes in filtered modules will be referred as \textit{internal} cotangent complexes $\dL^{int}$.
	\end{construction}

	\begin{df}\label{df: ap et afp}
	\begin{enumerate}[label=(\arabic*)]
	\item Let $R$ be a simplicial commutative ring. An $R$-module is \textit{perfect} if it is a compact object in $\m_R$. An $R$-module $M$ is \textit{almost perfect} if, for all $n\ge 0$, there exists a map $P_n\to M$ from perfect $R$-module $P_n$ such that the cofibre is $n$-connected. Let $\perf_R$ (or $\aperf_R$) be the full subcategory of $\m_R$ spanned by the (almost) perfect complexes over $R$.
	\item Given $A\in \scr_{R}$, $A$ is \textit{(almost) finitely presented} over $R$ if $\pi_0(A)$ is a finitely generated $\pi_0(R)$-algebra and $\dL_{A/R}$ is (almost) perfect over $A$.
	\end{enumerate}
	
	\end{df}	
	
	\subsection{Frobenius}
	An advantage of $\dalg_{\bF_p}$ compared with $\calg_{\bF_p}$ is that there is a genuine Frobenius. For simplicity, we only discuss the Frobenius of connective algebras. Recall that, in the \infcat\ $\scr_{\bF_p}\simeq \mathcal{P}_{\Sigma}(\mathrm{Poly}^{\omega}_{\bF_p})$, the identity functor $id$ admits an endotransformation, called \textit{the absolute Frobenius map},\[\frob:id_{\scr_{\bF_p}}\to id_{\scr_{\bF_p}}\] obtained by extending the Frobenius of (finitely generated) $\bF_p$-polynomials. More generally, fixing a base ring $R\in \scr_{\bF_p}$, there is a \textit{relative Frobenius twist} $(-)^{(1)}:\scr_{R}\to \scr_{R}$ defined by\[ B^{(1)}:=B\otimes_{R,\frob} R,\]together with a relative Frobenius map $\varphi_B:B^{(1)}\to B$.

	\begin{df}
		A simplicial commutative ring $R\in \scr_{\bF_p}$ is said to be \textit{perfect} if its Frobenius is an equivalence. The \textit{left perfection functor} is given by the assignment 
		\[(-)_{\perf}:A\mapsto \colim\big(A\xrightarrow{\frob}A\xrightarrow{\frob}A\xrightarrow{\frob} \ldots\big).\]
	\end{df}
	The left perfection is always an ordinary ring.
	\begin{lemma}\label{lemma: homotopy operation by Frob}
		For an arbitrary $A\in \scr_{\bF_p}$, the natural transformation $\frob$ induces
		\begin{itemize}
			\item the ordinary absolute Frobenius $(-)^{p}:\pi_0A\to \pi_0A$ in degree $0$;
			\item null map on $\pi_i$ for $i\ge1$.
		\end{itemize}
		In particular, the natural map $(A)_{\perf}\to (\pi_0 A)_{\perf}$ is an equivalence.
	\end{lemma}
	
	\begin{proof}
		The derived transformation $\frob$ induces a unary homotopy operation $\frob_*:\pi_i A\to \pi_iA$ of degree $0$ and weight $p$, which can be encoded by some cycle in $\pi_i\lsym^p_{\bF_p}(\bF_p[i])$ by the Yoneda lemma. However, the group $\pi_i\lsym_{\bF_p}(\bF_p[i])$ always vanishes for $i\ge 1$, \cite[Proposition 6.7]{bousfield1968operations}. This implies that $\frob_*$ is null at degree $i\ge 1$. When $i=0$, we go back to the case of classical commutative rings. Let $K$ be the fibre of $A\to \pi_0 A$. The transfer map $K\to K$ induced by $\frob$ is then null homotopic, so there is $\colim_n K\simeq 0$ and $(A)_{\perf}\simeq (\pi_0 A)_{\perf}$.
	\end{proof}

	\begin{lemma}\label{lemma: left perfection}
		The perfection functor $(-)_{\perf}$ is an \'etale hypersheaf taking values in $\mathrm{Ring}_{\bF_p}\simeq \scr^{\heartsuit}_{\bF_p}$, which is furthermore quasi-coherent on the small \'etale site.
	\end{lemma}
	\begin{proof}
		For every $R\in \scr$, let $\scr^{\acute{e}t}_R\subset\scr_R$ be the full subcategory spanned by \'etale $R\to S$.
		Suppose that $A\to B$ is an \'etale morphism in $\scr_{\bF_p}$. Consider the commuting square as follows, % https://q.uiver.app/#q=WzAsNCxbMCwwLCJBIl0sWzEsMCwiQSJdLFswLDEsIlxccGlfMEEiXSxbMSwxLCJcXHBpXzBBIl0sWzAsMSwiXFxmcm9iIl0sWzAsMl0sWzEsM10sWzIsMywiXFxmcm9iIiwyXV0=
		\[\begin{tikzcd}[ampersand replacement=\&]
			A \& A \\
			{\pi_0A} \& {\pi_0A.}
			\arrow["\frob", from=1-1, to=1-2]
			\arrow[from=1-1, to=2-1]
			\arrow[from=1-2, to=2-2]
			\arrow["\frob"', from=2-1, to=2-2]
		\end{tikzcd}\]Base change induces the following commutative diagram% https://q.uiver.app/#q=WzAsNCxbMCwwLCJcXHNjcl57XFxhY3V0ZXtlfXR9X3tBfSJdLFsxLDAsIlxcc2NyXntcXGFjdXRle2V9dH1fe0F9Il0sWzAsMSwiXFxzY3Jee1xcYWN1dGV7ZX10fV97XFxwaV8wQX0iXSxbMSwxLCJcXHNjcl57XFxhY3V0ZXtlfXR9X3tcXHBpXzBBfSJdLFswLDIsIlxcc2ltZXEiXSxbMiwzXSxbMSwzLCJcXHNpbWVxIiwyXSxbMCwxXV0=
		\[\begin{tikzcd}[ampersand replacement=\&]
		{\scr^{\acute{e}t}_{A}} \& {\scr^{\acute{e}t}_{A}} \\
		{\scr^{\acute{e}t}_{\pi_0A}} \& {\scr^{\acute{e}t}_{\pi_0A}}
		\arrow[from=1-1, to=1-2]
		\arrow["\simeq", from=1-1, to=2-1]
		\arrow["\simeq"', from=1-2, to=2-2]
		\arrow[from=2-1, to=2-2]
		\end{tikzcd}\]where the vertical equivalences are from the \'etale rigidity \cite[Theorem 7.5.1.11]{HA}. The relative Frobenius $i:\pi_0B\otimes_{\pi_0 A,\frob}\pi_0 A\to\pi_0B$ is an isomorphism by \cite[0EBS]{stack}, so $B\otimes_{A,\frob}A\to B$ is an equivalence by lifting $i$ along the right vertical equivalence. It means that $B\xrightarrow{\frob}B$ can be recovered by applying $B\otimes_A$ on $A\xrightarrow{\frob} A$. Therefore, we have $(B)_{\perf}\simeq B\otimes_A\colim(A\xrightarrow{\frob}A\xrightarrow{\frob} \ldots)\simeq B\otimes_A (A)_{\perf}$, i.e., $(-)_{\perf}$ is quasi-coherent on the small \'etale site and then an \'etale hypersheaf. Finally, the functor $(-)_{\perf}$ takes values in discrete rings by Lemma \ref{lemma: homotopy operation by Frob}.
	\end{proof}
	
	\begin{df}
		Suppose that $k\subset\bF_p$ is a field. For every $B\in \scr_k$, let $B^{(n)}:=B\otimes_{k,\frob^{\circ n}} k$ be the $n$-th Frobenius twist of $B$, and $\varphi^{n}_B:B^{(n)}\to B$ denote the $n$-th relative Frobenius. The \textit{right perfection $(B)^{\perf}_k$ of $B$ (relative to $k$)} refers to the limit algebra $\lim_n B^{(n)}$.
	\end{df}
	\begin{lemma}\label{lemma: no need for higher groups of Frob twists}
		The truncation gives rise to an equivalence in $\scr_k$ natural in $B$
		\[\lim_n B^{(n)}\simeq \lim_n \pi_0 B^{(n)}.\]
	\end{lemma}
	\begin{proof}
		Denote $F_n:=\mathrm{fib}(B^{(n)}\to \pi_0 B^{(n)})$. When $k$ is perfect, the relative Frobenius agrees with the absolute one, so it has been shown by Lemma \ref{lemma: homotopy operation by Frob} that $F_n\to F_{n-1}$ induces null maps on the level of $\pi_*$, so this map itself is also null-homotopic as $\m_{\bF_p}$ is hypercomplete. For a general field $k\supset \bF_p$, this transfer map is still null-homotopic, and it suffices to note that $\pi_*(B\otimes_{k,\frob}k)\cong\pi_*(B)\otimes_{k,\frob}k$ and that the relative Frobenius is $k$-linear.
		
		Then, consider the following fibre sequence to conclude this lemma
		\[0\simeq \lim_n F_n\to \lim_n B^{(n)}\to\lim_n \pi_0 B^{(n)}.\]\end{proof}\begin{alert}
		\begin{enumerate}[label=(\arabic*)]
			\item The existence of $\lim^1$ might prevent the right perfection from being discrete.
			\item The Frobenius of cosimplicial commutative rings induces $\mathscr{P}^0$, the Dyler--Lashof operation with no degree jump, so it is not zero in general.
		\end{enumerate}

		\end{alert}
	
	\subsection{Derived algebraic geometry}
	Derived algebraic geometry is built on simplicial commutative rings rather than ordinary commutative rings. See \cite{toen2008homotopical} for a general study. 
	\begin{recoll}
		Let $\daff:=(\scr)^{op}$ be the \infcat\ of affine derived schemes, where the \'etale morphisms of algebras induce a Grothendieck topology. Given a cocomplete \infcat\ $\mathcal{C}$, a $\mathcal{C}$-valued \'etale (hyper)sheaf $F$ is a functor $F:\daff^{op}\to \mathcal{C}$ such that, for every \'etale (hyper)covering $U^\bullet\to X$, the augmented cosimplicial diagram
		\[F(X)\to F(U^\bullet)\]exhibits $F(X)\simeq \tot F(U^\bullet)$ as the totalization of $ F(U^\bullet)$. If $\mathcal{C}$ is presentable, the inclusion of sheaves into presheaves has a localizing left adjoint called sheafification,
		\[\mathrm{PSH}(\daff,\mathcal{C})\to \mathrm{SH}^{\acute{e}t}(\daff,\mathcal{C}).\]
		The \infcat\ of \textit{derived stacks} is by definition $\dSt:=\mathrm{SH}^{\acute{e}t}(\daff,\mathcal{S})$ of \'etale sheaves in spaces. The derived stacks form an $\infty$-topos in the sense of \cite[\S6]{HTT}.
	\end{recoll}
	\begin{exa}
	\begin{enumerate}[label=(\arabic*)]
		\item There is a functor $\spec: \dalg^{op}\to \dSt$ that maps $R$ to $\map_{\dalg}(R,-)$. The restriction $\spec|_{\scr^{op}}$ agrees with the Yoneda functor of $\daff$. 
		\item Derived schemes, derived Deligne--Mumford stacks and derived Artin stacks are natural generalizations of their classical counterparts. See \cite[\S1.3, \S2.2 ]{toen2008homotopical} and \cite[\S20.6]{SAG}.
		\item There is a $\Catoo$-valued hypersheaf
		\[R\mapsto \perf_R.\]Let $\perf^{\simeq}_R\subset\perf_R$ be the maximal Kan subcomplex. We obtain a hypercomplete derived stack $\stperf$ such that $\stperf(R):=\perf^{\simeq}_R$. This is in fact a locally derived Artin stack following \cite[Theorem 1.1]{toen2007moduli}.
		\item Given $R\in \scr$, the cotangent complex $\dL_{-/R}$ defines a $\m_R$-valued \'etale hypersheaf on $\dSt_R:=(\dSt)_{/\spec R}$.
	\end{enumerate}
	\end{exa}	
	\begin{df}
		For a derived stack $X\in\dSt$, the derived \infcat\ of \textit{quasi-coherent sheaves over $X$} is defined by right Kan extension, i.e.,
		\[\QC(X):=\lim_{\spec R\to X} \m_R.\]
	\end{df}
	\begin{df}
		Given $R\in \scr$, we say $F\in \dSt_R$ has a \textit{(global) cotangent complex} if there exists an $\dL_{F/R}\in \QC(F)$ such that, for every point $x:\spec A\to F$ and $M\in \m_A$, there is a natural equivalence 
		\[\map_A(x^*\dL_{F/R},M)\simeq F(A\oplus M)\times_{F(A)}\{x\}.\]
	\end{df}
	\begin{theorem}\cite[Theorem 1.4.3.2]{toen2008homotopical}
		Let $F$ be a derived Artin stack over $R$. Then, $F$ admits a global cotangent complex $\dL_{F/R}$.
	\end{theorem}
	\begin{df}
		Suppose that $f:X\to Y$ is a morphism in $\dSt_R$ with a global cotangent complex, where $R\in \scr$. We say a morphism $f:X\to Y$ is \textit{formally \'etale, smooth or unramified} if the relative cotangent complex $\dL_{f}$ is contractible, vector bundle or connected, respectively\footnote{Caveat: Our notion of being formally unramified is \textbf{weaker} than the standard notion.}.
	\end{df}
	\subsection{Formal derived schemes}
	We recollect the basics of formal derived algebraic geometry here and prove some useful propositions for later use. The reader may see \cite[\S6]{lurie2004derived} and \cite[\S7, \S8]{SAG} \footnote{SAG is in the scope of spectral algebraic geometry, but most of the properties of formal geometry are shared with derived algebraic geometry. A recent note \cite{chough2024formal} explains the comparison between them.} for a detailed study. Thanks to the Frobenius map, we may have some shortcuts since we work over some $\bF_p$.
	
	\begin{df}\label{df: adic ring}
		Let $\mathrm{Ring}^{\adic}_{\bF_p}$ be the $1$-category of topological rings over $\bF_p$ whose topology is the adic topology of a finitely generated ideal of definition. Define the \infcat\ of adic simplicial commutative rings over $\bF_p$ as the pullback
		\begin{equation}
			\scr^{\adic}_{\bF_p}:=\scr\times_{\mathrm{Ring}_{\bF_p}}\mathrm{Ring}^{\adic}_{\bF_p},
		\end{equation}where the left leg is $\pi_0$. Given $A$ in $\scr^{\adic}_{\bF_p}$ and a chosen ideal of definition $(f_i)_{1\le i\le r}\subset\pi_0 A$, there is a pro-object $\big(A\sslashop (f^{p^n}_1,\ldots,f^{p^n}_r)\big)_{n\ge0}$ corepresenting a functor
		\begin{equation}
			\spf(A):=\underset{n\in\mathbb{N}}{\colim} \spec A\sslashop (f^{p^n}_1,\ldots,f^{p^n}_r) 		\end{equation} 
			called the \textit{formal spectrum of $A$.}

	\end{df}
	
	\begin{lemma}
		For every $A\in \scr^{\adic}_{\bF_p}$ and a choice of $(f_i)_{1\le i\le r}$, the induced functor $\spf(A)$ admits a natural transformation to $\spec(A)$. This natural transformation induces $(-1)$-truncated maps
		\[\spf(A)(B)\hookrightarrow\spec(A)(B),\]which identifies the left-hand side with the components of $\map_{\scr_{\bF_p}}(A,B)$ spanned by morphisms $F:A\to B$ such that the image of every $f_i$ is nilpotent in $\pi_0(B)$.
		
		Furthermore, the assignment $A\mapsto \spf A$ as above is independent from the choice of $(f_i)$. Hence, this gives rise to a functor
		\begin{equation*}
			\spf: 	\scr^{\adic,op}_{\bF_p}\to \mathrm{PSH}(\daff,\mathcal{S}).
		\end{equation*}
	\end{lemma}
	\begin{proof}
		
		The space $\big(\spec A\sslashop (f^{p^n}_1,\ldots,f^{p^n}_r)\big)(B)$ is given by the pullback $F_n$% https://q.uiver.app/#q=WzAsNCxbMSwxLCJcXE9tZWdhXntcXGluZnR5fUJee1xcdGltZXMgcn0iXSxbMCwxLCJcXG1hcF97XFxzY3Jfe1xcYkZfcH19KEEsQikiXSxbMSwwLCIqIl0sWzAsMCwiRl9uIl0sWzEsMF0sWzMsMV0sWzMsMl0sWzIsMCwiMCJdXQ==
		\[\begin{tikzcd}[ampersand replacement=\&]
			{F_n} \& {*} \\
			{\map_{\scr_{\bF_p}}(A,B)} \& {\Omega^{\infty}B^{\times r}.}
			\arrow[from=1-1, to=1-2]
			\arrow[from=1-1, to=2-1]
			\arrow["0", from=1-2, to=2-2]
			\arrow[from=2-1, to=2-2]
		\end{tikzcd}\]
		The transfer map $F_n\to F_{n+1}$ is induced by post-composing the absolute Frobenius of $B$ at the bottom right corner. Then, the space $\spf (A)\simeq \colim_n F_n$ is equivalent to the pullback 
		
		% https://q.uiver.app/#q=WzAsNCxbMSwxLCJcXE9tZWdhXntcXGluZnR5fShCKV97XFxwZXJmfV57XFx0aW1lcyByfSJdLFswLDEsIlxcbWFwX3tcXHNjcl97XFxiRl9wfX0oQSxCKSJdLFsxLDAsIioiXSxbMCwwLCJGX3tcXGluZnR5fSJdLFsxLDBdLFsyLDAsIjAiXSxbMywxXSxbMywyXV0=
		\begin{equation}\label{e: reformulate SPF}
			\begin{tikzcd}[ampersand replacement=\&]
				{F_{\infty}} \& {*} \\
				{\map_{\scr_{\bF_p}}(A,B)} \& {\Omega^{\infty}(B)_{\perf}^{\times r}.}
				\arrow[from=1-1, to=1-2]
				\arrow[from=1-1, to=2-1]
				\arrow["0", from=1-2, to=2-2]
				\arrow[from=2-1, to=2-2]
			\end{tikzcd}
		\end{equation}
		Thus this lemma follows from Lemma \ref{lemma: left perfection}.\end{proof}
	\begin{remark}\label{rk: adic ideal system}
		The definition of $\spf$ here is slightly different from \cite[\S6.1]{lurie2004derived}, since here the tower of algebras is given by modding out $f^{p^n}_i$ instead of $f^{2^n}_i$. This lemma shows that it causes no change. Similarly, let $S$ be a set of cardinal $r$, and $S\to A$ be a map sending $s_i$ to $f_i$. The pro-system of algebras
		\begin{equation}
			A\sslash(f^I||I|=n)\simeq A\otimes_{\mathbb{Z}[S]}Z[S]/J^{n}
		\end{equation}also pro-corepresents $\spf A$, where $J=(s_1,\ldots,s_r)\subset\mathbb{Z}[S]$ is the ideal generated by $S$.
	\end{remark}	
	\begin{prop}
		The functor $\spf(A)$ is an \'etale hypersheaf.
	\end{prop}
	\begin{proof}
		In diagram (\ref{e: reformulate SPF}), the three corners on the bottom or on the right are clearly \'etale hypersheaves in $B$. Therefore, so is the forth corner $\spf(A)(B)\simeq F_\infty$.
	\end{proof}
	\begin{df}[Derived adic completeness]
		Suppose that $A\in \scr^{\adic}_{\bF_p}$ and $I\subset\pi_0 A$ is an ideal of definition. A module $M\in\m_A$ is called $I$-adic complete if, for all $f\in I$, the mapping space
		\[\map_{A}(A[f^{-1}],M)\]is contractible. Let $\m^{\cpl{I}}_A\subset\m_A$ be the full subcategory of $I$-adic complete modules. The adic ring $A$ is said to be complete if $A\in \m^{\cpl{I}}_A\subset\m_A$ for some ideal of definition $I$.
	\end{df}

	\begin{prop}
		\begin{enumerate}[label=(\arabic*)]
			\item \cite[Variant 7.3.5.6]{SAG} The inclusion $\m^{\cpl{I}}_A\subset\m_A$ has a localizing left adjoint $(-)^\wedge_I$, which is compatible with the symmetric monoidal structure of $\m_A$. Therefore, there is a natural tensor product $\widehat{\otimes}$ on $\m^{\cpl{I}}_A$ so that $(-)^{\wedge}_I$ is symmetric monoidal. In particular, $A^\wedge_I$ is a complete adic algebra over $A$.
			\item \cite[Proposition 6.1.4]{lurie2004derived} Let $(f_i)_{1\le i\le r}$ be a chosen set of generators of $I$, and $A_n:=A\sslashop (f^{p^n}_1,\ldots,f^{p^n}_r)$. The pullback along $i:\spf A\to \spec A$ 
			\begin{equation}
				i^*:\QC_{\spec A}\to \QC_{\spf A}
			\end{equation}can be identified with $(-)^{\wedge}_I$ as symmetric monoidal functors.
		\end{enumerate}
	\end{prop}
	\begin{prop}\label{prop: Spf cotangent}
		Let $R$ be a simplicial commutative ring over $\bF_p$. Consider an adic ring $A$ over $R$, i.e. $(A,I)\in \scr^{\adic}_{R}\simeq \scr_R\times_{\scr_{\bF_p}}\scr^{\adic}_{\bF_p}$.
		As a derived stack over $R$, $\spf A$ has a global cotangent complex $\dL_{\spf A/R}$. Moreover, there are natural equivalences,
		\begin{equation}
			(\dL_{A/R})^\wedge_I\simeq \lim_n\dL_{A_n/R}\simeq \dL_{\spf A/R}.
		\end{equation}
	\end{prop}
	\begin{proof}
	For every map $x:\spec B \to \spf A$ over $R$ and every $B$-module $M$, the pullback of spaces
	\begin{equation}\label{e: cpl cotangent general case}
	\map_{\dSt_R}(\spec B\oplus M, \spf A)\times_{\map_{\dSt_R}(\spec B, \spf A)}\{x\}
	\end{equation}is computed as the colimit of the following spaces
	\begin{equation*}
	\map_{\dSt_R}(\spec B\oplus M, \spec A_n)\times_{\map_{\dSt_R}(\spec B, \spec A_n)}\{x\},
	\end{equation*}which is well-defined for large $n$ and is equivalent to $\map_{B}(x^*\dL_{A_n/R},M)$. Next, recall that there is a fibre sequence
	\begin{equation*}
		A_n\otimes_A\dL_{A/R}\to \dL_{A_n/R}\to \dL_{A_n/A},
	\end{equation*}where $\dL_{A_n/A}\simeq A\otimes_{R[x^{p^n}_1,\ldots,x^{p^n}_r]}\dL_{R/R[x^{p^n}_1,\ldots,x^{p^n}_r]}$. One can see that the transfer map between $\dL_{A_n/A}$ is null-homotopic and hence we have $(\dL_{A/R})^\wedge_I\simeq \lim_n\dL_{A_n/R}$. Meanwhile, the space in (\ref{e: cpl cotangent general case}) is equivalent to $\map_{B}(B\otimes_{A}\dL_{A/R},M)$ as there exists some factorization $A\to A_n \to B$. Hence, $\dL_{\spf A/R}$ exists and is the quasi-coherent sheaf corresponding to $(\dL_{A/R})^\wedge_I\simeq \lim_n\dL_{A_n/R}$.
	\end{proof}

	\section{Shifted symplectic structures via infinitesimal cohomology}\label{sec: 2 Shifted symplectic structure via infinitesimal cohomology}
	\subsection{Infinitesimal cohomology and foliations}

	\begin{construction}\label{c: inf coh, ordinary}
		There is an adjoint pair
		\begin{equation}
			F\mathbbb{\Pi}:\dalg^{\Delta^1}\rightleftarrows \fil\dalg:(F_0\to \gr_0),
		\end{equation}where the left adjoint $F\mathbbb{\Pi}$ maps $R\to B$ to its \textit{Hodge filtered derived infinitesimal cohomology} $\infcohfil{B}{R}$. See \cite{antieau2025filtrations} for details. We will mainly consider its derived completion
		\begin{equation}
			\infcohcplfil{B}{R}:=\lim_{n}F_{[0,n]}\mathbbb{\Pi}\mathstrut_{B/R}.
		\end{equation}
		In particular, we write $\infcohcpl{B}{R}:=|\infcohcplfil{B}{R}|$.
	\end{construction}
	\begin{lemma}\label{lemma: graded pieces of infcoh}
		For every $R\to B\in \dalg$, the graded pieces of $F\mathbbb{\Pi}\mathstrut_{B/R}$ form a freely generated graded derived algebra, more precisely,
		\[\gr\mathbbb{\Pi}\mathstrut_{B/R}\simeq \lsym_B\dL_{B/R}[-1].\]
	\end{lemma}
	\begin{proof}
		The functor $\gr\mathbbb{\Pi}(-/-)$ has a right adjoint
		\[G:\gr\dalg\xrightarrow{g}\fil\dalg\xrightarrow{(F_0\to\gr_0)}\dalg^{\Delta^1},\]where $g$ equips a graded algebra $A_*$ with trivial transfer maps. This right adjoint sends $A_*$ to $A_0\to A_0\oplus A_1[1]$ where $ A_0\oplus A_1[1]$ is a square-zero extension. Thus, a morphism $f:(R\to B)\to G(A_*)$ is determined by $B\to A_0$ together with a $B$-module map $\dL_{B/R}\to A_{1}[1]$. Equivalently, this amounts to a graded algebra map
		\[\lsym_B \dL_{B/R}[-1]\to A_*.\]
	\end{proof}

	\begin{prop}\label{prop: InfCoh is hypersheaf}
		Fixing $R\in \scr$, the functors $\widehat{F\mathbbb{\Pi}}\mathstrut_{-/R}$ and $F_{[0,n]}\mathbbb{\mathbbb{\Pi}}\mathstrut_{-/R}$ ($n\ge0$) are all \'etale hypersheaves over $\daff_R$ that take values in completely filtered derived algebras over $R$.
	\end{prop}
	\begin{proof}
		Since the structure sheaf $\cO$ and $\lsym^n_\cO\dL_{\cO/R}$ are \'etale hypersheaves, $\gr\mathbbb{\Pi}\mathstrut_{B/R}$ and $\gr_{[0,n]}\mathbbb{\Pi}\mathstrut_{B/R}$ are \'etale hypersheaves in graded algebras by Lemma \ref{lemma: graded pieces of infcoh}. Then, this proposition holds since $\widehat{F\mathbbb{\Pi}}\mathstrut_{-/R}$ and $F_{[0,n]}\mathbbb{\mathbbb{\Pi}}\mathstrut_{-/R}$ are complete.\end{proof}
	\begin{df}
		Given $X\in\dSt_R$, the completed (Hodge-filtered) derived infinitesimal cohomology of $X$ relative to $R$ is defined by Kan extension,
		\[\infcohcpl{X}{R}:=\lim\limits_{\spec(B)\to X}\infcohcpl{B}{R},\quad\quad\infcohcplfil{X}{R}:=\lim\limits_{\spec(B)\to X}\infcohcplfil{B}{R}.\]
	\end{df}
	\begin{exa}
		Let $A$ be an adic simplicial commutative ring over $R$, and $(A_n)_{n\ge0}$ be a pro-system of $R$-algebras corepresenting it as in Definition \ref{df: adic ring}. Then, we have the natural equivalence
		\begin{equation*}
			\infcohcplfil{\spf A}{R}\simeq \lim_{n\ge 0}\infcohcplfil{A_n}{R}.
		\end{equation*}
	\end{exa}
	
	We will also consider an internal notion of infinitesimal cohomology in the filtered algebra context, mainly for computational good. The reader may regard it as an abstract version of deformation to the normal cone of derived foliations, see \cite[\S1.3]{toen2015caracteres}.
	\begin{variant}\label{var: internal infinitesimal cohomology}
		Observe that in Construction \ref{c: inf coh, ordinary} the crucial input is the derived algebra context in the sense of \cite[\S4.2]{raksit2020hochschild}. If we replace the ordinary derived algebra context of $\m_{\mathbb{Z}}$ with that of $\fil\m_\mathbb{Z}$, we obtain a theory of \textit{internal Hodge-filtered derived infinitesimal cohomology} of filtered derived rings. More precisely, for every map $\mathcal{R}\to \mathcal{B}$ in $\fil\dalg$, it functorially assigns to this map a bi-filtered algebra over $\mathcal{R}$
		\begin{equation}
				(\mathcal{R}\to)F\mathbbb{\Pi}^{int}_{\mathcal{B}/\mathcal{R}},
		\end{equation}whose internal graded pieces form a free internal graded algebra
		\[\lsym_{\mathcal{B}}\dL^{int}_{\mathcal{B/\mathcal{R}}}[-1].\]
		We mainly consider its completion with respect to both filtrations $\what{F\mathbbb{\Pi}}^{int}_{\mathcal{B}/\mathcal{R}}$.
	\end{variant}
	We also consider a notion of higher algebraic foliations, the infinitesimal derived foliations\cite{toen2023infinitesimal}, in the disguise of \textit{foliation-like algebras} \cite{fu2024duality}. The foliation-like algebras have $\infcohcplfil{B}{R}$ as primary examples and help compute $\infcohcplfil{B}{R}$ in \S\ref{sec: 3.1 local exactness}. Assume that there is $R\to B$ in $\scr$ such that $B$ is coherent in the sense of \cite[Definition 7.2.4.16]{HA}. Foliation-like algebras live in the ambient \infcat\ defined by a cartesian square,
	% https://q.uiver.app/#q=WzAsNCxbMCwxLCJcXHtCXFx9Il0sWzAsMCwiXFxtYXRoY2Fse0R9XntcXGZpbH1fe0IvUn0iXSxbMSwxLCJcXGRhbGdfe1J9Il0sWzEsMCwiXFxmaWxcXGRhbGdfUiJdLFsxLDBdLFswLDJdLFszLDIsIlxcZ3IiLDJdLFsxLDNdLFsxLDIsIiIsMCx7InN0eWxlIjp7Im5hbWUiOiJjb3JuZXIifX1dXQ==
	\[\begin{tikzcd}[ampersand replacement=\&]
		{\mathcal{D}^{\fil}_{B/R}} \& {\fil\dalg_R} \\
		{\{B\}} \& {\dalg_{R}.}
		\arrow[from=1-1, to=1-2]
		\arrow[from=1-1, to=2-1]
		\arrow["\lrcorner"{anchor=center, pos=0.125}, draw=none, from=1-1, to=2-2]
		\arrow["\gr"', from=1-2, to=2-2]
		\arrow[from=2-1, to=2-2]
	\end{tikzcd}\]
	A \textit{foliation-like algebra} over $B$ (relative to $R$) is then a complete $\mathcal{A}\in\mathcal{D}^{\fil}_{B/R}$ such that $\gr \mathcal{A}$ is naturally equivalent to $\lsym_B\gr_1\mathcal{A}$ as a graded algebra over $B$.
	\begin{remark}
	The \textit{completed Hodge-filtered derived de Rham cohomology} $\what{DR}(B/R)$ also lies in $\mathcal{D}^{\fil}_{B/R}$, though it is not a foliation-like algebra. The filtered algebra $\infcohfil{B}{R}$ is initial in $\mathcal{D}^{\fil}_{B/R}$. Thus, there is a canonical comparison map
	\begin{equation}\label{e: Inf VS DeRham}
		c:\infcohcplfil{B}{R}\to \widehat{F\mathrm{DR}}\mathstrut_{B/R},
	\end{equation}
	whose graded pieces \[\gr c:\lsym_B(\dL_{B/R}[-1])\to \wedge^*_B(\dL_{B/R})[-*]\simeq \Gamma_B(\dL_{B/R}[-1])\]are computed as the norm map of permutation groups.
	\end{remark}
	
	\begin{prop}\label{prop: inf foliation basics}Let $R\to B$ be a morphism in $\scr$ and suppose that $B$ is coherent. Then:
		\begin{enumerate}[label=(\arabic*)]
			\item Let $\scr^{\wedge\mathrm{aft}}_{B/R}\subset\scr_{R//B}$ be the full subcategory of $R\to A\to B$ such that $\pi_0(A)\to \pi_0(B)$ is surjective with kernel $I$, $\pi_0(A)$ is $I$-adically complete, and $\dL_{B/A}$ is almost perfect over $B$. Then, $\infcohcplfil{-}{R}$ embeds $\scr^{\wedge\mathrm{aft}}_{R//B}$ fully faithfully into $\mathcal{D}^{\fil}_{B/R}$. Moreover, the essential image is spanned by the foliation-like algebras $\mathcal{A}$ for which $\gr_1\mathcal{A}$ is a connective almost perfect $B$-module.
			\item Let $\mathcal{A}^{\bullet}$ be a cosimplicial diagram of foliation-like algebras over $B$, such that $\gr\mathcal{A}^{\bullet}$ is a diagram in $\aperf_B$ and its totalization $L$ is also almost perfect over $B$. Then, the straightforward totalization $\tot\mathcal{A}^{\bullet}$ is naturally a foliation-like algebra and the $\gr_1$-piece is $L$.
		\end{enumerate}
	\end{prop}
	\begin{proof}
		The first statement is exactly \cite[Proposition 6.17]{Fu2025} The second statement is proved by the contravariant categorical equivalence in \cite[Theorem 4.25]{fu2024duality} and the fact that the forgetful functor of partition Lie algebroids reflexes geometric realizations.
	\end{proof}
	From this point of view, the completed filtered derived infinitesimal cohomology $\infcohcplfil{B}{A}$ is an ``integrable infinitesimal derived foliation''.
	
	\begin{exa}\label{rk: adic vs Infcoh}
		Let $A$ be a simplicial commutative ring. Choose a set of $r$ cycles $(s_1,\ldots,s_r):S\to A$ ($r=|S|$) and denote $B:=A\sslashop(s_1,\ldots,s_r)$. Then, there are equivalences
		\begin{equation*}
			\infcohcplfil{B}{A}\simeq \big(A\otimes_{\mathbb{Z}[S]}\infcohfil{\mathbb{Z}}{\mathbb{Z}[S]}\big)\mathstrut^\wedge\simeq \big(A\otimes_{\mathbb{Z}[S]}\mathbb{Z}\formalpower{S}_{\adic}\big)\mathstrut^\wedge,
		\end{equation*}whose $F_{[0,n]}$-component is $A_n:=A\sslashop(s^I||I|=n)$.
		
		Conversely, suppose that $\mathcal{A}$ is a foliation-like algebra over $B$ (relative to $R$) such that $\gr_1\mathcal{A}$ is a finite free-$B$ module (without shifting). Set $A:=|\cA|$. Then, $A\to B$ can be written as a derived quotient: Choose $B$-module generators $S\to \gr_1\mathcal{A}$ and lift these cycles to $F_1\cA\to A$. Note that 
		$\mathcal{A}\widehat{\otimes}_{\mathbb{Z}[S]}\mathbb{Z}\simeq B$ in $\widehat{\fil}\dalg$. Since $\mathbb{Z}$ is a perfect $\mathbb{Z}[S]$-module, $\mathcal{A}\otimes_{\mathbb{Z}[S]}\mathbb{Z}\simeq B$ is already complete. Then, we have $A\sslashop(s_1,\ldots,s_r)\simeq B$. Additionally, $A$ is a complete adic ring, and the pro-system $(F_{[0,n]}\mathcal{A})_{n\ge0}$ pro-corepresents $\spf A$.
	\end{exa}
	
	\subsection{Infinitesimal symplectic forms}
	In ordinary geometry, a symplectic form is a non-degenerate closed $2$-form. The subtlety of defining a symplectic form in derived geometry over general bases is around concerns what a closed form is.
	\begin{df}
		Let $R\to B$ be a morphism of simplicial commutative rings. The space of \textit{infinitesimal (resp. de Rham) closed $q$-forms of degree $n$} on $\spec B$ relative to $R$ is, by definition, the space of $R$-linear maps of the form
		\[\omega:R\to \infeqtrcohcpl{B}{R}{q}[q+n],\quad (\textit{resp. } \omega:R\to \what{F_q\mathrm{DR}}_{B/R}[q+n]).\]Alternatively, the spaces of closed forms are defined as the following $\infty$-loop spaces
		\begin{align*}
			\infform{q}{n}{R}{\spec B}&:=\Omega^{\infty}\infeqtrcohcpl{B}{R}{q}[q+n]\text{, or }\\
			\dRform{q}{n}{R}{\spec B}&:=\Omega^\infty\what{F_q\mathrm{DR}}_{B/R}[q+n]\text{ respectively}.
		\end{align*}
	\end{df}
	\begin{remark}
	Pantev--To\"en--Vaqui\'e--Vezzosi's foundational work \cite{pantev2013shifted} uses de Rham shifted closed forms. The morphism $c:\infcohcplfil{-}{R}\to \widehat{F\mathrm{DR}}\mathstrut_{-/R}$ (\ref{e: Inf VS DeRham}) induces a comparison map\begin{equation}
		|c|:\infform{q}{n}{R}{-}\to \dRform{q}{n}{R}{-}.
	\end{equation}
	Though $c$ is an equivalence in characteristic $0$, its two ends can differ dramatically even in the smooth case. In fact, the ordinary de Rham cohomology carries too many cycles in the mod $p$ case via the Cartier isomorphism, while Grothendieck's infinitesimal cohomology is rather restrictive following Ogus' work \cite{ogus1975cohomology}.
	\end{remark}
	\begin{construction}\label{c: infinitesimal SYMP}
		Define $\ordform{q}{n}{R}{\spec B}:=\Omega^{\infty+n}(\wedge^q_B\dL_{B/R})$, the \textit{space of (ordinary) $q$-forms of degree $n$ over $\spec B$ relative to $R$}. A $2$-form $\omega_0\in \ordform{q}{n}{R}{\spec B}$ of degree $n$ is adjoint to the $B$-linear map\[\Theta_{\omega_0}:\dT_{B/R}\xrightarrow{\simeq}\dL_{B/R}[n].\]We say that $\omega_0$ is \textit{non-degenerate} if $\Theta_{\omega_0}$ is an equivalence. Meanwhile, the natural map $F_qc$ induces
		\[\what{F_q\mathbbb{\Pi}}(-/R)\to \what{F_q\mathrm{DR}}(-/R)\to \wedge^q(\dL_{-/R})[-q].\]
		Therefore, every $\omega\in \infform{q}{n}{R}{\spec B}$ has an underlying $q$-form of degree $n$. In particular, every $\omega\in \infformtwo{n}{R}{\spec B}$ induces a $B$-module morphism $\Theta_\omega:\dT_{B/R}\to \dL_{B/R}[n]$.
		
		 An \textit{infinitesimal symplectic form} of degree $n$ over $\spec B$ (relative to $R$) is defined as a non-degenerate infinitesimal closed $2$-form of degree $n$. The space of infinitesimal symplectic forms of degree $n$ over $\spec B$ is defined by the cartesian diagram as follows,
		% https://q.uiver.app/#q=WzAsNCxbMSwwLCJcXG5kZm9ybXR3b3tufXtSfXstfSJdLFsxLDEsIlxcb3JkZm9ybXR3b3tufXtSfXstfSJdLFswLDEsIlxcaW5mZm9ybXR3b3tufXtSfXstfSJdLFswLDAsIlxcaW5mc3ltcHtufXtSfXstfSJdLFswLDEsIiIsMCx7InN0eWxlIjp7InRhaWwiOnsibmFtZSI6Imhvb2siLCJzaWRlIjoidG9wIn19fV0sWzIsMSwifGN8IiwyXSxbMywyXSxbMywwXSxbMywxLCIiLDEseyJzdHlsZSI6eyJuYW1lIjoiY29ybmVyIn19XV0=
		\[\begin{tikzcd}[ampersand replacement=\&]
			{\infsymp{n}{R}{-}} \& {\ndformtwo{n}{R}{-}} \\
			{\infformtwo{n}{R}{-}} \& {\ordformtwo{n}{R}{-}}
			\arrow[from=1-1, to=1-2]
			\arrow[from=1-1, to=2-1]
			\arrow["\lrcorner"{anchor=center, pos=0.125}, draw=none, from=1-1, to=2-2]
			\arrow[hook, from=1-2, to=2-2]
			\arrow["{|c|}"', from=2-1, to=2-2]
		\end{tikzcd}\]where $\ndformtwo{n}{R}{-}\subset \ordformtwo{n}{R}{-}$ is the full sub-presheaf of non-degenerate $2$-forms.
	\end{construction}
	\begin{remark}
		\begin{enumerate}[label=(\arabic*)]
			\item The four vertices in this square are all \'etale hypersheaves over $\daff_R:=(\scr_R)^{op}$ following Proposition \ref{prop: InfCoh is hypersheaf}.
			\item We can also consider $\symp^{\deRham}_R(-,n)$ the hypersheaf of \textit{de Rham symplectic forms} of degree $n$ as in \cite{pantev2013shifted}. The natural map $c$ induces $|c|:\infsymp{n}{R}{-}\to\symp^{\deRham}_R(-,n)$.
		\end{enumerate}

	\end{remark}
	Next, we extend this construction from affine case to global case by right Kan extension:
	\begin{df}\label{df: infinitesimal symp}
		Let $X$ be a derived stack over $R$. The space of \textit{infinitesimal symplectic forms of degree $n$} on $X$ (relative to $R$) is defined by
		\[\infsymp{n}{R}{X}:=\lim_{\spec(B)\to X}\infsymp{n}{R}{\spec B}.\]Defining $\ordformtwo{n}{R}{X}$, $\infformtwo{n}{R}{X}$, and $\ndformtwo{n}{R}{X}$ similarly, there remains a pullback
		% https://q.uiver.app/#q=WzAsNCxbMSwwLCJcXG5kZm9ybXR3b3tufXtSfXtYfSJdLFsxLDEsIlxcb3JkZm9ybXR3b3tufXtSfXtYfSJdLFswLDEsIlxcaW5mZm9ybXR3b3tufXtSfXtYfSJdLFswLDAsIlxcaW5mc3ltcHtufXtSfXtYfSJdLFswLDEsIiIsMCx7InN0eWxlIjp7InRhaWwiOnsibmFtZSI6Imhvb2siLCJzaWRlIjoidG9wIn19fV0sWzIsMSwifGN8IiwyXSxbMywyXSxbMywwXSxbMywxLCIiLDEseyJzdHlsZSI6eyJuYW1lIjoiY29ybmVyIn19XV0=
		\[\begin{tikzcd}[ampersand replacement=\&]
			{\infsymp{n}{R}{X}} \& {\ndformtwo{n}{R}{X}} \\
			{\infformtwo{n}{R}{X}} \& {\ordformtwo{n}{R}{X}.}
			\arrow[from=1-1, to=1-2]
			\arrow[from=1-1, to=2-1]
			\arrow["\lrcorner"{anchor=center, pos=0.125}, draw=none, from=1-1, to=2-2]
			\arrow[hook, from=1-2, to=2-2]
			\arrow["{|c|}"', from=2-1, to=2-2]
		\end{tikzcd}\]
	\end{df}
	As commented in \cite{pantev2013shifted}, the space of shifted symplectic forms is generally complicated, but the space of ordinary forms over derived Artin stacks can be made explicit. Since the proof of \cite[Proposition 1.14]{pantev2013shifted} is characteristic-free, we have:
	\begin{prop}\label{prop: global section of forms}
		Let $F\in \dArt_R$ be a derived Artin stack over $R$, in the sense of \cite[\S1.3]{toen2008homotopical}. Then, there is a natural equivalence
		\[\map_{\QC(F)}\big(\structuresheaf{F},(\wedge^{q}_{\structuresheaf{F}}\dL_{F/R}[n])\big)\xrightarrow{\simeq}\ordform{q}{n}{R}{F}.\] In particular, every $\omega\in \infformtwo{n}{R}{F}$ induces an $\cO_F$-linear map
		\[\Theta_\omega:\dT_{F/R}\to \dL_{F/R}[n],\]and $\omega$ is an infinitesimal symplectic form if and only if $\Theta_\omega$ is an equivalence.
	\end{prop}
	\begin{construction}\label{c: exact structure}
		The completed Hodge filtration gives rise to a fibre sequence
		\begin{equation}
			\infoneexcohcpl{-}{R}[1+n]\xrightarrow{d_{\Inf}}\what{F_2\mathbbb{\Pi}}_{-/R}[2+n]\xrightarrow{\incl} \infcohcpl{-}{R}[2+n].
		\end{equation}As an analogue to the classical de Rham cohomology theory, $\infoneexcohcpl{-}{R}[1+n]$ is called the chain complex of \textit{exact $2$-forms} of degree $n$. The \textit{classifying stack} of infinitesimal exact $2$-forms is defined as $\exformtwo{n}{R}{-}:=\Omega^{\infty}\infoneexcohcpl{-}{R}[1+n]$.
		
		The classifying stack of \textit{exact infinitesimal symplectic forms} of degree $n$ (relative to $R$) is then defined by the cartesian square % https://q.uiver.app/#q=WzAsNCxbMCwxLCJcXGV4Zm9ybXR3b3tufXtSfXstfSJdLFsxLDEsIlxcaW5mZm9ybXR3b3tufXtSfXstfSJdLFsxLDAsIlxcaW5mc3ltcHtufXtSfXstfSJdLFswLDAsIlxcaW5mc3ltcGV4e259e1J9ey19Il0sWzIsMV0sWzAsMSwifGRfe1xcSW5mfXwiLDJdLFszLDBdLFszLDJdLFszLDEsIiIsMSx7InN0eWxlIjp7Im5hbWUiOiJjb3JuZXIifX1dXQ==
		\[\begin{tikzcd}[ampersand replacement=\&]
			{\infsympex{n}{R}{-}} \& {\infsymp{n}{R}{-}} \\
			{\exformtwo{n}{R}{-}} \& {\infformtwo{n}{R}{-}.}
			\arrow[from=1-1, to=1-2]
			\arrow[from=1-1, to=2-1]
			\arrow["\lrcorner"{anchor=center, pos=0.125}, draw=none, from=1-1, to=2-2]
			\arrow[from=1-2, to=2-2]
			\arrow["{|d_{\Inf}|}"', from=2-1, to=2-2]
		\end{tikzcd}\]
	\end{construction}
	\begin{remark}
		We do not specify $\Inf$ in $\exformtwo{n}{R}{-}$, because $
		F_{[0,1]}c:\infoneexcohcpl{B}{R}\to F_{[0,1]}\mathrm{DR}_{B/R}$ is an equivalence. It suffices to note that both ends are the fibre of the universal differential
		\[d_{\deRham}\simeq d_{\Inf}:B\to \dL_{B/R}.\]
	\end{remark}
	This section is closed by a primary example of exact infinitesimal symplectic forms, which is analogous to the canonical symplectic form on the cotangent bundle in ordinary symplectic topology.
	\begin{exa}
		Let $X$ be a derived Deligne--Mumford stack finitely presented over $R\in\scr$. The tangent $\dT_{X/R}$ is perfect over $X$, so we may consider the cotrace map
		\[\ctr:\mathcal{O}_X\to \dT_{X/R}\otimes_{\mathcal{O}_X} \dL_{X/R}\]which is dual to $id_{\dT_{X/R}}$. Set $T_n:=\dT_{X/R}[-n]$ and $\cO_X[T_n]:=\lsym_{\cO_X}(\dT_{X/R}[-n])$. The relative non-connective spectrum $T^*_{X}[n]:=\spec_X(T_n)$ is then a derived Artin stack, and it is called the \textit{$n$-shifted cotangent stack} of $X$ (relative to $R$). There is a natural projection $p_n:T^*_X[n]\to X$, which gives rise to a fibre sequence $p^*_n\dL_{X/R}\to\dL_{T^*_X[n]/R}\to\dL_{T^*_X[n]/X}$. Then, the composite
		\begin{equation}
			\cO_X\xrightarrow{\ctr}T_n\otimes \dL_{X/R}[n]\to \cO_X[T_n]\otimes_{\cO_X}\dL_{X/R}[n]\to p_{n,*}p^*_n\dL_{X/R}[n]
		\end{equation}
		gives rise to a $1$-form $\Liouv_n$ of degree $n$ over $T^*_X[n]$, called \textit{the Liouville $1$-form (of degree $n$)}.
	\end{exa}
	\begin{prop}\label{prop: Liouv}
		The underlying $2$-form of $d_{\Inf}\Liouv_n$ is non-degenerate, so it is an infinitesimal symplectic form of degree $n$ over $T^*_X[n]$ relative to $R$.
	\end{prop}
	\begin{proof}
		The closedness of $d_{\Inf}\Liouv_n$ in $\infcohcplfil{T^*_X[n]}{R}$ is clear, since it is. It only remains to verify that the underlying $2$-form (still denoted as $d_{\Inf}\Liouv_n$) is non-degenerate. Without loss of generality, we can assume that $X=\spec S$ by passing to some \'etale covering of $X$. Now, set $\mathcal{Q}:=\infcohcplfil{S[T_n]}{R}$ and $\cA:=\infcohcplfil{S[T_n]}{S}$. There is a natural map of bi-filtered algebras
		\begin{equation}
			\infcohcplfil{S[T_n]}{R}\to \widehat{F\mathbbb{\Pi}}\mathstrut^{int}_{\cA/\mathcal{Q}},
		\end{equation}where $\infcohcplfil{S[T_n]}{R}$ has internal filtration degree $0$. Consider the connection map of the external filtration to obtain a square of filtered modules as follows (up to a shift),

		% https://q.uiver.app/#q=WzAsNCxbMCwwLCJcXGJpZ1xce1xcbGRvdHNcXHRvIDBcXHRvIFxcZExfe1NbVF9uXS9SfVstMl1cXGJpZ1xcfSJdLFswLDEsIlxcYmlnXFx7XFxsZG90c1xcdG8gMFxcdG9cXGxzeW1eMl97U1tUX25dfShcXGRMX3tTW1Rfbl0vUn1bLTFdKVxcYmlnXFx9Il0sWzEsMCwiXFxiaWdcXHtcXGxkb3RzXFx0byAwXFx0byBTW1Rfbl1cXG90aW1lc19TXFxkTF97Uy9SfVstMl1cXHRvIFxcZExfe1NbVF9uXS9SfVstMl1cXGJpZ1xcfSJdLFsxLDEsIlxcbHN5bV4yX3tTW1Rfbl19XFxiaWdcXHtcXGxkb3RzXFx0byAwXFx0byBTW1Rfbl1cXG90aW1lc19TXFxkTF97Uy9SfVstMV1cXHRvIFxcZExfe1NbVF9uXS9SfVstMV1cXGJpZ1xcfSJdLFswLDEsImRfe1xcSW5mfSJdLFswLDJdLFsxLDNdLFsyLDMsImRee2V4dH0iLDJdXQ==
		\[\begin{tikzcd}[ampersand replacement=\&, column sep=huge, scale cd=0.8]
			{\big\{\ldots\to 0\to \dL_{S[T_n]/R}[-2]\big\}} \& {\big\{\ldots\to 0\to S[T_n]\otimes_S\dL_{S/R}[-2]\to \dL_{S[T_n]/R}[-2]\big\}} \\
			{\big\{\ldots\to 0\to\lsym^2_{S[T_n]}(\dL_{S[T_n]/R}[-1])\big\}} \& {\lsym^2_{S[T_n]}\big\{\ldots\to 0\to S[T_n]\otimes_S\dL_{S/R}[-1]\to \dL_{S[T_n]/R}[-1]\big\}.}
			\arrow[from=1-1, to=1-2]
			\arrow["{d_{\Inf}}", from=1-1, to=2-1]
			\arrow["{d^{ext}}"', from=1-2, to=2-2]
			\arrow[from=2-1, to=2-2]
		\end{tikzcd}\]This induces an $S$-linear differential map
	\[\gr^{int}_1d^{ext}:S[T_n]\otimes_S\dL_{S/R}[n]\to S[T_n]\otimes_S (T_n\otimes_S\dL_{S/R})[n],\]which agrees with $d_{\Inf}:S[T_n]\to \dL_{S[T_n]/S}$, the universal differential of $S[T_n]$, tensored by the identity of $\dL_{S/R}[n]$ by Lemma \ref{lemma: integrable crystal}. Hence it induces an identity on the component $T_n\otimes_S\dL_{S/R}[n]$ on both sides. It implies that the following line
	\[\omega_n:R\to S\xrightarrow{\ctr}T_n\otimes_S \dL_{S/R}[n]\to S[T_n]\otimes_{S}\dL_{S/R}[n]\to \dL_{S[T_n]/R}[n]\xrightarrow{d_{\Inf}} \what{F_2\mathbbb{\Pi}}(S[T_n]/R)[n+2]\]exhibits a non-degenerate $2$-form over the non-connective algebra $S[T_n]$.
	Then, note that $S[T_n]\to p_{n,*}\cO_{T^*_X[n]}$ induces a comparison $\infcohcplfil{S[T_n]}{R}\to \infcohcplfil{T^*_X[n]}{R}$ and further a diagram% https://q.uiver.app/#q=WzAsNSxbMSwwLCJcXGRMX3tTW1Rfbl0vUn0iXSxbMSwxLCJwX3tuLCp9XFxkTF97VF4qX1hbbl0vUn0iXSxbMiwwLCJcXGxzeW1eMl97U1tUX25dfVxcZExfe1NbVF9uXS9SfSJdLFsyLDEsInBfe24sKn1cXGxzeW1eMl97XFxjT197VF4qX1hbbl19fVxcZExfe1ReKl9YW25dL1J9Il0sWzAsMCwiXFxvbWVnYV9uOlJbLW5dIl0sWzAsMV0sWzAsMiwiZF97XFxJbmZ9Il0sWzEsMywicF97biwqfWRfe1xcSW5mfSJdLFsyLDNdLFs0LDBdLFs0LDEsInBfe24sKn1cXExpb3V2X24iLDJdXQ==
	\[\begin{tikzcd}[ampersand replacement=\&,column sep=2cm]
		{\omega_n:R[-n]} \& {\dL_{S[T_n]/R}} \& {\lsym^2_{S[T_n]}\dL_{S[T_n]/R}} \\
		\& {p_{n,*}\dL_{T^*_X[n]/R}} \& {p_{n,*}\lsym^2_{\cO_{T^*_X[n]}}\dL_{T^*_X[n]/R}.}
		\arrow[from=1-1, to=1-2]
		\arrow["{p_{n,*}\Liouv_n}"', from=1-1, to=2-2]
		\arrow["{d_{\Inf}}", from=1-2, to=1-3]
		\arrow[from=1-2, to=2-2]
		\arrow[from=1-3, to=2-3]
		\arrow["{p_{n,*}d_{\Inf}}", from=2-2, to=2-3]
	\end{tikzcd}\]Moreover, it is straightforward to verify that $p^*_{n}T_n\simeq \dL_{T^*_X[n]/S}$, and then $\dL_{T^*_X[n]/R}\simeq p^*_n(T_n\oplus\dL_{S/R})$. Hence, the above diagram shows that $d_{\Inf}\Liouv_n$ is non-degenerate.
	\end{proof}
	\begin{lemma}\label{lemma: integrable crystal}
		Let $S\to S'$ be an arrow in $\scr_R$. Consider $\cA:=\infcohcplfil{S'}{S}$ and $\mathcal{Q}:=\infcohcplfil{S'}{R}$ , the algebraically integrable foliation-like algebras. Then, there is an equivalence of filtered $\cA$-modules
		\begin{equation}
			\gr^{int}_n\widehat{F\mathbbb{\Pi}}\mathstrut^{int}_{\cA/\mathcal{Q}}\simeq \lsym^n_{S}(\dL_{S/R}[-1])\widehat{\otimes}_{S} \cA(n).
		\end{equation}
	\end{lemma}
	
	\begin{proof}
		Denote $\mathcal{Q}_0:=\infcohcplfil{S}{R}$. Consider the cocartesian square in $\widehat{\fil}\dalg_R$
		% https://q.uiver.app/#q=WzAsNCxbMCwxLCJTIl0sWzEsMSwiXFxpbmZjb2hjcGxmaWx7Uyd9e1N9Il0sWzAsMCwiXFxpbmZjb2hjcGxmaWx7U317Un0iXSxbMSwwLCJcXGluZmNvaGNwbGZpbHtTJ317Un0iXSxbMCwxXSxbMiwwXSxbMiwzXSxbMywxXV0=
		\[\begin{tikzcd}[ampersand replacement=\&]
			{\infcohcplfil{S}{R}} \& {\infcohcplfil{S'}{R}} \\
			S \& {\infcohcplfil{S'}{S}}
			\arrow[from=1-1, to=1-2]
			\arrow[from=1-1, to=2-1]
			\arrow[from=1-2, to=2-2]
			\arrow[from=2-1, to=2-2]
		\end{tikzcd}\]This induces an equivalence of completed bi-filtered algebras
		\[\widehat{F\mathbbb{\Pi}}\mathstrut^{int}_{S/\mathcal{Q}_0}\widehat{\otimes}_{\mathcal{Q}_0}\mathcal{Q}\simeq \widehat{F\mathbbb{\Pi}}\mathstrut^{int}_{\cA/\mathcal{Q}}.\]Taking $\gr^{int}$-pieces, we obtain
		\[\gr^{int}_n\widehat{F\mathbbb{\Pi}}\mathstrut^{int}_{\cA/\mathcal{Q}}\simeq \lsym_S (\dL_{S/R}[-1](1))\otimes_S S\what{\otimes}_{\mathcal{Q}_0}\mathcal{Q}\simeq \lsym^n_{S}(\dL_{S/R}[-1])\widehat{\otimes}_{S} \cA(n).\]
	\end{proof}	
	\section{Formal integration of infinitesimal symplectic forms}\label{sec: Formal integration of infinitesimal symplectic forms}

	\subsection{Local Exactness modulo $\res_{JF}$}\label{sec: 3.1 local exactness}
	\begin{theorem}[Poincar\'e type lemma]\label{theorem: computing inf coh of B/k}
		Let $k$ be a field of characteristic $p>0$, and let $B\in \scr_k$ be an almost finitely presented algebra. Then, there are natural equivalences
		\begin{equation}
			\lim_n\pi_0 B^{(n)}\xleftarrow{\simeq}\lim_n B^{(n)}\xrightarrow{\simeq} \infcohcpl{B}{k},
		\end{equation}where $\lim\limits_{n}B^{(n)}$ is the right perfection relative to $k$.
	\end{theorem}
	
	Let $A\to B$ be a map of almost finitely presented simplicial commutative $k$-algebras, and suppose that $\pi_0(A)\twoheadrightarrow\pi_0(B)$ is surjective. We regard $\infcohcplfil{B}{A}$ as a \textit{derived adic completion} of $A$ along $A\to B$ in the spirit of Proposition \ref{prop: inf foliation basics}(1). In particular, denote the derived adic completion along the diagonal $B^{\otimes_k n}\to B$ by
	\begin{equation}\label{e: what is an adic fil of diagonal}
		B^{\widehat{\otimes}_k n}_{\adic}:= \widehat{F\mathbbb{\Pi}}(B/B^{\otimes_k n}).
	\end{equation}
	
	We are going to see that line (\ref{e: what is an adic fil of diagonal}) can serve as a cosimplicial resolution of $\infcohcpl{B}{k}$. Nevertheless, this resolution is too large to compute explicitly. For instance, $B\widehat{\otimes}_kB$ is the $\infty$-jet sheaf of $\spec B$, and its $\mathbb{E}_1$-Koszul duality encodes Grothendieck's differential operator algebra over $B$, which is known to be huge in the mod $p$ case, even when $B$ is smooth \cite{smith2006differential}. We can fix this problem by approximating $k\to B$ by the tower $B^{(n)}\to B$.

	\begin{lemma}\label{lemma: cot of frob tower}
		 The cotangent complex $\dL_{B/B^{\otimes_{B^{(n)}}m}}$ ($m\ge 1$) is given by
		 \begin{equation}\label{e: cot of Frob tower}
		 	\dL_{B/k}^{\oplus m-1}\oplus\Sigma \varphi^{n,*}_B\dL_{B^{(n)}/k}^{\oplus m-1}\simeq \dL_{B/k}^{\oplus m-1}\oplus \Sigma\frob^{\circ n,*}(\dL_{B/k}^{\oplus m-1}) .
		 \end{equation}
		 Additionally, the natural map $\dL_{B/B^{\otimes_{B^{(n+l)}}m}}\to \dL_{B/B^{\otimes_{B^{(n)}}m}}$ is given by the identity on the copies of $\dL_{B/k}$ and is null on $\varphi^{n+l,*}_B\dL_{B^{(n+l)}/k}^{\oplus m-1}$, where $l\ge1$.
	\end{lemma}
	\begin{proof}
		When $B$ is a finitely generated $k$-polynomial, there is always an identically null map 
		\[\frob^*\Omega_{B/k}\xrightarrow{0=d\frob}\Omega_{B/k}.\]Repeating this null map $n$ times and taking the derived natural transformation shows that the natural map
		\[\frob^{\circ n,*}\dL_{B/k}\to \dL_{B/k}\]is null-homotopic. Then, the formula (\ref{e: cot of Frob tower}) is given by repeatedly applying the fundamental fibre sequence of the cotangent complexes.
	\end{proof}
	\begin{lemma}\label{lemma: completeness of frob jet}
		For $n\ge 0$,  there is a natural equivalence
		\[B^{\otimes_{B^{(n)}}\bullet}\simeq |\what{F\mathbbb{\Pi}}(B/B^{\otimes_{B^{(n)}}\bullet})|.\]
		For this reason, we denote $\what{F\mathbbb{\Pi}}(B/B^{\otimes_{B^{(n)}}\bullet})$ by $B^{{\otimes}_{B^{(n)}}\bullet}_{\adic}$.
	\end{lemma}
	\begin{proof}
		The kernel of $\pi_0B^{\otimes_{B^{(n)}}\bullet}\to \pi_0B$ is killed by the ordinary $p^n$-power, and the cotangent $\dL_{B/B^{\otimes_{B^{(n)}}\bullet}}$ is almost perfect as a $B$-module. Then, this lemma follows from Proposition \ref{prop: inf foliation basics}.
	\end{proof}
	
	\begin{lemma}\label{lemma: compare cosimplicial resol s}
		There is a commuting square in $\widehat{\fil}\dalg_k$ as follows% https://q.uiver.app/#q=WzAsNCxbMCwwLCJcXGluZmNvaGNwbGZpbHtCfXtrfSJdLFswLDEsIlxcdG90IEJee1xcd2lkZWhhdHtcXG90aW1lc31fayBcXGJ1bGxldCsxfV97XFxhZGljfSJdLFsxLDAsIlxcaW5mY29oY3BsZmlse0J9e0JeeyhuKX19Il0sWzEsMSwiXFx0b3QgQl57XFxvdGltZXNfe0JeeyhuKX19IFxcYnVsbGV0KzF9X3tcXGFkaWN9Il0sWzAsMSwiXFxzaW1lcSIsMl0sWzIsMywiXFxzaW1lcSIsMl0sWzAsMl0sWzEsM11d
		\[\begin{tikzcd}[ampersand replacement=\&]
			{\infcohcplfil{B}{k}} \& {\infcohcplfil{B}{B^{(n)}}} \\
			{\tot B^{\widehat{\otimes}_k \bullet+1}_{\adic}} \& {\tot B^{\otimes_{B^{(n)}} \bullet+1}_{\adic},}
			\arrow[from=1-1, to=1-2]
			\arrow["\simeq"', from=1-1, to=2-1]
			\arrow["\simeq"', from=1-2, to=2-2]
			\arrow[from=2-1, to=2-2]
		\end{tikzcd}\]where the vertical arrows are equivalences, which is functorial in $n\ge 0$.
	\end{lemma}
	\begin{proof}
		There is a natural augmented cosimplicial diagram
		\[\infcohcplfil{B}{k}\to \widehat{F\mathbbb{\Pi}}(B/B^{\otimes_k \bullet+1}),\]whose $\gr_1$-piece is equivalent to the cogroupoid object generated by $\dL_{B/k}\to 0$ in $\aperf_B$. Thus, this is in fact a totalization diagram of foliation-like algebras, see Proposition \ref{prop: inf foliation basics} (2). Then, line \ref{e: what is an adic fil of diagonal} implies the equivalence on the left. The right equivalence can be proved in the same way. 
	\end{proof} 
	\begin{lemma}\label{lemma: approximating the jets}
		Fix $m\ge 1$. There is an equivalence $B^{\widehat{\otimes}_k m}_{\adic}\simeq \lim_n B^{\otimes_{B^{(n)}} m}_{\adic}$ in $\what{\fil}\dalg_k$.
	\end{lemma}
	\begin{proof}
		The diagram $B^{\otimes_{B^{(\bullet)}} m}_{\adic}$ is a tower of connective and almost perfect foliation-like algebras. Observe that their $\gr_1$-pieces form a tower of almost perfect complexes,
		\[\ldots\to \dL_{B/k}^{\oplus m-1}\oplus\Sigma \varphi^{n+1,*}_B\dL_{B^{(n+1)}/k}^{\oplus m-1}\to \dL_{B/k}^{\oplus m-1}\oplus\Sigma \varphi^{n,*}_B\dL_{B^{(n)}/k}^{\oplus m-1}\to \dL_{B/k}^{\oplus m-1}\oplus\Sigma \varphi^{n-1,*}_B\dL_{B^{(n-1)}/k}^{\oplus m-1}\to \ldots\]whose limit is the almost perfect $B$-module $\dL^{\oplus m-1}_{B/k}$, see Lemma \ref{lemma: cot of frob tower}. Hence, by Proposition \ref{prop: inf foliation basics} (2), $\lim_n B^{\otimes_{B^{(n)}} m}_{\adic}$ is also an almost perfect and connective foliation-like algebra, and the natural map from $B^{\widehat{\otimes}_k m}_{\adic}$ is an equivalence.
	\end{proof}
	
	\begin{lemma}\label{lemma: bastard}
		Suppose that $A\in\scr_k$ is almost finitely presented (Definition \ref{df: ap et afp}). Let $B$ be an $A$-algebra satisfying
		\begin{enumerate}[label=(\arabic*)]
			\item its underlying $A$-module is almost perfect,
			\item the map of sets $\spec(\pi_0 B)\to \spec (\pi_0 A)$ is a surjection.
		\end{enumerate}
		Then, the natural map of derived algebras $A\to \tot(B^{\otimes_A\bullet+1})$ is an equivalence.
		
	\end{lemma}
	\begin{proof}
	Let $\tau_{\le\bullet}A\to \tau_{\le\bullet}B$ be the Postnikov tower of $A\to B$. At level $0$, the map $\pi_0A\to \pi_0B$ is descendable in the sense of \cite[Definition 3.17]{mathew2016galois} by \cite[Proposition 11.25]{bhatt2017projectivity}. Since $\tau_{\le n}A\to \pi_0A$ is descendable \cite[Proposition 3.32]{mathew2016galois}, the map $\tau_{\le n}A\to \tau_{\le n}B$ is also descendable by \cite[Proposition 3.23]{mathew2016galois}. Hence we have an equivalence
	\begin{equation}\label{e: bastard truncated}
		\tau_{\le n}A\xrightarrow{\simeq}\tot\tau_{\le\bullet}B^{\otimes_{\tau_{\le\bullet}A}(\bullet+1)}.
	\end{equation} 
	Then, note that the tower $\tau_{\le\bullet}B^{\otimes_{\tau_{\le\bullet}A}(m+1)}$ stabilizes at every homotopy degree, and furthermore, the limit is $B^{\otimes_Am+1}$. Therefore, taking limits in (\ref{e: bastard truncated}) gives $A\simeq \tot B^{\otimes_A\bullet+1}$.
	\end{proof}

	\begin{lemma}\label{lemma: computing infcoh of B/B^(n)}
		The natural map $B^{(n)}\to \infcohcpl{B}{B^{(n)}}$ is an equivalence.
	\end{lemma}
	\begin{proof}
	Since $\pi_0(B^{(n)})$ is the classical $n$-th Frobenius twist of $\pi_0(B)$, Lemma \ref{lemma: bastard} gives $B^{(n)}\simeq \tot B^{\otimes_{B^{(n)}}\bullet +1}$. Then, conclude by Lemma \ref{lemma: compare cosimplicial resol s}.
	\end{proof}
	\begin{proof}[Proof of Theorem \ref{theorem: computing inf coh of B/k}]
		Consider the following diagram with a commutative square inside,% https://q.uiver.app/#q=WzAsNixbMCwwLCJcXGluZmNvaGNwbHtCfXtrfSJdLFswLDEsIlxcdG90IEJee1xcd2lkZWhhdHtcXG90aW1lc31fayBcXGJ1bGxldCsxfSJdLFsxLDAsIlxcbGltX25cXGluZmNvaGNwbHtCfXtCXnsobil9fSJdLFsxLDEsIlxcbGltX25cXHRvdCBCXntcXG90aW1lc197Ql57KG4pfX0gXFxidWxsZXQrMX0iXSxbMiwwLCJcXGxpbV9uQl57KG4pfSJdLFszLDAsIlxcbGltX25cXHBpXzBCXnsobil9Il0sWzAsMSwiXFxzaW1lcSwgXFx0ZXh0eyhcXHJlZntsZW1tYTogY29tcGFyZSBjb3NpbXBsaWNpYWwgcmVzb2wgc30pfSIsMl0sWzIsMywiXFxzaW1lcSwgXFx0ZXh0eyhcXHJlZntsZW1tYTogY29tcGFyZSBjb3NpbXBsaWNpYWwgcmVzb2wgc30pfSIsMl0sWzAsMl0sWzEsMywiXFxzaW1lcSwgXFx0ZXh0eyhcXHJlZntsZW1tYTogYXBwcm94aW1hdGluZyB0aGUgamV0c30pfSJdLFs0LDIsIlxcc2ltZXEsIFxcdGV4dHsoXFxyZWZ7bGVtbWE6IGNvbXB1dGluZyBpbmZjb2ggb2YgQi9CXihuKX0pfSIsMl0sWzQsNSwiXFxzaW1lcSwgKFxccmVme2xlbW1hOiBubyBuZWVkIGZvciBoaWdoZXIgZ3JvdXBzIG9mIEZyb2IgdHdpc3RzfSkiXV0=
		\[\begin{tikzcd}[ampersand replacement=\&]
			{\infcohcpl{B}{k}} \& {\lim_n\infcohcpl{B}{B^{(n)}}} \& {\lim_nB^{(n)}} \& {\lim_n\pi_0B^{(n)}} \\
			{\tot B^{\widehat{\otimes}_k \bullet+1}} \& {\lim_n\tot B^{\otimes_{B^{(n)}} \bullet+1}.}
			\arrow[from=1-1, to=1-2]
			\arrow["{\simeq, \text{(\ref{lemma: compare cosimplicial resol s})}}"', from=1-1, to=2-1]
			\arrow["{\simeq, \text{(\ref{lemma: compare cosimplicial resol s})}}"', from=1-2, to=2-2]
			\arrow["{\simeq, \text{(\ref{lemma: computing infcoh of B/B^(n)})}}"', from=1-3, to=1-2]
			\arrow["{\simeq, (\ref{lemma: no need for higher groups of Frob twists})}", from=1-3, to=1-4]
			\arrow["{\simeq, \text{(\ref{lemma: approximating the jets})}}", from=2-1, to=2-2]
		\end{tikzcd}\]
		I summerize the known equivalences and mark the number of the lemmas that they rely on.
		\end{proof}

			We now apply Theorem \ref{theorem: computing inf coh of B/k} to study the infinitesimal closed $2$-forms. By the fibre sequence \[\infoneexcohcpl{B}{R}\xrightarrow{d_{\Inf}}\infeqtrcohcpl{B}{R}{2}[1]\xrightarrow{\incl}\infcohcpl{B}{R}[1],\]an exact structure of some $\omega\in \infformtwo{-1}{R}{\spec B}$ amounts to a null-homotopy $\incl\omega\sim 0$. The group $\pi_{-1}\infcohcpl{B}{R}\cong \lim^1_n \pi_0B^{(n)}$ then contains the obstructions to $\omega$ being exact. We shall see that these obstructions can be viewed heuristically as formal series in the Frobenius.
			\begin{exa}\label{exa: jet frobenius}
				Let $S:=k[x_1,\ldots,x_t]$ be a polynomial ring. The group $\pi_{-1}\infcohcpl{S}{k}$ is a quotient of $\prod_{n\ge0}k[x^{p^n}_1,\ldots,x^{p^n}_t]$, where two sequences $(f_n(X^{p^n}))$ and $(g_n(X^{p^n}))$ are equivalent if and only if $f_n=g_n$ for all but finite $n$. The cosimplicial resolution of $\infcohcplfil{S}{k}$ in Lemma \ref{lemma: compare cosimplicial resol s} is given by
				\[k[x_i]\xrightarrow{\partial_0}\big(k[x_i]\otimes_kk[x_i]\big)^\wedge_{\Delta_1}\xrightarrow{\partial_1}\big(k[x_i]\otimes_kk[x_i]\otimes_kk[x_i]\big)^\wedge_{\Delta_2}\to \cdots\]carrying the degreewise adic filtration, where $\Delta_j$ is the kernel of $S^{\otimes_k j+1}\to S$. In this context, the cycle $x^{p^n}_i\otimes 1-1\otimes x^{p^n}_i\in \pi_{-1}\infeqtrcohcpl{S}{k}{p}$ equals $(\mathscr{P}^0)^{\circ n}(v_i)$, where $v_i=x_i\otimes 1-1\otimes x_i\in \pi_{-1}\infeqtrcohcpl{S}{k}{1}$, and $\mathscr{P}^0$ is the Dyer--Lashof operation of degree $0$. Thus every cycle of $\pi_{-1}\infcohcpl{S}{k}$ can be represented by an infinite sum\footnote{If $k$ is perfect, the coefficients $a_n$ are redundant.}
				\[\sum_{n}a_n\big(\mathscr{P}^0\big)^{\circ n}(h_n)\] where $h_n$ is a polynomial in $v_i$. Intuitively, the group $\pi_{-1}\infcohcpl{S}{k}$ consists of classes of formal power series in Frobenius in the $\infty$-jet $ S\what{\otimes}_k S$.
			\end{exa}
			
			\begin{construction}\label{c: jF power map}
				Given a finitely generated polynomial ring $S$ over $k$, the \v{C}ech nerve $S^{\otimes_k\bullet+1}$ carries a degreewise ordinary adic filtration given by the augmentation $S^{\otimes_k\bullet+1}\to S$ to the constant diagram $S$. Denote by $C(S)$ this filtered cosimplicial ring. The degreewise completion $\widehat{C}(S)$ is an explicit model for the resolution of $\infcohcplfil{S}{k}$ in Lemma \ref{lemma: compare cosimplicial resol s}, since $k\to S$ is smooth. For every $n\ge 1$, there is a map of cosimplicial commutative rings natural in $S$
				\[S^{(n)}\to F_{[0,p^n-1]}C(S),\]where $S^{(n)}$ is a constant diagram. Hence we obtain an inclusion map
				\begin{equation}
					\mathbbb{i}_n:B^{(n)}\to F_{[0,p^{n}-1]}\mathbbb{\Pi}_{B/k}
				\end{equation}
				for every $B\in\scr_k$ by left Kan extension.  Furthermore, there is a family of maps
				\[\mathbbb{j}_n:B^{(n)}\xrightarrow{\mathbbb{i}}F_{[0,p^{n}-1]}\mathbbb{\Pi}_{B/k}\xrightarrow{d_{\Inf}}\what{F_{p^n}\mathbbb{\Pi}}_{B/k}[1],\mathrlap{\quad (n\ge 1)}\]that plays an analogous role to $(\mathscr{P}^0)^{\circ n}(v_i)$ in Example \ref{exa: jet frobenius}. By projecting to $F_2$, we obtain
				\begin{equation}
				\cJ:\prod_{n\ge 1}B^{(n)}\xrightarrow{ \mathbbb{j}_1+\mathbbb{j}_2+\cdots}\infeqtrcohcpl{B}{k}{2}[1],
				\end{equation}which we call the \textit{jet-Frobenius (JF) power map}.
			\end{construction} 
			\begin{remark}
				When $B$ is almost finitely presented over $k$, the composite
				\[\prod_{n\ge1}\pi_0 B^{(n)}\xrightarrow{\pi_0\cJ}\pi_{-1}\infeqtrcohcpl{B}{k}{2}\xrightarrow{\incl}\pi_{-1}\infcohcpl{B}{k}\]agrees with the map $\prod_{n\ge1}\pi_0 B^{(n)}\twoheadrightarrow \lim^1_n\pi_0 B^{(n)}$ in Milnor's sequence via Theorem \ref{theorem: computing inf coh of B/k}.
			\end{remark}
			\begin{df}
			Fix $\omega\in\infformtwo{-1}{k}{\spec B}$ an infinitesimal closed $2$-form of degree $-1$. We say a cycle $\res_{JF}\in \pi_{-1}\infeqtrcohcpl{B}{k}{2}$ is a \textit{jet-Frobenius residue} of $\omega$ if
			\begin{enumerate}[label=(\arabic*)]
				\item there exists some $g\in\prod\limits_{\mathbf{n\ge2}}\pi_0 B^{(n)}$ such that $\pi_0\cJ(g)=\res_{JF}$;\vspace{-1em}
				\item we have $\incl(\omega)=\incl(\res_{JF})$ in $\pi_{-1}\infcohcpl{B}{k}$.
			\end{enumerate}Here we say that $g$ represents $\res_{JF}$.
			\end{df}
			
	\begin{corollary}\label{cor: wJF infinite degree}
		Suppose that $B$ is almost of finite presentation over $k$, and there is an infinitesimal symplectic form $\omega$ of degree $-1$. Then there exist a jet-Frobenius residue $\res_{JF}$ and an exact structure $\lambda\in \infsympex{-1}{\spec B}{k}$ such that $\omega\sim d_{\Inf}\lambda+\res_{JF}$ in $\infsymp{-1}{\spec B}{k}$.
	\end{corollary}
	\begin{proof}
		The kernel of $\prod_{n\ge1}\pi_0 B^{(n)}\twoheadrightarrow \lim^1_n\pi_0 B^{(n)}$ contains all finite sequences, so it restricts to a surjection $\prod_{n\ge N}\pi_0 B^{(n)}\twoheadrightarrow \lim^1_n\pi_0 B^{(n)}$ for every $N$. Hence there exists a $g\in \prod_{n\ge2}\pi_0 B^{(n)}$ such that $\incl\cJ(g)=\incl\omega$. Since $\res_{JF}:=\cJ(g)$ lifts to $\what{F_{p^2}\mathbbb{\Pi}}_{B/k}$, the underlying $2$-form of $\res_{JF}$ is null-homotopic, and therefore $\omega-\res_{JF}$ remains infinitesimally symplectic.
	\end{proof}
	
	The jet-Frobenius residue is not unique. We will return to this point in \S\ref{s5.3: Algebraization of Darboux charts} and discuss its influence on the construction of Darboux charts.

	\subsection{Lagrangian foliations and formal quotients}\label{s:3.2 Lagrangian foliation}
	The subsection aims to show that, Zariski locally, the $(-1)$-shifted infinitesimal symplectic derived schemes support integrable Lagrangian foliations.
	\begin{df}
		Let $L$ be a perfect complex over $B\in\scr$. Assume that $\omega\in \pi_n(\wedge^2_B L)$ is a $2$-form of degree $-n$. Then we define:
		\begin{enumerate}[label=(\arabic*)]
			\item A \textit{distribution} of $\omega$ is a morphism $l:L\to E$ in $\perf_B$.
			\item This distribution is \textit{isotropic} if the composition
			\[B[n]\xrightarrow{\omega}\wedge^2_BL\xrightarrow{\wedge^2 l} \wedge^2_B E\]is null-homotopic. An \textit{isotropy structure} of $l$ is a 2-cell as follows % https://q.uiver.app/#q=WzAsNSxbMCwwLCJFXlxcdmVlIl0sWzEsMCwiTF5cXHZlZSJdLFsyLDAsIkxbLW5dIl0sWzIsMSwiRVstbl0iXSxbMCwxLCIwIl0sWzAsMSwibF5cXHZlZSJdLFsxLDIsIlxcVGhldGFfXFxvbWVnYSJdLFsyLDMsImwiXSxbMCw0XSxbNCwzXSxbNCwyLCJcXHNpbWVxIiwxLHsic3R5bGUiOnsiYm9keSI6eyJuYW1lIjoibm9uZSJ9LCJoZWFkIjp7Im5hbWUiOiJub25lIn19fV1d
			\[\begin{tikzcd}[ampersand replacement=\&]
				{E^\vee} \& {L^\vee} \& {L[-n]} \\
				0 \&\& {E[-n]}
				\arrow["{l^\vee}", from=1-1, to=1-2]
				\arrow[from=1-1, to=2-1]
				\arrow["{\Theta_\omega}", from=1-2, to=1-3]
				\arrow["l", from=1-3, to=2-3]
				\arrow["\simeq"{description}, draw=none, from=2-1, to=1-3]
				\arrow[from=2-1, to=2-3]
			\end{tikzcd}\]
			\item Suppose that $\omega$ is non-degenerate. An isotropy structure of $l$ is \textit{Lagrangian} if the induced map $\mathrm{cofib}(l^\vee)\to E[-n]$ is an equivalence.
		\end{enumerate}
	\end{df}
	\begin{remark}[Yoga of Lagrangian isotropy]\label{rk: Lagrangian isotropy as a non-degenerated matrix}
	Let $\omega$ be as in (3) of this definition. Suppose that the isotropy structure is determined by a lift $\tilde{\omega}$ of $\omega$ in $F:=\mathrm{fib}(\wedge^2l)$. There is a fibre sequence \[\wedge^2\mathrm{fib}(l)\to F\xrightarrow{r} E\otimes_B \mathrm{fib}(l).\] Then, the composite $r\tilde{\omega}$ is dual to the equivalence $\mathrm{cofib}(l^\vee)\to E[-n]$.
	
	Indeed, the natural transformation $\wedge^2_B(-)\to (-)^{\otimes_B 2}$ by skew-symmetrization induces a morphism of fibre sequences% https://q.uiver.app/#q=WzAsNixbMCwwLCJGIl0sWzAsMSwiXFxtYXRocm17ZmlifShsKVxcb3RpbWVzX0JFIl0sWzEsMSwiTFxcb3RpbWVzX0JFIl0sWzEsMCwiXFx3ZWRnZV4yX0JMIl0sWzIsMCwiXFx3ZWRnZV4yX0IgRSJdLFsyLDEsIkVcXG90aW1lc19CRSwiXSxbMCwxLCJyIl0sWzEsMl0sWzAsM10sWzMsNF0sWzMsMl0sWzQsNV0sWzIsNV1d
	\[\begin{tikzcd}[ampersand replacement=\&]
		\mathllap{(B[n]\xlongrightarrow{\tilde{\omega}})\quad}F \& {\wedge^2_BL} \& {\wedge^2_B E} \\
		{\mathrm{fib}(l)\otimes_BE} \& {L\otimes_BE} \& {E\otimes_BE,}
		\arrow[from=1-1, to=1-2]
		\arrow["r", from=1-1, to=2-1]
		\arrow[from=1-2, to=1-3]
		\arrow[from=1-2, to=2-2]
		\arrow[from=1-3, to=2-3]
		\arrow[from=2-1, to=2-2]
		\arrow[from=2-2, to=2-3]
	\end{tikzcd}\]where the middle downward arrow is the composition $\wedge^2_BL\to L\otimes_B L\to L\otimes_B E$. The image of $\tilde{\omega}$ in $L\otimes_B E$ is dual to $l\Theta_\omega$, and its image in $E\otimes_B E$ is dual to $l\Theta_\omega l^\vee\simeq 0$. Hence the composition $r\tilde{\omega}$ is dual to $\mathrm{cofib}(l^\vee)[n]\simeq E$. 		
	\end{remark}
	\begin{construction}[Favorable Lagrangian distribution]\label{c: Lagrangian distribution}
		Let $B$ be a finitely presented $k$-algebra, and fix a schematic point $x\in \spec \pi_0B$. The image of the ordinary differential
		\[\pi_0 B\xrightarrow{d_{\deRham}}\Omega_{\pi_0B/k}\otimes_B \kappa(x)\] spans the target after base change to $\kappa(x)$, where $\kappa(x)$ denotes the residue field. Assume that $\dim_{\kappa(x)}\Omega_{\pi_0B/k}\otimes_B \kappa(x)=r$. After replacing $\spec B$ by a Zariski neighborhood of $x$ (still denoted $\spec B$), we can choose $b_i\in \pi_0 B$ ($1\le i\le r$) so that $\pi_0\dL_{B/k}$ is spanned by $d_{\deRham}b_i$.  For later convenience, take an $r$-tuple $T=(t_i)$ and consider the commutative diagram
		\begin{equation}\label{e:Lag dist}
			\begin{tikzcd}[ampersand replacement=\&]
				{k\linsp{T}} \& {B\linsp{T}\mathrlap{\ \ =:E}} \\
				B \& {\dL_{B/k}.}
				\arrow[from=1-1, to=1-2]
				\arrow["{(b_i)}"', from=1-1, to=2-1]
				\arrow["{(d_{\deRham}b_i)}", from=1-2, to=2-2]
				\arrow["{d_{\deRham}}"', from=2-1, to=2-2]
			\end{tikzcd}
		\end{equation}

		 Assume that there exists $\omega_0\in\ndformtwo{-1}{k}{B}$. Then, the cotangent complex $\dL_{B/k}$ has Tor-amplitude $[0,1]$, and $C:=\mathrm{fib}(E\to \dL_{B/k})$ has Tor-amplitude $0$. The map induced by $\omega_0$
		\[\theta:C^\vee\to \dT_{B/k}[1]\xrightarrow{\Theta_{\omega_0}}\dL_{B/k}\to C[1]\]is null homotopic for degree reason. \textit{By fixing a null homotopy} $\theta\sim0$, we obtain a dotted filler% https://q.uiver.app/#q=WzAsNCxbMCwxLCJFIl0sWzEsMSwiXFxkTF97Qi9rfSwiXSxbMSwwLCJcXGRUX3tCL2t9WzFdIl0sWzAsMCwiQ15cXHZlZSJdLFsyLDEsIlxcVGhldGFfe1xcb21lZ2FfMH0iXSxbMCwxXSxbMywyXSxbMywwLCIiLDIseyJzdHlsZSI6eyJib2R5Ijp7Im5hbWUiOiJkYXNoZWQifX19XV0=
		\[\begin{tikzcd}[ampersand replacement=\&]
			{C^\vee} \& {\dT_{B/k}[1]} \\
			E \& {\dL_{B/k}.}
			\arrow[from=1-1, to=1-2]
			\arrow[dashed, from=1-1, to=2-1]
			\arrow["{\Theta_{\omega_0}}", from=1-2, to=2-2]
			\arrow[from=2-1, to=2-2]
		\end{tikzcd}\]This produces  an isomorphism of $\kappa(x)$-vector spaces $C^\vee\otimes_B\kappa(x)\cong E\otimes_B\kappa(x)$. Therefore, we can replace $C$ by $E^\vee[1]$ after shrinking $\spec B$ around $x$ in Zariski topology.
		% https://q.uiver.app/#q=WzAsNCxbMCwxLCJFIl0sWzEsMSwiXFxkTF97Qi9rfSJdLFsxLDAsIlxcZFRfe0Iva31bMV0iXSxbMiwxLCJFXlxcdmVlWzFdIl0sWzAsMV0sWzAsMl0sWzIsM10sWzEsM10sWzIsMSwiXFxzaW1lcSIsMV1d
		\[\begin{tikzcd}[ampersand replacement=\&]
			\& {\dT_{B/k}[1]} \\
			E \& {\dL_{B/k}} \& {E^\vee[1]}
			\arrow["\simeq"{description}, from=1-2, to=2-2]
			\arrow[from=1-2, to=2-3]
			\arrow[from=2-1, to=1-2]
			\arrow[from=2-1, to=2-2]
			\arrow[from=2-2, to=2-3]
		\end{tikzcd}\]
		This exhibits $\dL_{B/k}\to E^\vee[1]$ as a Lagrangian distribution.
	\end{construction}
	\begin{df}Let $\infcohcplfil{B}{k}\to \cA$ denote a foliation-like algebra, and suppose that there exists $\omega\in \infsymp{-n}{k}{\spec B}$. We say that $\cA$ is a \textit{Lagrangian infinitesimal foliation} with respect to $\omega$ if the following conditions hold:
	\begin{enumerate}[label=(\arabic*)]
		\item the cycle $k[n-2]\to \infeqtrcohcpl{B}{k}{2}\to F_2\cA$ is null-homotopic;
		\item the distribution $\dL_{B/k}\to \dL_\cA$ is Lagrangian for the underlying $2$-form of $\omega$.
	\end{enumerate}
		\end{df}
	\begin{exa}\label{exa: f&lambda}
		Let $\rho:\infcohcplfil{B}{k}\to \cA$ be a foliation-like algebra, and let $\lambda$ be an exact $2$-form of degree $n$. An \textit{isotropic structure} of $\rho$ with respect to $d_{\Inf}\lambda$ amounts to a null homotopy $\rho d_{\Inf}\lambda\sim 0$, or equivalently, a dotted lifting as follows
		% https://q.uiver.app/#q=WzAsMyxbMCwxLCJrWy1uLTFdIl0sWzEsMSwiRl97WzAsMV19XFxjQSJdLFsxLDAsInxcXGNBfCJdLFsyLDFdLFswLDEsIlxccmhvXFxsYW1iZGEiLDJdLFswLDIsImYiLDAseyJzdHlsZSI6eyJib2R5Ijp7Im5hbWUiOiJkYXNoZWQifX19XV0=
		\[\begin{tikzcd}[ampersand replacement=\&]
			\& {|\cA|} \\
			{k[-n-1]} \& {F_{[0,1]}\cA.}
			\arrow[from=1-2, to=2-2]
			\arrow["f", dashed, from=2-1, to=1-2]
			\arrow["{\rho\lambda}"', from=2-1, to=2-2]
		\end{tikzcd}\]
	\end{exa}
	\begin{prop}\label{prop: Lag Foliation}
		Following Construction \ref{c: Lagrangian distribution}, there is a Lagrangian infinitesimal foliation
		\begin{equation}
			\infcohcplfil{B}{k}\to \infcohcplfil{B}{k[T]}(=:\ccA)
		\end{equation}enhancing the distribution $\dL_{B/k}\to E^\vee[1]$. Its formal quotient $A:=|\ccA|$ satisfies that
		\begin{enumerate}[label=(\arabic*)]
			\item there are $r$-cycles $k\linsp{S}\to F_1\ccA\to A$ lifting a basis of $\gr_1\ccA\simeq E^\vee$, which identifies $B$ as the derived quotient $A\sslashop(S)\simeq A\otimes_{k[S]}k$;
			\item $A$ is an adic ring, and the pro-system $A_n:=F_{[0,n]}\ccA$ pro-corepresents $\spf(A)$. The cotangent complex $\dL_{\spf A/k}$ is given by $\Omega_{k[T]/k}\otimes_{k[T]} A$.
		\end{enumerate}
	\end{prop}
	\begin{proof}
			The left vertical arrow in (\ref{e:Lag dist}) defines a map $k[T]\to B$ in $\scr_k$, whose cotangent complex is the cofibre $\mathrm{cofib}(E\to \dL_{B/k})\simeq E^\vee[1]$. This map induces an integrable foliation-like algebra $\ccA:=\infcohcplfil{B}{k[T]}$. Note that $\ccA$ is complete and connective by definition. Hence every $-1$ cycle is null homotopic for $F_{\ge2}\ccA$, showing that $\ccA$ is isotropic with respect to $\omega\in \infsymp{-1}{k}{\spec B}$. The isotropy structure turns out to be Lagrangian following Construction \ref{c: Lagrangian distribution}. Thus we obtain the derived quotient $A\sslashop(S)\simeq B$ and $\spf A\simeq \colim\spec A_n$, as in Example \ref{rk: adic vs Infcoh}.

			We now compute $\dL_{\spf(A)/k}$. The shifted tangent $\dT_{B/k[T]}[1](\to \dT_{B/k}[1])$ is a perfect \textit{partition Lie algebroid} and $A$ can be identified with its Chevalley--Eilenberg complex \cite[Proposition 4.5, 4.23]{fu2024duality}. Since $A\to B$ is a derived surjection, there is $\dL_{B/k[T]}\simeq \dL_{B/A}$ by \cite[Theorem 2.25]{brantner2025formal}. In other words, there is $\dL_{A/k[T]}\otimes_A B\simeq 0$ and furthermore $\dL_{A/k[T]}\otimes_A A_n\simeq 0$, since $A_n$ can be rebuilt from finite steps of $k$-linear extensions by free $B$-modules. Hence, there is
			\[\Omega_{k[T]/k}\otimes_{k[T]}A_n\simeq \dL_{A/k}\otimes_{A}A_n.\]Passing to the limits, we obtain $\Omega_{k[T]/k}\otimes_{k[T]}A\simeq \dL_{\spf(A)/k}$.
	\end{proof}

	\section{Darboux charts}\label{sec: 4 Proof}
	
	This section is devoted to cooking up the formal integration results in \S\ref{sec: Formal integration of infinitesimal symplectic forms} into a minimal formal Darboux chart. Every such formal Darboux chart can be algebraized after a mild modification. Here we specify the inputs, and later we will only adjust these inputs for algebraization.
	
	\begin{inputdata}\label{input}
		Let $B$ be a finitely presented simplicial $k$-algebra, and fix $\omega\in \infsymp{\spec B}{k}{-1}$. In this section, we work with the following chosen data:
		\begin{enumerate}[label=(\arabic*)]
			\item A jet-Frobenius residue $\res_{JF}\in \pi_{-1}\infformtwo{\spec (B)}{k}{-1}$ such that $\omega-\res_{JF}$ admits a lift $\lambda\in\infsympex{\spec(B)}{k}{-1}$.
			\item An integrable Lagrangian foliation $\ccA:=\infcohcplfil{B}{k[T]}$, where $E:=(\gr_1\ccA)^\vee$ can be identified with the free $B$-module $B\linsp{T}$ and fits into the diagram \[\begin{tikzcd}[ampersand replacement=\&]
				\& {\dT_{B/k}[1]} \\
				E \& {\dL_{B/k}} \& {E^\vee[1].}
				\arrow["\simeq"{description}, from=1-2, to=2-2]
				\arrow[from=1-2, to=2-3]
				\arrow[from=2-1, to=1-2]
				\arrow[from=2-1, to=2-2]
				\arrow[from=2-2, to=2-3]
			\end{tikzcd}\]
			\item A fixed isotropy structure $f\in A(:=|\ccA|)$ of $\ccA$ with respect to $d_{\Inf}\lambda$.
		\end{enumerate}
			In particular, we have fixed a cycle $\widetilde{f}:=(\lambda,f)$ of $\widetilde{A}:=A\times_{A_1}\infoneexcohcpl{B}{k}$, where $A_1:=F_{[0,1]}\ccA$.
		
	\end{inputdata}

  These data always exist after passing to a Zariski neighborhood of some schematic point $x\in\spec(B)$, see Corollary \ref{cor: wJF infinite degree}, Construction \ref{c: Lagrangian distribution}, and Proposition \ref{prop: Lag Foliation}.
  
  Our approach to proving the main theorem is that, we will show the following square is a commutative and Cartesian
	% https://q.uiver.app/#q=WzAsNCxbMCwwLCJcXHNwZWMgQiJdLFsxLDAsIlxcc3BmIEEiXSxbMCwxLCJcXHNwZiBBIl0sWzEsMSwiVF4qX3tcXHNwZiBBfSJdLFsxLDMsIjAiXSxbMiwzLCJkZiIsMl0sWzAsMl0sWzAsMV1d
	\[\begin{tikzcd}[ampersand replacement=\&]
		{\spec B} \& {\spf A} \\
		{\spf A} \& {T^*_{\spf A}}
		\arrow[from=1-1, to=1-2]
		\arrow[from=1-1, to=2-1]
		\arrow["0", from=1-2, to=2-2]
		\arrow["df"', from=2-1, to=2-2]
	\end{tikzcd}\]where
	$T^*_{\spf A}:=\mathbb{V}_{\spf A}(\dT_{\spf A/k})$ is the cotangent space, which then exhibits $\spec B$ as the derived critical locus of $f:\spf A\to \mathbb{A}^{1}_k$. Then, we identify $\omega$ with the infinitesimal symplectic form of degree $-1$ determined by $\dcrit(f)$, see \cite[\S2.2]{pantev2013shifted} and Example \ref{exa: dcrit} for the adaptation to the mod $p$ case.
	
	Our proof is inspired by the works of Brav--Bussi--Joyce and Bouaziz--Grojnowski \cite{brav2019darboux, bouaziz2013d}. However, in positive characteristics, there is no access to models of derived foliations using graded mixed cdgas, which prevents us from writing down an explicit shifted potential. Instead, we carry out a detailed study of square-zero extensions.
	
	\subsection{Square-zero extensions}\label{sec: 4.1 sqz extensions}
	\begin{construction}[Continuing Construction \ref{c: non-connective cotangent complex}]\label{c: square-zero extension} The projection $\mathcal{M}\to \dalg$ has a natural section given by $B\mapsto \dL_B[-1]$. Let $\mathcal{M}_{\dL[-1]/}$ be the relative slice category, which is classified by the $\Catoo$-valued functor $B\mapsto \m_{B,\dL_{B}[-1]/}$. Therefore, the forgetful functor\begin{equation*}
			\mathcal{M}_{\dL[-1]/}\to \mathcal{M}
		\end{equation*}
		is a relative right adjoint over $\dalg$ and is fibrewise monadic. Meanwhile, let $\fil_{[0,1]}\dalg\subset\fil\dalg$ be the full subcategory of complete $A$ such that $\gr A\in \gr_{[0,1]}\dalg$, and consider the left adjoint $F\Pi$ of $\gr:\fil_{[0,1]}\dalg\to \gr_{[0,1]}\dalg$ (relative to $\dalg$ by $\gr_0$). The relative left adjoint $F\Pi$ sends $(B,M)$ to $\infoneexcohcpl{B}{\mathbb{Z}}\oplus M((1))$ as a fraction of Lemma \ref{lemma: graded pieces of infcoh}, whose graded pieces form $(B,\dL_{B}[-1]\oplus M)$. Therefore, $\gr\circ F\Pi$ defines monads, fibrewise over $\dalg$, that coincide with those determined by $\mathcal{M}_{\dL[-1]/}$. Therefore, we can identify $\mathcal{M}_{\dL_{[-1]}/}$ with $\fil_{[0,1]}\dalg$ via the assignment
		\begin{equation*}
			(R,\rho:\dL_R[-1]\to M)\mapsto\big( M[-1]\to R\oplus_{\rho }M\big).
		\end{equation*}
		We call $\mathcal{M}_{\dL[-1]/}\simeq \fil_{[0,1]}\dalg$ the \infcat\ of \textit{square-zero extensions of derived rings}.
		
		Given a base ring $R\in \scr$, the theory of square-zero extensions over $R$ is encoded by the slice category $\fil_{[0,1]}\dalg_R\simeq (\mathcal{M}_{\dL[-1]/})_{(R,0)/}$. In this context, a square-zero extension is equivalent to a pair $(B,\rho:\dL_{B/R}[-1]\to M)$ or as an $R$-linear differential $d:B\to M$.
		
	\end{construction}
	\begin{exa}\label{exa: squarezero}
		Here are the tricks that help to obtain or detect square-zero extensions
		\begin{enumerate}[label=(\arabic*)]
			\item	A filtered derived algebra is an abundant supply of square-zero extensions. For every natural numbers $n\ge 0$, there is a functor
			\begin{equation*}
				\fil\dalg_k\to \fil_{[0,1]}\dalg_k
			\end{equation*}sending a filtered algebra $\mathcal{A}=(\ldots F_2\mathcal{A}\to F_1\mathcal{A}\to F_0\mathcal{A})$ to $F_{[n-1,n]}\mathcal{A}\to F_{[0,n]}\mathcal{A}$, in other words the square-zero extension of $F_{[0,n-1]}\mathcal{A}\to \gr_n \mathcal{A}[1]$.
			\item Given an arrow $S\to R$ in $\dalg_k$, and an $R$-linear map $\rho:\dL_{R/k}\to N[1]$, the following top-left corner defined by pullback in $\dalg_k$ is a square-zero extension by $S\to R\xrightarrow{\rho d_{\deRham}}N[1]$
			% https://q.uiver.app/#q=WzAsNixbMCwwLCJTXFxvcGx1c197XFxyaG99TiJdLFswLDEsIlJcXG9wbHVzX3tcXHJob31OIl0sWzEsMCwiUyJdLFsxLDEsIlIiXSxbMiwxLCJOWzFdIl0sWzIsMCwiTlsxXSJdLFszLDRdLFswLDFdLFsxLDNdLFsyLDNdLFswLDJdLFsyLDVdLFs1LDQsIiIsMSx7ImxldmVsIjoyLCJzdHlsZSI6eyJoZWFkIjp7Im5hbWUiOiJub25lIn19fV0sWzAsMywiIiwxLHsic3R5bGUiOnsibmFtZSI6ImNvcm5lciJ9fV1d
			\[\begin{tikzcd}[ampersand replacement=\&]
				{S\oplus_{\rho}N} \& S \& {N[1]} \\
				{R\oplus_{\rho}N} \& R \& {N[1].}
				\arrow[from=1-1, to=1-2]
				\arrow[from=1-1, to=2-1]
				\arrow["\lrcorner"{anchor=center, pos=0.125}, draw=none, from=1-1, to=2-2]
				\arrow[from=1-2, to=1-3]
				\arrow[from=1-2, to=2-2]
				\arrow[equals, from=1-3, to=2-3]
				\arrow[from=2-1, to=2-2]
				\arrow[from=2-2, to=2-3]
			\end{tikzcd}\]
			\item For every $R\in\dalg_k$ and $R$-linear maps $\dL_{R/k}\xrightarrow{\rho}N[1]\xrightarrow{l}M[1]$, then $R\oplus_{\rho}N$ is exhibited as the square-zero extension of $R\oplus_{l\rho}M$ by $\mathrm{fib}(l)$.
		\end{enumerate}
	\end{exa}

	\begin{exa}[Following \S\ref{s:3.2 Lagrangian foliation}]\label{exa: An}
		The algebras $A_n$, $F_{[0,n]}\mathbbb{\Pi}\mathstrut_{B/k}$ are square-zero extensions of $A_{n-1}$ and $F_{[0,n-1]}\mathbbb{\Pi}\mathstrut_{B/k}$ by $\lsym^{n}_B E^\vee$ and $\lsym^{n}_B (\dL_{B/k}[-1])$, respectively. The $B$-linear map $\dL_{B/k}\to E^{\vee}[1]$ exhibits $\infoneexcohcpl{B}{k}$ as a square-zero extension of $A_1$ by $E[-1]$. Additionally, $\widetilde{A}_n:=A_n\times_{A_1}\infoneexcohcpl{B}{k}$ can be determined as the extension of $A_n$ by $E$ for $n\in[1,\infty]$. Observe that $E\simeq B\otimes_A\dL_{\spf A/k}$. We are going to to see that the differential
		\[A_n\xrightarrow{d} E\]defining $\widetilde{A}_n$ is actually a ``quotient'' of the universal differential $d_{\deRham}:A\to \dL_{\spf A/k}$.
	\end{exa}
	We now approximate $\infcohcplfil{B}{k}$ by the internal square-zero extensions of filtered algebras within the context of Input \ref{input}.
		\begin{variant}[Following Variant \ref{var: internal infinitesimal cohomology}]  Write $\mathcal{Q}:=\infcohcplfil{B}{k}$ for short. There is a natural map of fully complete bi-filtered algebras
		\begin{equation*}
			\mathcal{Q}\to \widehat{F\mathbbb{\Pi}}\mathstrut^{int}_{\ccA/\mathcal{Q}}.
		\end{equation*}
		The internal filtration degree $[0,1]$ part of $\widehat{F\mathbbb{\Pi}}\mathstrut^{int}_{\ccA/\mathcal{Q}}$ is then a square-zero extension determined by
		\begin{equation}\label{e: internal diff}
			\mathbbb{d}:\ccA\to \dL^{int}_{\ccA/\mathcal{Q}}.
		\end{equation}Moreover, there is a morphism of complete filtered algebras $\mathcal{Q}\to \mathrm{fib}(\mathbbb{d})$.
	\end{variant}
	Next, we show that $\mathbbb{d}$ can be identified with an adic filtration on $d_{\deRham}:A\to \dL_{\spf A/k}$, built in an \textit{ad hoc} way. Recall that there exists some map $k[S]\to A$ such that $A\otimes_{k[S]}k\simeq B$, see Proposition \ref{prop: Lag Foliation} (1). Let $I=(S)\subset k[S]$ be the ideal of the origin. There is an ordinary differential induced by the quotient of the de Rham differential of $k[S]$,
	\[k[S]/I^{n+1}\xrightarrow{d}\Omega_{k[S]/k}\otimes_{k[S]}k[S]/I^n.\]For all $n\ge 1$, this quotient induces the following differential $d_n$
	\begin{align*}
		d_n:A_n\simeq A\otimes_{k[S]}k[S]/I^{n+1}\xrightarrow{d_{\deRham}}&\dL_{\spf A/k}\otimes_{A}A_n\bigoplus_{\Omega_{k[S]/k}\otimes A_n}\dL_{k[S]/I^{n+1}/k}\otimes_{k[S]/I^{n+1}}A_n\\
		\to& \dL_{\spf A/k}\otimes_{A}A_{n-1} \bigoplus_{\Omega_{k[S]/k}\otimes A_{n-1}}\Omega_{k[S]/k}\otimes A_{n-1}\\
		\simeq& \dL_{\spf A/k}\otimes_{A}A_{n-1},
	\end{align*}
	where the first arrow is the universal differential of $A_n$. These differentials $d_n$ ($n\ge1$) form an inverse system and give rise to an internal differential of $\ccA$
	\begin{equation}\label{e: adic diff}
		d_{\adic}:\ccA\to \dL_{\spf A/k}\otimes_A \ccA(1).
	\end{equation}
	\begin{prop}\label{prop: two diff one geometry}
		The internal differential $\mathbbb{d}$ in (\ref{e: internal diff}) is homotopic to $d_{\adic}:\ccA\to \dL_{\spf A/k}\otimes_A \ccA(1)$ in Line (\ref{e: adic diff}), as objects in $\fil_{[0,1]}\fil\dalg$, the \infcat\ of (truncated) bi-filtered algebras.
	\end{prop}
	\begin{proof}
	The $\gr_0$-piece of $\mathrm{fib}(d_{\adic})$ is equivalent to $B$, so there is an arrow $\mathcal{Q}\to \mathrm{fib}(d_{\adic})$ by the universal property of $\mathcal{Q}\simeq \infcohcplfil{B}{k}$. Thus, the ad hoc differential $d_{\adic}$ is $\mathcal{Q}$-linear, and there is a natural factorization
	\[d_{\adic}:\ccA\xrightarrow{\mathbbb{d}}\dL^{int}_{\ccA/\mathcal{Q}}\xrightarrow{f}\dL_{\spf A/k}\otimes_{A}\ccA(1).\]It is then sufficient to show $f$ is an equivalence of $\ccA$-modules.
	Recall that $\ccA:=\infcohcplfil{B}{k[T]}$ is integrable by construction, so $\dL^{int}_{\ccA/\mathcal{Q}}\simeq\dL_{\spf A/k}\otimes_A \ccA(1)$ by Lemma \ref{lemma: integrable crystal}. At the same time, there is $\dL_{\spf A/k}\otimes_A \ccA(1)\simeq\Omega_{k[T]}\otimes_{k[T]}\ccA(1)$ by Proposition \ref{prop: Lag Foliation}, which means that both ends of $f$ are finite free filtered modules (up to a weight shift) over $\ccA$. By Nakayama's lemma, the question is reduced to prove that the base change $f\otimes B:E(1)\to E(1)$ is an equivalence.

	Now, note that the external graded pieces of $\mathcal{Q}\to \widehat{F\mathbbb{\Pi}}\mathstrut^{int}_{\ccA/\mathcal{Q}}$ agree with $\gr\mathcal{Q}\to\widehat{F\mathbbb{\Pi}}\mathstrut^{int}_{\gr\ccA/\gr\mathcal{Q}}$, where the target is the internal infinitesimal cohomology of graded algebras. Since $\gr\mathcal{Q}\to \gr\ccA$ is given by applying $\lsym_B$ to $\dL_{B/k}[-1]\to E^\vee$, the map $\gr\mathcal{Q}\to \gr^{int}\mathbbb{\Pi}^{int}_{\gr\ccA/\gr\mathcal{Q}}$ is given by
	\[\lsym_B\big(\dL_{B/k}[-1](1)^{ext}\big)\to \lsym_B\big(E^\vee(1)^{ext}\oplus E[-1](1)^{ext}(1)^{int}\big).\]
	In particular, the map
	$\gr_1\mathcal{Q}\to \gr_1F^{int}_{[0,1]}\mathbbb{\Pi}^{int}_{\ccA/\mathcal{Q}}$ can be identified with the identity of $\dL_{B/k}[-1]$. Then, we have the equivalence of filtered algebras $F_{[0,1]}\mathcal{Q}\simeq  F_{[0,1]}F^{int}_{[0,1]}\mathbbb{\Pi}^{int}_{\ccA/\mathcal{Q}}$, where we only remember the external filtration of the target. This equivalence exhibits $F_{[0,1]}\mathbbb{\Pi}_{B/k}$ as the internal square-zero extension of $F_{[0,1]}\mathbbb{d}$. Meanwhile, $F_{[0,1]}\mathbbb{\Pi}_{B/k}$ is also the square-zero extension of $d_{\deRham}:B\to \dL_{B/k}$. Therefore, there is a commutative map
	% https://q.uiver.app/#q=WzAsNixbMCwwLCJcXGNjQSJdLFswLDEsIlxcZExee2ludH1fe1xcY2NBL1xcbWF0aGNhbHtRfX0iXSxbMSwwLCJGX3tbMCwxXX1cXGNjQSJdLFsxLDEsIkUoMSkiXSxbMiwwLCJCIl0sWzIsMSwiXFxkTF97Qi9rfSgxKSwiXSxbMCwyXSxbMiwzLCJGX3tbMCwxXX1cXG1hdGhiYmJ7ZH0iXSxbMCwxLCJcXG1hdGhiYmJ7ZH0iLDJdLFsxLDNdLFsyLDRdLFszLDVdLFs0LDUsImRfe1xcZGVSaGFtfSJdXQ==
	\[\begin{tikzcd}[ampersand replacement=\&]
		{\mathllap{(k[T]\to)\quad}\ccA} \& {F_{[0,1]}\ccA} \& B \\
		{\dL^{int}_{\ccA/\mathcal{Q}}} \& {E(1)} \& {\dL_{B/k}(1),}
		\arrow[from=1-1, to=1-2]
		\arrow["{\mathbbb{d}}"', from=1-1, to=2-1]
		\arrow[from=1-2, to=1-3]
		\arrow["{F_{[0,1]}\mathbbb{d}}", from=1-2, to=2-2]
		\arrow["{d_{\deRham}}", from=1-3, to=2-3]
		\arrow[from=2-1, to=2-2]
		\arrow[from=2-2, to=2-3]
	\end{tikzcd}\]where $E\to \dL_{B/k}$ is as in Construction \ref{c: Lagrangian distribution}. Hence the elements of $T$ are sent to a basis of $E$ by this diagram. Meanwhile, $T$ is sent to a basis of $E$ by $A\to A_1\xrightarrow{d_1}E$. Thus, $|f|\otimes B$ and then $|f|$ are equivalences. However, both ends of $f$ are purely in weight $1$, so $f$ is an equivalence as well.
	\end{proof}
	We now add further examples of square-zero extensions after Example \ref{exa: An}.
	\begin{prop}\label{prop: ccB}
		Let $\ccC$ denote the square-zero extension by $\mathbbb{d}$ in (\ref{e: internal diff}) and set $C_n:=F_{[0,n]}\ccC$. Then, $C_2$ fits into two sequences of square-zero extensions:
		\begin{enumerate}[label=(\arabic*)]
			\item The maps $C_2\to \widetilde{A}_2\to A_2$ exhibit the first two terms as square-zero extensions of $A_2$ by the $A_2$-modules $A_1\otimes_A \dL_{\spf_A/k}\to B\otimes_A \dL_{\spf A/k}$ respectively.
			\item If $\operatorname{char}(k)>2$, there are square-zero extensions of $\infoneexcohcpl{B}{k}$
			\[\inftwoexcohcpl{B}{k}\to C_2\to \widetilde{A}_2\] by  $\lsym^2_B(\dL_{B/k}[-1])\to M \to \lsym^2_BE^\vee$. Denote $H:E^\vee\to E$ the connection map of the Lagrangian distribution in Input \ref{input}. Here $M$ is the fibre of \[\lsym^2_B E^\vee\xrightarrow{a\cdot b\mapsto H(a)\otimes b+b\otimes H(a)}E^\vee\otimes_B E.\]
			
		\end{enumerate}
	\end{prop}
	\begin{proof}
		By passing to the graded pieces, we can apply Proposition \ref{prop: two diff one geometry}.
	\end{proof}
	\begin{alert}
		It is crucial to assume $p>2$ in (2). Otherwise, $\lsym^2_{B}(\dL_{B/k}[-1])$ would contain terms arising from Dyer--Lashof operations.
	\end{alert}
	\subsection{Minimal formal Darboux charts}\label{sec: 4.2 Minimal formal Darboux}
	Recall that $\spf A\simeq \colim_n \spec A_n$. The cotangent space $T^*_{\spf A}$ (relative to $k$) is the relative linear stack of the tangent complex $\dT_{\spf A/k}$. More explicitly, we have\begin{equation}\label{e: square of df Dcrit}
		T^*_{\spf A}\simeq \mathbb{V}_{\spec A_n}(\dT_{\spf A/k}\otimes_AA_n)\simeq \colim_n \spec \lsym_{A_n}(A_n\otimes_{k[T]}T_{k[T]/k}).
	\end{equation}
	
	\begin{prop}\label{prop: dCrit stk}
		In the context of Input \ref{input}, the cycle $\widetilde{f}=(\lambda,f)\in \pi_0\widetilde{A}$ induces a cartesian square in $\dSt_k$,
		\begin{equation}\label{e: prop: dCrit stk}
		% https://q.uiver.app/#q=WzAsNCxbMCwwLCJcXHNwZWMgQiJdLFsxLDAsIlxcc3BmIEEiXSxbMCwxLCJcXHNwZiBBIl0sWzEsMSwiVF4qX3tcXHNwZiBBfSJdLFsxLDMsIjAiXSxbMiwzLCJkX3tcXGRlUmhhbX1mIiwyXSxbMCwyXSxbMCwxXSxbMCwzLCIiLDAseyJzdHlsZSI6eyJuYW1lIjoiY29ybmVyIn19XV0=
		\begin{tikzcd}[ampersand replacement=\&]
			{\spec B} \& {\spf A} \\
			{\spf A} \& {T^*_{\spf A}}
			\arrow[from=1-1, to=1-2]
			\arrow[from=1-1, to=2-1]
			\arrow["\lrcorner"{anchor=center, pos=0.125}, draw=none, from=1-1, to=2-2]
			\arrow["0", from=1-2, to=2-2]
			\arrow["{d_{\deRham}f}"', from=2-1, to=2-2]
		\end{tikzcd}
		\end{equation}
	\end{prop}
	\begin{proof}
	Unveiling the definition, an $R$-point of $T^*_{\spf A}$ can be represented by $A\to A_n\to R$ for some $n$ together with an $R$-linear map $\dT_{\spf A/k}\otimes_A R\to R$. The commutativity amounts to showing that $d_{\deRham}f:\dT_{\spf A/k}\to A$ factors through $F_{\ge1}\ccA$, or dually, that the cycle $d_{\deRham}f$ is sent to $0$ by $\dL_{\spf A/k}\to E$. This holds because we have $(\lambda,f)\in \pi_0\widetilde{A}$ lifting $f$ along the fibre sequence
	\begin{equation*}
		\widetilde{A}\to A\to E.
	\end{equation*}
	For later use, let $\widetilde{d_{\deRham}f}$ be the lift in $\mathrm{fib}(\dT_{\spf A}\to E)$. 
	
	Next, we prove that the square (\ref{e: square of df Dcrit}) is cartesian. By Proposition \ref{prop: ccB} (2), there is a span of square-zero extensions of derived algebras,% https://q.uiver.app/#q=WzAsOSxbMCwwLCJcXGluZnR3b2V4Y29oY3Bse0J9e2t9Il0sWzAsMSwiXFxpbmZ0d29leGNvaGNwbHtCfXtrfSJdLFsxLDAsIlxcaW5mb25lZXhjb2hjcGx7Qn17a30iXSxbMSwxLCJcXHdpZGV0aWxkZXtBfV8yIl0sWzEsMiwiXFx3aWRldGlsZGV7QX1fMiJdLFswLDIsIkNfMiJdLFsyLDIsIkVcXG90aW1lc19CRV5cXHZlZVxcbWF0aHJsYXB7XFxxdWFkKFxcc2ltZXEgXFxncl8yXFxkTF57aW50fV97XFxtYXRoY2Fse1F9L1xcY2NBfSl9Il0sWzIsMSwiXFxtYXRocm17ZmlifShFXFxvdGltZXNfQiBFXlxcdmVlXFx0byBcXHdlZGdlXjIgRV5cXHZlZSkiXSxbMiwwLCJcXGxzeW1fQihcXGRMX3tCL2t9Wy0xXSlbMV0iXSxbMCwxLCIiLDAseyJsZXZlbCI6Miwic3R5bGUiOnsiaGVhZCI6eyJuYW1lIjoibm9uZSJ9fX1dLFszLDQsIiIsMCx7ImxldmVsIjoyLCJzdHlsZSI6eyJoZWFkIjp7Im5hbWUiOiJub25lIn19fV0sWzAsMl0sWzMsMl0sWzEsM10sWzEsNV0sWzUsNF0sWzQsNl0sWzcsNl0sWzMsN10sWzIsOCwiZF97XFxJbmZ9Il0sWzcsOF1d
	\[\begin{tikzcd}[ampersand replacement=\&]
		{\inftwoexcohcpl{B}{k}} \& {\infoneexcohcpl{B}{k}} \& {\lsym_B(\dL_{B/k}[-1])[1]} \\
		{\inftwoexcohcpl{B}{k}} \& {\widetilde{A}_2} \& {\mathrm{fib}(E\otimes_B E^\vee\to \wedge^2 E^\vee)} \\
		{C_2} \& {\widetilde{A}_2} \& {E\otimes_BE^\vee\mathrlap{\quad(\simeq \gr_2\dL^{int}_{\mathcal{Q}/\ccA})}}
		\arrow[from=1-1, to=1-2]
		\arrow[equals, from=1-1, to=2-1]
		\arrow["{d_{\Inf}}", from=1-2, to=1-3]
		\arrow[from=2-1, to=2-2]
		\arrow[from=2-1, to=3-1]
		\arrow[from=2-2, to=1-2]
		\arrow[from=2-2, to=2-3]
		\arrow[equals, from=2-2, to=3-2]
		\arrow[from=2-3, to=1-3]
		\arrow[from=2-3, to=3-3]
		\arrow[from=3-1, to=3-2]
		\arrow[from=3-2, to=3-3]
	\end{tikzcd}\]
	The chosen cycle $(\lambda,f)$ can be sent to $\widetilde{A}_2$, then lands in the three right vertices in this diagram. Its image in the top-right corner is the underlying non-degenerate $2$-form of $\omega$, while its image in the middle term is given by the Lagrangian isotropy determined by $f$, so we see that the image of the cycle $(\lambda, f)$ in $E\otimes_B E^\vee$ induces an equivalence $E^\vee\simeq E^\vee$, see Remark \ref{rk: Lagrangian isotropy as a non-degenerated matrix}. Meanwhile, $C_2$ and $\widetilde{A}_2$ are also square-zero extensions of $A_2$, and the lowest line can be extended to 
	% https://q.uiver.app/#q=WzAsNixbMCwwLCJDXzIiXSxbMCwxLCJcXHdpZGV0aWxkZXtBfV8yIl0sWzEsMCwiQV8yIl0sWzEsMSwiQV8yIl0sWzIsMCwiQV8xXFxvdGltZXNfQlxcZExfe1xcc3BmIEEva30iXSxbMiwxLCJFIl0sWzAsMV0sWzEsM10sWzIsMywiIiwyLHsibGV2ZWwiOjIsInN0eWxlIjp7ImhlYWQiOnsibmFtZSI6Im5vbmUifX19XSxbMCwyXSxbMiw0LCJkXzIiLDJdLFszLDVdLFs0LDVdXQ==
	\[\begin{tikzcd}[ampersand replacement=\&]
		{C_2} \& {A_2} \& {A_1\otimes_B\dL_{\spf A/k}} \\
		{\widetilde{A}_2} \& {A_2} \& E
		\arrow[from=1-1, to=1-2]
		\arrow[from=1-1, to=2-1]
		\arrow["{d_2}"', from=1-2, to=1-3]
		\arrow[equals, from=1-2, to=2-2]
		\arrow[from=1-3, to=2-3]
		\arrow[from=2-1, to=2-2]
		\arrow[from=2-2, to=2-3]
	\end{tikzcd}\]by Proposition \ref{prop: ccB} (1). It implies that the image of $(\lambda, f)$ into $E\otimes_B E^\vee$ is homotopic to the image of $\widetilde{d_{\deRham}f}$ under $F_1\ccA\otimes_A\dL_{\spf A/k}\to \gr_1\otimes_A\dL_{\spf A/k}\simeq E^\vee\otimes_B E$. Subsequently, in the square induced by the lift $\widetilde{d_{\deRham}f}$ % https://q.uiver.app/#q=WzAsNCxbMCwwLCJcXGRUX3tcXHNwZiBBL2t9Il0sWzEsMCwiRl8xXFxjY0EiXSxbMCwxLCJcXG1hdGhsbGFwe0VcXHNpbWVxfUJcXG90aW1lc19BXFxkVF97XFxzcGYgQS9rfSJdLFsxLDEsIlxcZ3JfMVxcY2NBXFxtYXRocmxhcHtcXHNpbWVxIEUsfSJdLFswLDJdLFswLDEsIlxcd2lkZXRpbGRle2Rfe1xcZGVSaGFtfWZ9Il0sWzEsM10sWzIsMywiXFxzaW1lcSJdXQ==
	\begin{equation}\label{e: T --> F_1A}
		\begin{tikzcd}[ampersand replacement=\&]
			{\dT_{\spf A/k}} \& {F_1\ccA} \\
			{\mathllap{E\simeq}B\otimes_A\dT_{\spf A/k}} \& {\gr_1\ccA\mathrlap{\simeq E,}}
			\arrow["{\widetilde{d_{\deRham}f}}", from=1-1, to=1-2]
			\arrow[from=1-1, to=2-1]
			\arrow[from=1-2, to=2-2]
			\arrow["\simeq", from=2-1, to=2-2]
		\end{tikzcd}
	\end{equation}the bottom horizontal arrow is an equivalence of $B$-modules. Then consider % https://q.uiver.app/#q=WzAsNSxbMCwwLCJcXGNjQVtcXGRUX3tcXHNwZiBBL2t9XFxvdGltZXNfQVxcY2NBKDEpXSJdLFsxLDAsIlxcY2NBIl0sWzAsMSwiXFxjY0EiXSxbMSwxLCJcXG1hdGhjYWx7RH0iXSxbMiwxLCJCIl0sWzAsMSwiMCJdLFswLDIsIlxcd2lkZXRpbGRle2Rfe1xcZGVSaGFtfX0iLDJdLFsyLDNdLFsxLDNdLFszLDRdXQ==
	\[\begin{tikzcd}[ampersand replacement=\&]
	{\ccA[\dT_{\spf A/k}\otimes_A\ccA(1)]} \& \ccA \\
	\ccA \& {\mathcal{D}} \& B,
	\arrow["0", from=1-1, to=1-2]
	\arrow["{\widetilde{d_{\deRham}}}"', from=1-1, to=2-1]
	\arrow[from=1-2, to=2-2]
	\arrow[from=2-1, to=2-2]
	\arrow[from=2-2, to=2-3]
	\end{tikzcd}\]where the square is a pushout of filtered algebras. Since $\dT_{\spf A/k}$ is finitely free over $A$ (Proposition \ref{prop: Lag Foliation}), $\mathcal{D}$ is computed by the Koszul complex and hence is complete. At the same time, the bottom equivalence in (\ref{e: T --> F_1A}) implies that $\gr\mathcal{D}\simeq B$. Finally, we have $A{\otimes_{0,A[\dT_{\spf A/k}],d_{\deRham}f}} A\simeq B$ so that (\ref{e: prop: dCrit stk}) is cartesian.
	\end{proof}
	
	Although $\spf A$ is not a derived Deligne-Mumford stack, it is a formally smooth affine formal derived scheme over $k$, so there is also a Liouville $1$-form (without shifting) and its differential $d_{\Inf}\Liouv$ is symplectic by applying exactly the same reasoning as in Proposition \ref{prop: Liouv}.
	
	\begin{prop}\label{prop: loop}
		Let $\mathrm{Liouv}$ be the tautological $1$-form on $T^*_{\spf A}$, and let $\lambda_f$ be the loop
		\begin{equation}\label{e: loop of exact structure}
			0=\iota^*0^*\mathrm{Liouv}\sim \iota^*(df)^*\mathrm{Liouv}\sim 0,
		\end{equation}where the middle homotopy is induced by the \textnormal{commutativity} in Proposition \ref{prop: dCrit stk} and the other homotopies are given by the \textnormal{isotropic structures}. {\normalfont Then}, $\lambda_f$ is equivalent to $\lambda$, the chosen exact symplectic form of degree $-1$ in Input \ref{input}.
	\end{prop}
	\begin{proof}%second proof
		There is a commutative diagram by Proposition \ref{prop: ccB} and \ref{prop: dCrit stk},
		% https://q.uiver.app/#q=WzAsMTAsWzEsMiwiXFxpbmZvbmVleGNvaGNwbHtUXipfe1xcc3BmIEF9fXtrfSJdLFsyLDMsIlxcaW5mb25lZXhjb2hjcGx7XFxzcGYgQX17a30iXSxbMiwxLCJcXGluZm9uZWV4Y29oY3Bse1xcc3BmIEF9e2t9Il0sWzIsMCwiXFxkTF97XFxzcGYgQS9rfVstMV0iXSxbMCwyLCJrWy0xXSJdLFsxLDAsIjAiXSxbMiw0LCJcXGRMX3tcXHNwZiBBL2t9Wy0xXSJdLFsxLDQsIkFbLTFdIl0sWzMsMiwiXFxpbmZvbmVleGNvaGNwbHtCfXtrfSJdLFs0LDIsIkVbLTFdLiJdLFswLDEsImRmXioiLDJdLFswLDIsIjBeKiJdLFszLDJdLFs0LDVdLFs1LDNdLFs0LDAsIlxcTGlvdXYiLDJdLFs2LDFdLFs3LDYsImRfe1xcZGVSaGFtfSJdLFs0LDcsImYiXSxbMiw4XSxbMSw4XSxbMyw5XSxbNiw5XSxbOSw4XV0=
		\[\begin{tikzcd}[ampersand replacement=\&]
			\& 0 \& {\dL_{\spf A/k}[-1]} \\
			\&\& {\infoneexcohcpl{\spf A}{k}} \\
			{k[-1]} \& {\infoneexcohcpl{T^*_{\spf A}}{k}} \&\& {\infoneexcohcpl{B}{k}} \& {E[-1].} \\
			\&\& {\infoneexcohcpl{\spf A}{k}} \\
			\& {A[-1]} \& {\dL_{\spf A/k}[-1]}
			\arrow[from=1-2, to=1-3]
			\arrow[from=1-3, to=2-3]
			\arrow[from=1-3, to=3-5]
			\arrow[from=2-3, to=3-4]
			\arrow[from=3-1, to=1-2]
			\arrow["\Liouv"', from=3-1, to=3-2]
			\arrow["f", from=3-1, to=5-2]
			\arrow["{0^*}", from=3-2, to=2-3]
			\arrow["{df^*}"', from=3-2, to=4-3]
			\arrow[from=3-5, to=3-4]
			\arrow[from=4-3, to=3-4]
			\arrow["{d_{\deRham}}", from=5-2, to=5-3]
			\arrow[from=5-3, to=3-5]
			\arrow[from=5-3, to=4-3]
		\end{tikzcd}\]This diagram exhibits a $2$-cell (between two null-homotopic paths from $k[-1]$ to $F_{[0,1]}\mathbbb{\Pi}_{B/k}$) which encodes $\lambda_f$. Meanwhile, zooming into $k[-1]\to E[-1]$, there is a commutative diagram
		% https://q.uiver.app/#q=WzAsNixbMiwyLCJFWy0xXSJdLFsyLDEsIjAiXSxbMSwyLCJBXzFbLTFdIl0sWzEsMSwiXFxpbmZvbmVleGNvaGNwbHtCfXtrfVstMV0iXSxbMCwyLCJBWy0xXSJdLFswLDAsImtbLTFdIl0sWzIsMCwiZF8xIl0sWzMsMl0sWzMsMV0sWzEsMF0sWzMsMCwiIiwxLHsic3R5bGUiOnsibmFtZSI6ImNvcm5lciJ9fV0sWzQsMl0sWzUsNCwiZiJdLFs1LDFdLFs1LDMsIiIsMix7InN0eWxlIjp7ImJvZHkiOnsibmFtZSI6ImRhc2hlZCJ9fX1dXQ==
		\[\begin{tikzcd}[ampersand replacement=\&]
			{k[-1]} \\
			\& {\infoneexcohcpl{B}{k}[-1]} \& 0 \\
			{A[-1]} \& {A_1[-1]} \& {E[-1],}
			\arrow[dashed, from=1-1, to=2-2]
			\arrow[from=1-1, to=2-3]
			\arrow["f", from=1-1, to=3-1]
			\arrow[from=2-2, to=2-3]
			\arrow[from=2-2, to=3-2]
			\arrow["\lrcorner"{anchor=center, pos=0.125}, draw=none, from=2-2, to=3-3]
			\arrow[from=2-3, to=3-3]
			\arrow[from=3-1, to=3-2]
			\arrow["{d_1}", from=3-2, to=3-3]
		\end{tikzcd}\]where the dotted filler is $\lambda$. It means that $\lambda_{f}\simeq \lambda$.\end{proof}
	
	\subsection{Proof of the main theorem}\label{s5.3: Algebraization of Darboux charts}
	This subsection aims to prove Theorem \ref{main theorem} and to discuss how the choice of $\res_{JF}$ influences the Darboux charts.
	
	\begin{construction}[Refining Input \ref{input}]\label{c: Refining}Fix a schematic point $x\in \spec B$. Note that $|\spec B|=|\spec A_2|=|\spf A|$ so that $x$ can be regarded as a point of $\spec A_2$. After passing to a Zariski neighborhood (still denoted as $\spec A_2)$, there is a factorization by \cite[18.4.7]{EGAIV4}
			\[\spec\pi_0A_2\xrightarrow{\text{close immersion}}\spec R\xrightarrow{\text{\'etale}} \spec k[T].\]By shrinking $\spec B$ and $\spf A$ correspondently, we obtain a natural commutative diagram
			
			% https://q.uiver.app/#q=WzAsNCxbMCwwLCJrW1RdIl0sWzAsMSwiUiJdLFsxLDAsIkEiXSxbMSwxLCJcXHBpXzBBXzIiXSxbMCwyXSxbMiwzXSxbMCwxXSxbMSwzLCIiLDAseyJzdHlsZSI6eyJoZWFkIjp7Im5hbWUiOiJlcGkifX19XSxbMSwyXV0=
			\[\begin{tikzcd}[ampersand replacement=\&]
				{k[T]} \& A \\
				R \& {\pi_0A_2,}
				\arrow[from=1-1, to=1-2]
				\arrow[from=1-1, to=2-1]
				\arrow[from=1-2, to=2-2]
				\arrow[from=2-1, to=1-2]
				\arrow[two heads, from=2-1, to=2-2]
			\end{tikzcd}\]where the lift $R\to A$ exists uniquely since $k[T]\to R$ is \'etale. Every isotropic structure $f$ of $\ccA$ can be modified modulo $\pi_0F_3\ccA$ so that it can be lifted to $f_R \in R$. There is then a natural map $\spec B\to \dcrit(f_R)$ given by 
			\[	R\sslashop (\partial_{t}f_R|t\in T)\to B.\]
			
			Let $J_R\subset R$ be the ideal generated by $\partial_{t}f_R$ for all $t\in T$. The radical $\sqrt{J_R}=\cap_i^n\mathfrak{p}_i$ is a finite intersection of the minimal primes containing $J_R$. By prime avoidance, we may choose a $g\in R$ such that $g\in \mathfrak{p}_i$ precisely when $\spec R/\mathfrak{p}_i$ does not contain any irreducible component of $\spec( \pi_0 B)^{red}$. We now replace $R$ by $R[g^{-1}]$ (still denoted as $R$ later).

	\end{construction}
	\begin{prop}\label{prop: filtered R vs ccA}
		Let $\what{R}_{J_R}$ be the ordinary $J_R$-adic completion. There is a natural equivalence\[\widehat{R}_{J_R}\simeq A,\]hence $A$ is actually discrete. Furthermore, the map $R\sslashop(\partial_t R)\to B$ is also an equivalence.
	\end{prop}
	\begin{proof}
		Set $I_R:=\ker(R\to \pi_0 B)$. Every irreducible component of $\spec \pi_0B$ is given by a $\mathfrak{q}\supset I_R$ in $R$. Then there is a faithfully flat ring morphism $R_{\mathfrak{q}}\to (R_{\mathfrak{q}})^\wedge_{I_R}$. However, the map $R\to (R_{\mathfrak{q}})^\wedge_{I_R}$ factors through $A_{\mathfrak{q}}\to  (R_{\mathfrak{q}})^\wedge_{I_R}$, which implies that $I_R(R_{\mathfrak{q}})^\wedge_{I_R}=J_R(R_{\mathfrak{q}})^\wedge_{I_R}$. Thus, we have $I_RR_{\mathfrak{q}}=J_RR_{\mathfrak{q}}$. This implies that $\sqrt{I_R}=\sqrt{J_R}$ and then $\what{R}_{J_R}=\what{R}_{I_R}$. The map $\spf(\what{R}_{J_R})\hookrightarrow \spec R$ is formally \'etale. Hence there is $\dL_{B/R}\simeq \dL_{B/\what{R}_{J_R}}$, which gives rise to an equivalence of foliation-like algebras
		\[\ccA\simeq \infcohcplfil{B}{R}\xrightarrow{\simeq}\infcohcplfil{B}{\what{R}_{J_R}}.\]We now obtain $\what{R}_{J_R}\simeq |\infcohcplfil{B}{\what{R}_{J_R}}|\simeq A$ by Proposition \ref{prop: inf foliation basics} (1). This implies that
		\[R\sslashop(\partial f_R)\simeq \what{R}_{J_R}\sslashop(\partial f_R)\simeq \gr_0\infcohcplfil{B}{\what{R}_{J_R}}\simeq B.\]\end{proof}
	
		{\theoremstyle{plain}
			\newtheorem*{restate}{Theorem \ref{main theorem}}
			\begin{restate}
				Let $k$ be a field of characteristic $p>2$, and let $X$ be a derived scheme over $k$ locally of finite presentation. Fix $
					\omega\in \infsymp{-1}{k}{X}$. For each schematic point $x\in |X|$, there exist data $(i,j,f,\phi,\res_{JF},\eta)$ called Darboux chart:
				% https://q.uiver.app/#q=WzAsNixbMCwwLCJYIl0sWzEsMCwiVSJdLFsyLDAsIlxcbWF0aGNhbHtVfSJdLFsxLDEsIlxcbWF0aGNhbHtVfSJdLFsyLDEsIlReKlxcbWF0aGNhbHtVfSJdLFszLDAsIlxcbWF0aGJie0F9XnsxfV9rIl0sWzEsMCwiaSxcXDtcXHRleHR7XFxub3JtYWxmb250IG9wZW59IiwyLHsic3R5bGUiOnsidGFpbCI6eyJuYW1lIjoiaG9vayIsInNpZGUiOiJib3R0b20ifX19XSxbMSwyLCJqLFxcO1xcdGV4dHtcXG5vcm1hbGZvbnQgY2xvc2VkfSIsMCx7InN0eWxlIjp7InRhaWwiOnsibmFtZSI6Imhvb2siLCJzaWRlIjoidG9wIn19fV0sWzEsMywiaiIsMCx7InN0eWxlIjp7InRhaWwiOnsibmFtZSI6Imhvb2siLCJzaWRlIjoidG9wIn19fV0sWzMsNCwiMCJdLFsyLDQsImRmIl0sWzEsNCwiIiwwLHsic3R5bGUiOnsibmFtZSI6ImNvcm5lciJ9fV0sWzIsNSwiZiJdXQ==
				\[\begin{tikzcd}[ampersand replacement=\&, column sep=1.8cm]
					X \& U \& {\mathcal{U}} \& {\mathbb{A}^{1}_k} \\
					\& {\mathcal{U}} \& {T^*\mathcal{U}}
					\arrow["{i,\;\text{\normalfont open}}"', hook', from=1-2, to=1-1]
					\arrow["{j,\;\text{\normalfont closed}}", hook, from=1-2, to=1-3]
					\arrow["j", hook, from=1-2, to=2-2]
					\arrow["\lrcorner"{anchor=center, pos=0.125}, draw=none, from=1-2, to=2-3]
					\arrow["f", from=1-3, to=1-4]
					\arrow["df", from=1-3, to=2-3]
					\arrow["0", from=2-2, to=2-3]
				\end{tikzcd}\]
				\begin{enumerate}[label=(\arabic*)]
					\item $i:U=\spec B\hookrightarrow X$ an open immersion through which $x$ factors;
					\item a closed immersion $j:U\to \mathcal{U}$ into a smooth variety $\mathcal{U}$ over $k$, a function $f:\mathcal{U}\to \bA^1_k$;
					\item an equivalence $\phi:U\simeq \dcrit(\mathcal{U},f)$;
					\item a cycle $\res_{JF}\in \pi_{-1}\infeqtrcohcpl{B}{k}{2}$ with null-homotopic underlying $2$-form;
					\item a homotopy $i^*\omega\stackrel{\eta}{\sim}\phi^*\omega_{f}+\res_{JF}$ in $\infsymp{-1}{k}{U}$, where $\omega_{f}$ is defined by the derived critical locus structure of $\dcrit(\mathcal{U},f)$ (Example \ref{exa: dcrit}).
				\end{enumerate}
				If $X$ is a derived Deligne--Mumford stack, the same statement holds except that {\normalfont Zariski} should be replaced by {\normalfont \'etale} everywhere.
			\end{restate}
		}
	\begin{proof}[Proof of Theorem \ref{main theorem}]
		
		Given $\omega\in \infsymp{-1}{k}{X}$ and $x\in |X|$, we can choose a Zariski neighborhood $U:=\spec B\xrightarrow{i}X$ that is small enough, so that Construction \ref{c: Lagrangian distribution} and \ref{c: Refining} (the first part regarding $R$) work at the same time. Therefore, there are closed immersions% https://q.uiver.app/#q=WzAsMyxbMCwwLCJVIl0sWzEsMCwiXFxzcGYgQSJdLFsxLDEsIlxcc3BlYyBSXFxtYXRocmxhcHs9OlxcbWF0aGNhbHtVfSx9Il0sWzAsMSwiaiciLDAseyJzdHlsZSI6eyJ0YWlsIjp7Im5hbWUiOiJob29rIiwic2lkZSI6InRvcCJ9fX1dLFswLDIsImoiLDIseyJzdHlsZSI6eyJ0YWlsIjp7Im5hbWUiOiJob29rIiwic2lkZSI6InRvcCJ9fX1dLFsxLDIsIlxcaW90YSJdXQ==
		\[\begin{tikzcd}[ampersand replacement=\&]
			U \& {\spf A} \\
			\& {\spec R\mathrlap{=:\mathcal{U},}}
			\arrow["{j'}", hook, from=1-1, to=1-2]
			\arrow["j"', hook, from=1-1, to=2-2]
			\arrow["\iota", from=1-2, to=2-2]
		\end{tikzcd}\]where $\spec R$ is smooth over $k$ and $A\simeq |\infcohcplfil{B}{R}|$, see Proposition \ref{prop: Lag Foliation}. Then, we can choose a jet-Frobenius residue $\res_{JF}$ and an exact structure $\lambda\in \infsympex{-1}{k}{\spec B}$ such that $i^*\omega\sim \res_{JF}+d_{\Inf}\lambda$ (Corollary \ref{cor: wJF infinite degree}). As in Construction \ref{c: Refining}, there exists a function $f\in R$ such that its image $\hat{f}$ in $A$ exhibits an isotropic structure of $\lambda$, and there is (after further shrinking $\mathcal{U}$) 
		\[\phi:U\stackrel{\hat{\phi}}{\simeq} \dcrit(\spf A,\hat{f})\simeq \dcrit(\mathcal{U},f),\]see Proposition \ref{prop: dCrit stk} and \ref{prop: filtered R vs ccA}. Furthermore, Proposition \ref{prop: loop} identifies $d_{\Inf}\lambda$ with the infinitesimal symplectic form $\hat{\phi}^*\omega_{\hat{f}} $ induced by the derived critical locus structure of $\dcrit(\spf A,\hat{f})$. Since $\iota$ is formally \'etale, there is $\phi^*\omega_{f}\simeq \hat{\phi}^*\omega_{\hat{f}}$.
	\end{proof}
	\begin{prop}\label{prop: perturbation of Res}
		Let $(X,\omega)$ be a $(-1)$-shifted infinitesimal symplectic derived scheme as in Theorem \ref{main theorem}. Suppose that $(i,j,f,\phi,\res_{JF},\eta)$ is a Darboux chart around $x\in |X|$. Given another Frobenius residue $\res'_{JF}$ of $\omega$, there exists a Darboux chart $(i,j,f',\phi,\res'_{JF},\eta')$ where
		\begin{enumerate}[label=(\arabic*)]
			\item the immersions $X\xhookleftarrow{i}\spec B\xhookrightarrow{j}\spec R$ are unchanged;
			\item there is $d_{\deRham}(f-f')=0$ so that $\phi$ is untouched:
			% https://q.uiver.app/#q=WzAsNCxbMCwwLCJcXHNwZWMgQiJdLFsxLDAsIlxcc3BlYyBSIl0sWzAsMSwiXFxzcGVjIFIiXSxbMSwxLCJUXipfe1xcc3BlYyBSfSJdLFswLDEsImoiXSxbMCwyLCJqIiwyXSxbMSwzLCJkZlxcc2ltZXEgZGYnIl0sWzIsMywiMCIsMl0sWzAsMywiIiwxLHsic3R5bGUiOnsibmFtZSI6ImNvcm5lciJ9fV1d
			\[\begin{tikzcd}[ampersand replacement=\&]
				{\spec B} \& {\spec R} \\
				{\spec R} \& {T^*_{\spec R};}
				\arrow["j", from=1-1, to=1-2]
				\arrow["j"', from=1-1, to=2-1]
				\arrow["\lrcorner"{anchor=center, pos=0.125}, draw=none, from=1-1, to=2-2]
				\arrow["{df\simeq df'}", from=1-2, to=2-2]
				\arrow["0"', from=2-1, to=2-2]
			\end{tikzcd}\]
			\item if $\res_{JF}$ and $\res_{JF}'$ are both in the image of $\prod_{n\ge N+1}\pi_0B^{(n)}$ under $\pi_0\cJ$, then $f'$ can be chosen in the way that $f-f'\in R^{(N)}$.
		\end{enumerate}
	\end{prop}
	\begin{proof}
	
	Choose $g,g'\in \prod_{n\ge N+1}\pi_0B^{(n)}$ representing $\res_{JF}$ and $\res_{JF}'$, respectively, and let $\psi\in \pi_0B^{(N+1)}$ be the sum of all components of $g-g'$, which is a finite sum. Then, consider the following commutative diagram,
	% https://q.uiver.app/#q=WzAsNyxbMCwwLCJCXnsoTisxKX0iXSxbMCwxLCJCXnsoMSl9Il0sWzEsMSwiRl97WzAsMV19XFxtYXRoYmJie1xcUGl9X3tCL2t9Il0sWzIsMSwiQV8yIl0sWzMsMSwiUiJdLFszLDAsIlJeeyhOKX0iXSxbMiwwLCJBXzJeeyhOKX0iXSxbMSwyLCJcXG1hdGhiYmJ7aX1fMSIsMl0sWzIsM10sWzQsMywiIiwwLHsic3R5bGUiOnsiaGVhZCI6eyJuYW1lIjoiZXBpIn19fV0sWzUsNF0sWzAsMV0sWzUsNiwiIiwyLHsic3R5bGUiOnsiaGVhZCI6eyJuYW1lIjoiZXBpIn19fV0sWzYsM10sWzAsNl1d
	\[\begin{tikzcd}[ampersand replacement=\&]
		{B^{(N+1)}} \&\& {A_2^{(N)}} \& {R^{(N)}} \\
		{B^{(1)}} \& {F_{[0,1]}\mathbbb{\Pi}_{B/k}} \& {A_2} \& R
		\arrow[from=1-1, to=1-3]
		\arrow[from=1-1, to=2-1]
		\arrow[from=1-3, to=2-3]
		\arrow[two heads, from=1-4, to=1-3]
		\arrow[from=1-4, to=2-4]
		\arrow["{\mathbbb{i}_1}"', from=2-1, to=2-2]
		\arrow[from=2-2, to=2-3]
		\arrow[two heads, from=2-4, to=2-3]
	\end{tikzcd}\]
	where the factorization $\mathbbb{i}_1$ is by \ref{c: jF power map}, and $R^{(N)}\to A_2^{(N)}$ is a surjection on $\pi_0$.
	Hence there is $\delta\in R^{(N)}$ whose image in $\pi_0A_2$ agrees with that of $\psi$. Finally, set $f':=f+\delta$.
	\end{proof}
	
	\section{Constructions of Examples}\label{sec: 5 exa}
	We now reproduce \cite[\S2]{pantev2013shifted} over general bases. The reasoning of AKSZ construction and Lagrangian intersection is formal and identical to the rational case. Interesting discussions are in \S\ref{sec: 5.3 perf}, where we construct a \textit{second Chern character $\ch_2$} for perfect complexes in characteristic $p>2$. This leads to a \textit{de Rham symplectic form} of degree $2$ over $\stperf$, the derived moduli stack of perfect complexes.
	
	\begin{remark}
		One can replace each $\mathbbb{\Pi}$ with the derived de Rham cohomology $\mathrm{DR}$ in \S5.1 and \S5.2, and all arguments remain valid for de Rham shifted symplectic forms.
	\end{remark}
	\subsection{AKSZ formalism}

	\begin{df}
		Let $X$ be a derived stack over $R\in\scr$. We say that $X$ is  \textit{$\cO$-compact} if, for every $A\in\scr_R$, the following hold:
		\begin{enumerate}[label=(\arabic*)]
			\item $\cO_{X_A}\in \QC(X_A)^{\omega}$ is a compact object;
			\item $\globalsection(X_A,-):=\hom_{X_A}(\cO_{X_A},-)$ sends perfect complexes over $X_A$ into $\perf_A\subset\m_A$.
		\end{enumerate}
	\end{df}
	\begin{prop}
		Given an $\cO$-compact derived stack $X$ over $R$,  there is a natural transformation of functors $\dSt^{op}_R\to \widehat{\fil}\dalg_R$,
		\begin{equation}\label{e: o-compact}
			\widehat{F\mathbbb{\Pi}}(-\times X/R)\to \widehat{F\mathbbb{\Pi}}(-/R)\otimes_R\globalsection(X,\cO_X).
		\end{equation}
	\end{prop}	
	\begin{proof}
		Since $\globalsection(X,\cO_X)$ is perfect over $R$, the functor $\globalsection(X,\cO_X)\otimes_R-$ preserves all small limits. Therefore, both sides of Line (\ref{e: o-compact}) are right Kan extensions of their restrictions to $\daff_R^{op}$. For every $A\in\scr_R$, there is a natural map of sheaves on $X$ taking value in $\fil\dalg_R$,
		\[F\mathbbb{\Pi}(A\otimes_R \cO_X/R)\simeq  \infcohfil{A}{R}\otimes_R F\mathbbb{\Pi}(\cO_X/R)\to \infcohfil{A}{R}\otimes_R\cO_X.\]
		Taking completed global sections, we get
		\[\infcohcplfil{X_A}{R}\to \widehat{\globalsection}(X,\infcohfil{A}{R}\otimes_R\cO_X).\]
		Meanwhile, since $X$ is $\cO$-compact, $\globalsection(X,-)$ commutes with $\infcohfil{A}{R}\otimes_R-$ by (1) and $(-)^{\wedge}$ commutes with $\globalsection(X,\cO_X)\otimes_R-$ by (2). Thus, there is a map functorial in $A\in\scr_R$ as follows,
		\[	\widehat{F\mathbbb{\Pi}}(X_A/R)\to \widehat{F\mathbbb{\Pi}}(A/R)\otimes_R\globalsection(X,\cO_X).\]
	\end{proof}
	Next, Pantev--To\"en--Vaqui\'e--Vezzosi requires $X$ to have a formal version of Poincar\'e duality:
	\begin{df}
		Let $X$ be an $\mathcal{O}$-compact derived stack over $R$. An \textit{$\cO$-orientation of degree $d$} is a morphism in $\m_R$
		\[\eta:\globalsection(X,\cO_X)\to R[-d]\]such that, for every $A\in\scr_R$ and $E\in \perf_{X_A}$, there is an $A$-linear equivalence
		\[-\cap\eta:\globalsection(X_A,E)\xrightarrow{\simeq} \globalsection(X_A,E^\vee)^\vee[-d].\]
	\end{df}
	\begin{exa}
		The Calabi--Yau $d$-varieties over $R$ are examples of $\cO$-compact stacks with an $\cO$-orientation of degree $d$. The $\cO$-orientation $[X]$ is given by the isomorphism $K_X\cong \cO_X$.
	\end{exa}
	\begin{construction}
		Let $X$ be an $\cO$-compact derived stack over $R\in\scr$ such that there is an $\cO$-orientation $\eta$ of degree $d$. Then, there is a formal integration
		\begin{equation}\label{e: integrating volume forms}
			\int_{\eta}:\widehat{F\mathbbb{\Pi}}(-\times X/R)\xrightarrow{\ref{e: o-compact}} \widehat{F\mathbbb{\Pi}}(-/R)\otimes_R\globalsection(X,\cO_X)\xrightarrow{id\otimes \eta}\widehat{F\mathbbb{\Pi}}(-/R)[-d].
		\end{equation}
	\end{construction}
	Here is a variant of \cite[Theorem 2.5]{pantev2013shifted}:
	\begin{theorem}\label{theorem: AKSZ}
		Let $F$ be a derived (locally) Artin stack over $R\in \scr$, and let $\omega\in \widehat{F_2\mathbbb{\Pi}}\mathstrut_{F/R}[n+2]$ be a shifted symplectic form of degree $n$. Suppose that $X$ is a derived $R$-stack that is $\cO$-compact and carries a $(-d)$-$\cO$-orientation\begin{equation}
			\eta:\globalsection(X,\cO_X)\to R[-d].
		\end{equation}If $M:=\map_{\dSt_R}(X,F)$ is again a derived (locally) Artin stack over $R$, then $M$ admits an infinitesimal symplectic form of degree $(n-d)$,
		\begin{equation}
			\int_{\eta}\omega\in\infsymp{n-d}{R}{M}.
		\end{equation}
	\end{theorem}
	\begin{proof}
	It is sufficient to show that $\int_\eta\omega$ is non-degenerate. For every $x:\spec A\to M$, the tangent complex at $x$ is given by
	\[x^*\dT_{M/R}\simeq \globalsection(X_A, f^*_x\dT_{F/R})\]following \cite[Corollary 2.2.6.14]{toen2008homotopical}, where $f_x:X_A\to F$ is adjoint to $x$. The underlying $2$-form of $\int_{\eta}\omega$ is then determined by
	\begin{align*}
		\wedge^2\globalsection(X_A,f^*_x\dT_{F/R})\xrightarrow{\text{the lax monoidality of }\globalsection}&\globalsection(X_A, \wedge^2_{\cO_{X_A}}f^*_x\dT_{F/R})\\
		\xrightarrow{\globalsection(X_A,\omega_A)}&\globalsection(X_A,\cO_{X_A})[n] \\
		\xrightarrow{\eta}&A[n-d].
	\end{align*}The induced map $x^*\dT_{M/R}\to x^*\dL_{M/R}[n-d]$ is homotopic to
	\[\globalsection(X_A,f^*_x\dT_{F/R})\xrightarrow{\globalsection(\Theta_\omega),\simeq}\globalsection(X_A,f^*_x\dL_{F/R})[n]\xrightarrow{-\cap\eta,\simeq}\globalsection(X_A,f^*_x\dT_{F/R})^\vee[n-d].\]
	\end{proof}

	\subsection{Lagrangian intersection}
	\begin{df}
		Let $f:X\to F$ be a morphism of derived Artin stacks over $R$, and let $\omega\in \infsymp{n}{R}{F}$ be an infinitesimal symplectic form of degree $n$. \begin{enumerate}[label=(\arabic*)]
		\item The space of \textit{isotropic structures} on $f$ relative to $\omega$ is defined by the cartesian square 
		% https://q.uiver.app/#q=WzAsNCxbMCwwLCJcXGlzb3R7Zn17XFxvbWVnYX0iXSxbMSwwLCJcXHswXFx9Il0sWzAsMSwiXFx7Zl4qXFxvbWVnYVxcfSJdLFsxLDEsIlxcaW5mZm9ybXR3b3tufXtSfXtYfSJdLFswLDJdLFswLDFdLFsxLDNdLFsyLDNdLFswLDMsIiIsMSx7InN0eWxlIjp7Im5hbWUiOiJjb3JuZXIifX1dXQ==
		\[\begin{tikzcd}[ampersand replacement=\&]
			{\isot{f}{\omega}} \& {\{0\}} \\
			{\{f^*\omega\}} \& {\infformtwo{n}{R}{X}}
			\arrow[from=1-1, to=1-2]
			\arrow[from=1-1, to=2-1]
			\arrow["\lrcorner"{anchor=center, pos=0.125}, draw=none, from=1-1, to=2-2]
			\arrow[from=1-2, to=2-2]
			\arrow[from=2-1, to=2-2]
		\end{tikzcd}\]
		
		\item If $\isot{f}{\omega}$ is nonempty, every object $h$ therein has an underlying isotropic distribution $f^*\dL_{F/R}[-1]\to \dL_{X/R}[-1]$. The isotropic structure $h$ is \textit{Lagrangian} if the induced map $\Theta_h:\dT_{f}\to \dL_X[n-1]$ is an equivalence.
		\end{enumerate}
	\end{df}
		
		\begin{theorem}\cite[Theorem 2.9]{pantev2013shifted} Given a cospan $X\xrightarrow{f}F\xleftarrow{g}Y$ of derived Artin stacks over $R$, the data of $\omega\in\infsymp{n}{R}{F}$, Lagrangian structures $h,k$ for $f,g$ jointly determine an infinitesimal symplectic form $R_{\omega,f,g}\in\infsymp{n-1}{R}{X\times_F Y}$.
		\end{theorem}
		\begin{proof}
			Set $Z:=X\times_Z Y$, and set $pr_X$, $pr_Y$ as the projections. Then, there is a loop in $\infformtwo{n}{R}{Z}$
			\[0\stackrel{h}{\sim}pr_X^*f^*\omega\sim pr_Y^*g^*\omega \stackrel{k}{\sim}0,\]i.e., a cycle $R_{\omega,h,k}\in\infformtwo{n-1}{R}{Z}$. It only remains to show that $R_{\omega,h,k}$ is non-degenerated. Recall that the isotropic structures of the underlying distributions are given by the commutative squares as follows,
		% https://q.uiver.app/#q=WzAsOCxbMCwxLCJmXipcXGRMX3tGL1J9W25dIl0sWzAsMCwiZl4qXFxkVF97Ri9SfSJdLFsxLDAsIlxcZFRfe2Z9WzFdIl0sWzEsMSwiXFxkTF97WC9SfVtuXSwiXSxbMywwLCJnXipcXGRUX3tGL1J9Il0sWzMsMSwiZ14qXFxkTF97Ri9SfVtuXSJdLFs0LDAsIlxcZFRfe2d9WzFdIl0sWzQsMSwiXFxkTF97WS9SfVtuXS4iXSxbMSwwLCJmXipcXG9tZWdhLFxcc2ltZXEiLDJdLFsxLDJdLFsyLDMsIlxcVGhldGFfaCxcXHNpbWVxIl0sWzAsM10sWzQsNSwiZ14qXFxvbWVnYSxcXHNpbWVxIiwyXSxbNCw2XSxbNSw3XSxbNiw3LCJcXFRoZXRhX2ssXFxzaW1lcSJdXQ==
		\[\begin{tikzcd}[ampersand replacement=\&]
			{f^*\dT_{F/R}} \& {\dT_{f}[1]} \&\& {g^*\dT_{F/R}} \& {\dT_{g}[1]} \\
			{f^*\dL_{F/R}[n]} \& {\dL_{X/R}[n],} \&\& {g^*\dL_{F/R}[n]} \& {\dL_{Y/R}[n].}
			\arrow[from=1-1, to=1-2]
			\arrow["{f^*\omega,\simeq}"', from=1-1, to=2-1]
			\arrow["{\Theta_h,\simeq}", from=1-2, to=2-2]
			\arrow[from=1-4, to=1-5]
			\arrow["{g^*\omega,\simeq}"', from=1-4, to=2-4]
			\arrow["{\Theta_k,\simeq}", from=1-5, to=2-5]
			\arrow[from=2-1, to=2-2]
			\arrow[from=2-4, to=2-5]
		\end{tikzcd}\] Hence, there is an equivalence of cofibre sequences
			
			% https://q.uiver.app/#q=WzAsNixbMCwxLCJcXHBpXipcXGRMX3tGL1J9W25dIl0sWzAsMCwiXFxwaV4qXFxkVF97Ri9SfSJdLFsxLDAsInByX1heKlxcZFRfe2Z9WzFdXFxvcGx1cyBwcl9ZXipcXGRUX3tnfVsxXSJdLFsxLDEsInByX1heKlxcZExfe1gvUn1bbl1cXG9wbHVzIHByX1leKlxcZExfe1kvUn1bbl0iXSxbMiwwLCJcXGRUX3taL1J9WzFdIl0sWzIsMSwiXFxkTF97Wi9SfVtuXSJdLFsxLDAsIlxccGleKlxcb21lZ2EsXFxzaW1lcSIsMl0sWzEsMl0sWzAsM10sWzIsMywicHJeKl9YXFxUaGV0YV9oK3ByXipfWVxcVGhldGFfayxcXHNpbWVxIl0sWzMsNV0sWzIsNF0sWzQsNSwiXFxUaGV0YV97XFxvbWVnYSxoLGt9Il1d
			\[\begin{tikzcd}[ampersand replacement=\&]
				{\pi^*\dT_{F/R}} \& {pr_X^*\dT_{f}[1]\oplus pr_Y^*\dT_{g}[1]} \& {\dT_{Z/R}[1]} \\
				{\pi^*\dL_{F/R}[n]} \& {pr_X^*\dL_{X/R}[n]\oplus pr_Y^*\dL_{Y/R}[n]} \& {\dL_{Z/R}[n],}
				\arrow[from=1-1, to=1-2]
				\arrow["{\pi^*\omega,\simeq}"', from=1-1, to=2-1]
				\arrow[from=1-2, to=1-3]
				\arrow["{pr^*_X\Theta_h+pr^*_Y\Theta_k,\simeq}", from=1-2, to=2-2]
				\arrow["{\Theta_{\omega,h,k}}", from=1-3, to=2-3]
				\arrow[from=2-1, to=2-2]
				\arrow[from=2-2, to=2-3]
			\end{tikzcd}\]where $\pi:=f\circ pr_X\simeq g\circ pr_Y$.
		\end{proof}
	\begin{exa}\label{exa: dcrit}
		Let $X$ be a Deligne--Mumford derived stack over $R$. There is a canonical infinitesimal symplectic form $\omega:=d_{\Inf}\Liouv$ of degree $0$ on $T^*_X$ (Proposition \ref{prop: Liouv}). Every function $f\in \pi_0\globalsection(X,\cO_X)$ induces a Lagrangian structure of
		\[df:X\to T^*_X.\]In particular, the derived critical locus $\dcrit(X,f)$, defined as the pullback
		% https://q.uiver.app/#q=WzAsNCxbMCwwLCJcXGRjcml0KFgsZikiXSxbMSwwLCJYIl0sWzAsMSwiWCJdLFsxLDEsIlReKl9YIl0sWzAsMl0sWzAsMV0sWzEsMywiZGYiXSxbMiwzLCIwIiwyXSxbMCwzLCIiLDEseyJzdHlsZSI6eyJuYW1lIjoiY29ybmVyIn19XV0=
		\[\begin{tikzcd}[ampersand replacement=\&]
			{\dcrit(X,f)} \& X \\
			X \& {T^*_X,}
			\arrow[from=1-1, to=1-2]
			\arrow[from=1-1, to=2-1]
			\arrow["\lrcorner"{anchor=center, pos=0.125}, draw=none, from=1-1, to=2-2]
			\arrow["df", from=1-2, to=2-2]
			\arrow["0"', from=2-1, to=2-2]
		\end{tikzcd}\]
		admits a $(-1)$-shifted infinitesimal symplectic form.
		
		The Lagrangian structure $X\xrightarrow{df}T^*_X$ contains \textbf{strictly more data} than the underlying closed immersion. For instance, when $X = \spec k[x]$ over a field $k$ of characteristic $p>0$, we have $d_{\mathrm{dR}}(f(x) + g(x^p)) = d_{\mathrm{dR}}f(x)$ for all $f,g \in k[x]$, yet they induce distinct Lagrangian structures.
	\end{exa}
	\subsection{2-Shifted de Rham symplectic form on $\stperf$}\label{sec: 5.3 perf}
	This subsection is devoted to constructing a de Rham symplectic form of degree $2$ over $\stperf$, the moduli stack of perfect complexes, in characteristic $p>2$.
	
	In characteristic $0$, the natural $2$-shifted symplectic form over $\stperf$ was constructed in \cite[\S2.3]{pantev2013shifted} using Chern character. Recall that To\"en--Vezzosi constructed a categorical Chern character in \cite{toen2015caracteres}
	\begin{equation}
		\ch:\stperf\to \Omega^\infty\HC^{-}(-/R)
	\end{equation}as a morphism of \'etale hypersheaves over $\daff_R$ for a fixed base ring $R\in\scr$, where $\HC^-$ is the \textit{negative cyclic homology} relative to $R$. Additionally, this morphism is additive and multiplicative, see \cite[\S2.4, \S2.5]{toen2015caracteres}. If $R$ is a $\mathbb{Q}$-algebra, there is a well-known HKR equivalence (see \cite{toen2011algebres})
	\begin{equation}\label{e: Q-coeff Ch}
		\HC^-(-/R)\simeq \bigoplus_{n\ge0} \what{F_n\mathrm{DR}}(-/R)[2n],
	\end{equation}
	where both ends are \'etale hypersheaves over $\daff_R$. Furthermore, this categorical Chern character is compatible with the classical one in de Rham cohomology and can be determined by powers of the Atiyah class, see \cite[Appendix B]{toen2015caracteres}\cite[Theorem 2.12]{pantev2013shifted}. The Atiyah class description shows that the second piece $\ch_2\in \pi_{-4}\what{F_2\mathrm{DR}}(\stperf/R)$ is a $2$-shifted symplectic form.
	
	However, it is obscure what $\ch$ is like away from $\mathbb{Q}$. One crucial difficulty is that line (\ref{e: Q-coeff Ch}) does not hold. Fortunately, there is still an HKR theorem that works in general characteristics.
	
	\begin{recoll}
		Let $R$ be a simplicial commutative ring. By \cite{antieau2019periodic,moulinos2022universal,raksit2020hochschild}, the HKR theorem over general bases states that there exists an \textit{HKR filtration} on $\HH(-/R)$ enhancing it to a functor valued in complete derived algebras
		\[F\HH(-/R):\scr_R\to\what{\fil}\dalg_R\]and equipped with a filtered $S^1$-action. Moreover, $\gr\HH(A/R)$ together with the graded circle action\footnote{A graded $S^1$-action is equivalent to an $R\oplus R[1]$-module structure.} can be identified with the Hodge completed derived de Rham cohomology
		\[\big(\wedge^*_A\dL_{A/R}[+*],d_{\deRham}\big)\]
		in the ambient category of \textit{graded mixed derived algebras} in the sense of \cite{moulinos2022universal}, or equivalently, \textit{$h^+$-derived algebras} \cite{raksit2020hochschild}. In addition, the graded pieces of $F\HC^-:=F\HH^{hS^1}$ are given by
		\begin{equation}
			\gr_n\HC^-(A/R)\simeq \what{F_n\mathrm{DR}}_{A/R}[2n],
		\end{equation}although the HKR filtration no longer splits.
	\end{recoll}
	We fix the issue by showing that, in characteristic $p>2$, the truncated HKR filtration on Hochschild homology splits functorially.
	\begin{prop}\label{prop: partial splitting of HH}
		
		Suppose that $R\in\mathrm{Ring}_{\bF_p}$ with $p>2$. There is an equivalence of $R$-modules functorial in $A\in\scr_R$,
		\begin{equation}
			F_{[a,a+p-2]}\HH(A/R)\simeq \bigoplus_{a\le n\le a+p-2}\bigwedge^n_A\dL_{A/R}[n].
		\end{equation}
	\end{prop}
	
	\begin{remark}
		This proposition holds for every $R\in\scr_{\bF_p}$ by left Kan extension, since both sides define sifted-colimit-preserving functors from $\scr_{\bF_p}^{\Delta^1}$.
	\end{remark}

	The proof requires us to dive into the details of \cite{moulinos2022universal}.
	\begin{recoll}
		Consider the ordinary ring scheme of $p$-typical Witt vectors
		\[\stwitt_{\pinfty}:\mathrm{Ring}_{R}\to \mathrm{Ring}_{\witt_{\pinfty}(R)}.\]There is a natural ring scheme endomorphism $\frob:\stwitt_{\pinfty}\to \stwitt_{\pinfty}$. Let $\Fix$ and $\Ker$ be the kernels of $\frob-id$ and $\frob$, respectively. Moulinos--Robalo--To\"en recovered $\HH(A/R)$ and $\gr\HH(A/R)$ as the algebras of global sections of the mapping stacks
		\begin{equation}
			\map_{\dSt_R}(B\Fix,\spec A),\ \ \ \ \ \map_{\dSt_R}(B\Ker,\spec A),
		\end{equation}\cite[Theorem 5.4.1]{moulinos2022universal}. Their construction of the HKR filtration relies on a filtered stack $\mathrm{H}_{\pinfty}\to [\bA^1/\bG_m]$ interpolating $\Fix$ and $\Ker$ \cite[Construction 2.3.4]{moulinos2022universal}:
		There is a natural $\bG_m$-action on $\stwitt_{\pinfty}$ given by multiplication by the Teichm\"uller character $[-]$. Consider the trivial group scheme $\stwitt\times \bA^1$ over $\bA^1$. There is a group scheme endomorphism defined by, for every $S\in \mathrm{Ring}_R$,
		\[\mathscr{G}_p:\witt_{\pinfty}(S)\times S\to \witt_{\pinfty}(S)\times S: (f, a)\mapsto (\frob(f)-[a^{p-1}f],a).\]The kernel $\ker(\mathscr{G}_p)$ is closed under the diagonal action of $\bG_m$ on $\stwitt\times\bA^1$, and the quotient $\mathrm{H}_{\pinfty}:=[\ker(\mathscr{G}_p)/\bG_m]$ is the desired filtered stack in \cite[Definition 2.3.7]{moulinos2022universal}. Furthermore, the filtered loop stack \cite[Definition 5.1.1]{moulinos2022universal}
		\begin{equation}\label{e: filtered loop stack}
			\map_{\dSt_R}(B\mathrm{H}_{\pinfty},\spec A)\to [\bA^1/\bG_m]
		\end{equation}encodes the HKR filtration by the Rees construction.
	\end{recoll}
	
	\begin{proof}[Proof of Proposition \ref{prop: partial splitting of HH}]Regard $\bF^\times_p$ as a constant sheaf. There are inclusions $\bF^\times_p\hookrightarrow \bG_m\xhookrightarrow{[-]}\stwitt_{\pinfty}$ of group schemes over $R$. Then the finite group $\bF^\times_p$ acts on $\stwitt_{\pinfty}$ by multiplying the Teichm\"uller characters. Since $\frob\in\mathrm{End}(\witt_{\pinfty}(A))$ can be written as $(a_n)_{n\in\mathbb{N}}\mapsto (a_n^p)_{n\in\mathbb{N}}$ in the standard coordinates, the action of $\bF^\times_p$ commutes with $\frob$. Then, $\bF^\times_p$ acts on $\ker(\mathscr{G}_p)$ by  acting on \textbf{only} the first component $\stwitt_{\pinfty}$, and this action descends to a relative action on $B\mathrm{H}_{\pinfty}\to[\bA^1/\bG_m]$. The induced action on the filtered loop stack (\ref{e: filtered loop stack}) defines an $\bF^\times_p$-action on the HKR-filtered $\HH$\[F\HH(A/R)\in\fil\dalg_R.\]
		
	We claim that each $\gr_n\HH(A/R)$ is pure of weight $n$ (mod $p-1$) under the $\bF_p^\times$-action. Let $R[\epsilon]$ be the ordinary ring of dual numbers. The point $[\epsilon]\in \Ker(R[\epsilon])$ is given by the Teichm\"uller character. This gives rise to a commuting square of stacks
	% https://q.uiver.app/#q=WzAsNCxbMCwwLCJcXHNwZWMoUltcXGVwc2lsb25dKSJdLFswLDEsIlxcc3BlYyhSKSJdLFsxLDAsIlxcc3BlYyhSKSJdLFsxLDEsIkJcXEtlciJdLFswLDFdLFsxLDNdLFswLDJdLFsyLDNdLFsxLDIsIltcXGVwc2lsb25dLFxcc2ltZXEiLDEseyJzdHlsZSI6eyJib2R5Ijp7Im5hbWUiOiJub25lIn0sImhlYWQiOnsibmFtZSI6Im5vbmUifX19XV0=
	\[\begin{tikzcd}[ampersand replacement=\&]
		{\spec(R[\epsilon])} \& {\spec(R)} \\
		{\spec(R)} \& {B\Ker.}
		\arrow[from=1-1, to=1-2]
		\arrow[from=1-1, to=2-1]
		\arrow[from=1-2, to=2-2]
		\arrow["{[\epsilon],\simeq}"{description}, draw=none, from=2-1, to=1-2]
		\arrow[from=2-1, to=2-2]
	\end{tikzcd}\]Let $\bF^\times_p$ act on $R[\epsilon]$ by $\lambda:r_1+r_2\epsilon\mapsto r_1+\lambda r_2\epsilon$. The above square is then $\bF^\times_p$-equivariant. Now consider the mapping stacks to $\spec A$. There is an $\bF^\times_p$-equivariant equivalence
	\[\map_{\dSt_R}(B\Ker,\spec A)\simeq \Omega_*\map_{\dSt_R}(\spec (R[\epsilon]),\spec A)\simeq  \spec \lsym_A(\dL_{A/R}[1])\]by \cite[equation (120)]{moulinos2022universal}, where $\dL_{A/R}[1]$ is of weight $1$.
	\end{proof}
	Next, we split $\ch:\stperf_R\to \Omega^\infty\HC^-(-/R)$ for lower degrees using algebraic K-theory.
	\begin{construction}
		It is more convenient to consider a more general context. Let $\mathcal{C}\mathrm{at}^{\perf}_{\infty}$ be the \infcat\ of the small idempotent-complete stable \infcats\ and exact functors. The universal property of connective K-theory induces a commutative triangle natural in $\mathcal{C}\in \mathcal{C}\mathrm{at}^{\perf}_{\infty}$% https://q.uiver.app/#q=WzAsMyxbMCwwLCJcXG1hdGhsbGFwe1xcbWF0aGNhbHtDfV57XFxzaW1lcX1cXHNpbWVxfSBcXG1hcF97XFxtYXRoY2Fse0N9XFxtYXRocm17YXR9XntcXHBlcmZ9X3tcXGluZnR5fX0oXFxTcF57XFxvbWVnYX0sXFxtYXRoY2Fse0N9KSJdLFsyLDAsIlxcT21lZ2Fee1xcaW5mdHl9XFxtYXRocm17VEhIfShcXG1hdGhjYWx7Q30pIl0sWzEsMSwiXFxPbWVnYV57XFxpbmZ0eX1LKFxcbWF0aGNhbHtDfSkiXSxbMCwyXSxbMCwxXSxbMiwxLCJcXE9tZWdhXlxcaW5mdHkoXFx0cmFjZV97XFxvcGVyYXRvcm5hbWV7XFx0ZXh0aXR7RGVubmlzfX19KSIsMl1d
		\[\begin{tikzcd}[ampersand replacement=\&]
			{\mathllap{\mathcal{C}^{\simeq}\simeq} \map_{\mathcal{C}\mathrm{at}^{\perf}_{\infty}}(\Sp^{\omega},\mathcal{C})} \&\& {\Omega^{\infty}\mathrm{THH}(\mathcal{C})} \\
			\& {\Omega^{\infty}K(\mathcal{C})}
			\arrow[from=1-1, to=1-3]
			\arrow[from=1-1, to=2-2]
			\arrow["{\Omega^\infty(\trace_{\operatorname{\textit{Dennis}}})}"', from=2-2, to=1-3]
		\end{tikzcd}\]where $\trace_{\operatorname{\textit{Dennis}}}$ is the Dennis trace, and the left arrow is induced by the unit map $\mathbb{S}\to K(\mathbb{S})\simeq K(\Sp^{\omega})$. See \cite[\S10]{blumberg2013universal} for details. Fix an $R\in \scr$. When $\mathcal{C}=\perf_A$, there is a natural map
		\[\mathrm{THH}(\perf_A)\simeq\mathrm{THH}(A)\to \HH(A/R).\]Moreover, the composite $\stperf(A)\to \Omega^{\infty}K(A)\to \Omega^{\infty}\HH(A/R)$ agrees with the trace map in \cite{toen2015caracteres}. Since the Dennis trace is $S^1$-invariant, the Chern character factorizes into
		\begin{equation}\label{e: K-theory and Dennis trace}
			\ch:\stperf(A)\to \Omega^\infty K(A)\to \Omega^\infty\HC^-(A/R).
		\end{equation}
		
		Let $\mon$ and $\grp$ be the \infcats\ of $\mathbb{E}_\infty$-monoids and $\mathbb{E}_\infty$-groups, respectively. There is an equivalence $\grp\simeq \Sp_{\ge0}$, and $K(A)$ can be regarded as an $\mathbb{E}_\infty$-group. The connective K-theory $K(A)$ is the group completion of the $\mathbb{E}_\infty$-monoid $(\vect^{\omega,\simeq}_A,\oplus)$. Hence, there is an equivalence
		\begin{equation}\label{restricting Dennis trace}
			\map_{\mathrm{PSH}(\daff_R,\grp)}(K,F)\simeq \map_{\mathrm{SH}^{\acute{e}t}(\mathrm{Aff}_R,\mon)}(\stvect,F),
		\end{equation}
		where $F$ is any sheaf of $\mathbb{E}_\infty$-groups, and $\stvect\simeq\coprod_{n\ge0}B\GL_n$  is the classifying stack of vector bundles. In particular, $\ch$ is determined by an additive map $\stvect\to \HC^-(-/R)$.
	\end{construction}
	 
	Write  $\map_{\mathrm{SH}^{\acute{e}t}(\daff_R,\mon)}$ as $\map^{\oplus}$ for short.
	\begin{prop}\label{prop: HC to HH}
	For an ordinary ring $R$, the comparison map	\[\map^{\oplus}(\stvect,F_{[a,b]}\HC^-(-/R))\xrightarrow{\simeq}\map^{\oplus}(\stvect,F_{[a,b]}\HH(-/R))\]is an equivalence of $0$-truncated spaces, i.e., merely sets.
	\end{prop}
	\begin{alert}
		We regard $\HH(-/R)$ and $\HC^-(-/R)$ as sheaves on $\daff_R$. When $X$ is a genuine derived Artin stack, $\HH(X/R)$ is rather the global sections of the \textit{formal completion} of the loop stack $\map(S^1,X)$ along the locus of constant loops.
	\end{alert}
	\begin{proof}
	By Lemma \ref{lemma: operadic app to monoid}, it suffices to show that, for every $m\ge 0$, the mapping spaces from $\stvect^{\times m}$ to $F_{[a,b]}\HC^-$ and to $F_{[a,b]}\HH$ are equivalent and $0$-truncated. 
	
	By \cite[\S10]{totaro2018hodge}, the Hodge-to-de Rham spectral sequence of $B\GL_n$ degenerates and there is an isomorphism of graded commutative rings
	\begin{equation}\label{e: bhatt and Totaro}
		\pi_*\what{\mathrm{DR}}(B\GL_n/\mathbb{Z})\simeq \mathbb{Z}[c_1,\ldots,c_n],
	\end{equation}
	where the $i$-th Chern class $c_i$ is of degree $-2i$. The HKR theorem \cite[Theorem 1.2.1]{raksit2020hochschild} then implies a non-canonical equivalence of filtered simplicial commutative rings
	\[F\HH(B\GL_n/\mathbb{Z})\simeq \mathbb{Z}\formalpower{c_1,\ldots,c_n},\]where $c_i$ is of degree $0$ and weight $i$. After base change to $R$, and taking the connective cover of $S^1$-invariants, the composite
	\[\tau_{\ge0}F\HC^-(B\GL_n/R)\to F\HC^-(B\GL_n/R)\to F\HH(B\GL_n/R)\]is an equivalence of discrete objects in $\fil\scr_R$. The underlying spaces of the $F_{[a,b]}$-pieces are
	\begin{equation}\label{i am dead}
		\map_{\dSt_R}(B\GL_n,F_{[a,b]}\HC^-(-/R))\xrightarrow{\simeq}\map_{\dSt_R}(B\GL_n,F_{[a,b]}\HH(-/R)),
	\end{equation}which is equivalent to $F_{[a,b]}R\formalpower{c_1,\ldots,c_n}$. Since (\ref{e: bhatt and Totaro}) is finite in each weight, the K\"unneth formula implies that $\pi_*\what{\mathrm{DR}}(B\GL_{n_1}\times\cdots\times B\GL_{n_m}/\mathbb{Z})$ is given by a straightforward tensor product over $\mathbb{Z}$. Hence we can prove a statement similar to (\ref{i am dead}) for any product $B\GL_{n_1}\times\cdots\times B\GL_{n_m}$, and then for any $\stvect^{\times m}$.
	\end{proof}

		\begin{lemma}\label{lemma: operadic app to monoid}
		Assume that $\mathcal{C}$ is an \infcat\ with finite products. Let $X$ and $Y$ be $\mathbb{E}_\infty$-monoids in $\mathcal{C}$. If every $\map_{\mathcal{C}}(X^{\times n},Y)$ is discrete, then the mapping space of $\mathbb{E}_\infty$-monoids $\map_{\mon(\mathcal{C})}(X,Y)$	is also discrete and can be identified with a subset of $\pi_0\map_{\mathcal{C}}(X,Y)$.	
		\end{lemma}
		\begin{proof}
			Following \cite[Proposition 2.4.2.5]{HA}, it suffices to consider the \textit{cartesian $\infty$-operad} $q:\mathcal{C}^\times\to N(\mathrm{Fin}_*)$ determined by $\mathcal{C}$ in the sense of \cite[\S2.1.1, \S2.4.1]{HA}. The $\mathbb{E}_\infty$-monoids are by definition the operadic sections of $q$. Regarding $Y$ as such a section, it gives rise to a slice $\infty$-operad $\mathcal{C}^{\times}_{/Y}\to \mathcal{C}^{\times}$ by \cite[Theorem 2.2.2.4]{HA}, whose fibre over $[n]\in N(\mathrm{Fin}_*)$ is $(\mathcal{C}_{/Y})^{\times n}$. Then, the mapping space $\map_{\mon(\mathcal{C})}(X,Y)$ is equivalent to the space of the $\infty$-operadic sections $f$
			% https://q.uiver.app/#q=WzAsNSxbMCwxLCJOKFxcbWF0aHJte0Zpbn1fKikiXSxbMSwxLCJcXG1hdGhjYWx7Q31ee1xcdGltZXN9Il0sWzEsMCwiXFxtYXRoY2Fse0N9XntcXHRpbWVzfV97L1l9Il0sWzEsMiwiTihcXG1hdGhybXtGaW59XyopIl0sWzAsMCwiXFxtYXRoY2Fse0R9XntcXG90aW1lc30iXSxbMCwxLCJYIl0sWzIsMV0sWzEsMywicSIsMl0sWzAsMywiIiwxLHsibGV2ZWwiOjIsInN0eWxlIjp7ImhlYWQiOnsibmFtZSI6Im5vbmUifX19XSxbNCwwXSxbNCwyXSxbNCwxLCIiLDEseyJzdHlsZSI6eyJuYW1lIjoiY29ybmVyIn19XSxbMCw0LCJmIiwxLHsib2Zmc2V0IjotMywic3R5bGUiOnsiYm9keSI6eyJuYW1lIjoiZGFzaGVkIn19fV1d
			\[\begin{tikzcd}[ampersand replacement=\&]
				{\mathcal{D}^{\otimes}} \& {\mathcal{C}^{\times}_{/Y}} \\
				{N(\mathrm{Fin}_*)} \& {\mathcal{C}^{\times}} \\
				\& {N(\mathrm{Fin}_*)}
				\arrow[from=1-1, to=1-2]
				\arrow[from=1-1, to=2-1]
				\arrow["\lrcorner"{anchor=center, pos=0.125}, draw=none, from=1-1, to=2-2]
				\arrow[from=1-2, to=2-2]
				\arrow["f"{description}, bend left, shift left=3, dashed, from=2-1, to=1-1]
				\arrow["X", from=2-1, to=2-2]
				\arrow[equals, from=2-1, to=3-2]
				\arrow["q"', from=2-2, to=3-2]
			\end{tikzcd}\]
			where $\mathcal{D}^\otimes$ is defined by a cartesian square. One object in $\mathcal{D}^{\otimes}_{[n]}$ is an $n$-tuple $(f_1,\ldots,f_n)$ of maps $f_i:X\to Y$. The mapping space $\map_{\mathcal{D}^{\otimes}}\big((f_1,\ldots,f_n),g\big)$ is the space of paths between $\mu_{Y,n}\circ (f_1,\ldots,f_n)$ and $g\circ \mu_{X,n}$ in $\map_{\mathcal{C}}(X^{\times n},Y)$, where $\mu_{X,n}$ and $\mu_{Y,n}$ are the structure maps of $X$ and $Y$. This means that $\mathcal{D}^{\otimes}$ is a discrete operad, and $\map_{\mon(\mathcal{C})}(X,Y)\simeq \calg(\mathcal{D})$ is discrete.\end{proof}
	\begin{corollary}\label{cor: split ch}
		For any $R\in\scr_{\bF_p}$, the truncated Chern character \[\ch:\stperf\to \Omega^{\infty}F_{[0,p-1]}\HC^-(-/R)\] {\normalfont canonically splits} as the sum of the images of certain additive maps
		\[\ch_i:\stperf\to\Omega^\infty F_{[i,p-1]}\HC^-(-/R).\]
	\end{corollary}
	\begin{proof}
		When $R=\bF_p$, $F_{[1,p-1]}\HH(-/\bF_p)$ splits naturally as $\bigoplus\limits_{1\le n\le p-1}\wedge^n\dL_{-/\bF_p}[n]$ by Proposition \ref{prop: partial splitting of HH}. At the same time, the HKR filtration induces an augmentation $F_{[0,p-1]}\HH(-/\bF_p)\simeq \bF_p\oplus F_{[1,p-1]}\HH$. Thus, the truncated trace map $\stvect\to \Omega^{\infty}F_{[0,p-1]}\HH(-/R)$ splits naturally. Consequently, there is a natural splitting of $\ch$ by Proposition \ref{prop: HC to HH}, (\ref{e: K-theory and Dennis trace}) and (\ref{restricting Dennis trace}). For general $R\to \scr_{\bF_p}$, it suffices to post-compose $F_{[i,p-1]}\HC^-(-/\bF_p)\to F_{[i,p-1]}\HC^-(-/R)$.
	\end{proof}
	We now aim to describe the underlying forms of $\ch_i$ via the Atiyah class.
	\begin{construction}[Atiyah class]\label{c: atiyah}
		Let $\cE$ be the universal perfect complex over $\stperf$. Recall that the global cotangent complex of $\stperf$ is
		\begin{equation}
			\dL_{\stperf}\simeq \cE^\vee\otimes_{\cO}\cE[-1]
		\end{equation}by \cite[Corollary 3.17]{toen2007moduli}. The \textit{universal Atiyah class} is defined as
		\begin{equation}
			\atiyah_{\cE}:\cE\xrightarrow{\trace\otimes id_{\cE}} \cE\otimes_\cO\cE^\vee\otimes_\cO\cE\simeq \cE\otimes_\cO \dL_{\stperf}[1].
		\end{equation}
		This Atiyah class is dual to the map $\cO\to \cE^\vee\otimes_\cO\cE\otimes_\cO\dL_{\stperf}[1]$, and it gives rise to a class $\trace(\atiyah_{\cE}):\cO\to \dL_{\stperf}[1]$ by post-composing with the evaluation of $\cE$, which is essentially the trace map $\trace_\cE:\cO\to \cE^\vee\otimes_\cO\cE$. More generally, we can consider the $m$-th iteration of the Atiyah class
		\[\atiyah^m_\cE:\cE\to\cE\otimes_\cO\dL_{\stperf}[1]\to \ldots \cE\otimes_\cO(\dL_{\stperf}[1])^{\otimes m}\to \cE\otimes_\cO\wedge^m\dL_{\stperf}[m],\]and its trace map $\trace(\atiyah_{\cE}^m)$, which agrees with $\trace_{\cE}^{m}:\cO\xrightarrow{\lsym^m(\trace_{\cE})}\lsym^m_\cO(\cE^\vee\otimes_\cO\cE)$.
	\end{construction}
	Recall that the HKR theorem implies that $\gr_i\HC^-(-/R)\simeq \what{F_i\mathrm{DR}}(-/R)[i]$. Therefore, we can define the \textit{$i$-th de Rham Chern character} $\ch^{\deRham}_i$ as the composite
	\[\stperf\xrightarrow{\ch_i}\Omega^\infty F_{[i,p-1]}\HC^-(-/R)\to\Omega^\infty \what{F_i\mathrm{DR}}(-/R)[i].\vspace{-1em}\]
	\begin{prop}\label{prop: ch_i as trace map}
		Suppose that $R\in\scr_{\bF_p}$ and $p>2$. The underlying $i$-form of $\ch^{\deRham}_i$ ($1\le i\le p-1$) is given by $\frac{1}{i!}\trace(\atiyah^i_\cE)$.
	\end{prop}
	\begin{proof}
		The loop stack $\cL\stperf:=\map(S^1,\stperf)$ classifies the automorphisms of perfect complexes. There is a natural projection $\varpi:\cL\stperf\to\stperf$ and a universal automorphism
		\begin{equation}
			a:\varpi^*\cE\to\varpi^*\cE.
		\end{equation}By To\"en--Vezzosi's construction in \cite[\S2.3]{toen2015caracteres}, the underlying cycle of $\ch$ in $\HH(\stperf)$ is determined by $\trace(a)$. At the same time, the identity of $\cE$ induces a section   $c:\stperf\to \cL\stperf$ of $\varpi$. The first order infinitesimal neighborhood of $c$ is an affine stack relative to $\stperf$ as follows
		\begin{equation}
			\spec_{\cO}(\cO\oplus \cE\otimes\cE^\vee)\xrightarrow{emb}\cL\stperf.
		\end{equation}The restriction $emb^*a$ is adjoint to the following map over $\stperf$ (along $\varpi\circ emb$)
		\begin{equation}
			\cE\xrightarrow{id_\cE\oplus \atiyah_\cE}\cE\otimes_\cO(\cO\oplus\cE^\vee\otimes_\cO\cE).
		\end{equation}
		
		After base change to $R$, we further restrict to the open substack $i_n:B\GL_n\to \stperf$. Let $\cV:=i^*_n\cE$ be the universal vector bundle of rank $n$. There is a map of filtered stacks over $B\GL_n$
		\begin{equation}\label{e: deformation to tangent GL_n}
			\widetilde{emb}_p:\spec_\cO(\lsym^{\le p-1}_\cO\cV^\vee\otimes_\cO \cV)\times[\bA^1/\bG_m]\to \map_{\dSt_R}(B\mathrm{H}_{\pinfty},B\GL_n).
		\end{equation} by Proposition \ref{prop: partial splitting of HH}, whose generic fibre
		\begin{equation}
			\widetilde{emb}_p:\spec_\cO(\lsym^{\le p-1}_\cO\cV^\vee\otimes_\cO \cV)\to \cL B\GL_n
		\end{equation}
		is precisely the $(p-1)$-th infinitesimal neighborhood of $c|_{B\GL_n}:B\GL_n\to \cL B\GL_n$. Moreover, the truncation $q(\ch|_{B\GL_n})$ by $q:\HH(B\GL_n/R)\to F_{[0,p-1]}\HH(B\GL_n/R)$ is homotopic to $\trace(\widetilde{emb}_p^*a_\cV)$, where $a_\cV$ is the universal equivalence of $\cV$.
		\begin{lemma}\label{lemma: rescaling of trace}
			Let $pr_i:F_{[0,p-1]}\HH(B\GL_n/R)\to \globalsection(B\GL_n,\wedge^i\dL_{B\GL_n/R})[i]$ ($0\le i\le p-1$) be the projection in Proposition \ref{prop: partial splitting of HH}. There is a homotopy $pr_i\trace(\widetilde{emb}_p^*a^l_\cV)\simeq l^i\trace(\widetilde{emb}_p^*a_\cV)$.
		\end{lemma}
		\begin{proof}[Proof of Lemma \ref{lemma: rescaling of trace}]
			The assignment $a_\cV\mapsto a^l_\cV$ corresponds to the $l$-fold covering $S^1\to S^1$, which can be modeled as $\times l:B\mathbb{Z}\to B\mathbb{Z}$ on the level of homotopy types. Meanwhile, the natural additive map $\mathbb{Z}\to \stwitt_{\pinfty}$ sending $1$ to $(1,0,0,\ldots)$ (standard coordinates) exhibits $B\Fix$ as the \textit{affinization} of $B\mathbb{Z}$ over $R$, see \cite{toen2006champs} and \cite[Proposition 3.3.2]{moulinos2022universal}. This means that $a_\cV\mapsto a^l_\cV$ is induced by the monoid action of $\bF_p$ on $\stwitt$ extending the $\bF^\times_p$-action in the proof of Proposition \ref{prop: partial splitting of HH}. Observe that this monoid action can be extended to compatible actions on both sides of Line (\ref{e: deformation to tangent GL_n}), where the action on the left is constant with respect to $[\bA^1/\bG_m]$. However, we have seen in the proof of Proposition \ref{prop: partial splitting of HH} that the $\bF_p$-action on the special fibre is the rescaling of the cotangent complex.
		\end{proof}
		Consider the map over $B\GL_n$ adjoint to $\widetilde{emb}_p^*a_\cV$,
		\[\cV\to\cV\otimes_{\cO}\lsym^{\le p-1}_\cO(\cV^\vee\otimes_\cO\cV).\]Let $h_i$ be the $i$-th term of this map.
		The discussion of the Atiyah class has shown that the first two terms are $id_{\cV}$ and $\atiyah_{\cV}$. Suppose that we have shown that $h_j\simeq\frac{1}{j!}\atiyah^j_\cV$ for all $j<i\le p-1$. Then, we consider the composite $\widehat{emb}^*_pa^2_{\cV}$, whose $i$-th term is homotopic to
		\[2h_i+\sum_{0<j<i}\frac{1}{(i-j)!j!}\atiyah_{\cV}^j\atiyah_{\cV}^{i-j}.\]However, Lemma \ref{lemma: rescaling of trace} tells us that this term is homotopic to $2^ih_i$. Thus, we have $h_i\simeq \frac{1}{i!}\atiyah^i_{\cV}$.
	\end{proof}
	\begin{theorem}\label{theorem: perf with symplectic structure}
		Let $R$ be a simplicial commutative ring of characteristic $p>2$. There exists a second Chern character $\ch^{\deRham}_2\in\pi_{-4} \what{F_2\mathrm{DR}}(\stperf/R)$, which defines a {\normalfont de Rham symplectic form of degree $2$} on $\stperf$, the derived classifying stack of perfect complexes over $R$. 
	\end{theorem}
	\begin{proof}
		Since $p>2$, there is a closed $2$-form $\ch^{\deRham}_2$ of degree $2$ in de Rham cohomology, constructed by splitting the categorical Chern character via the HKR theorem and Corollary \ref{cor: split ch}. Its underlying $2$-form is given by $\frac{1}{2}\atiyah^2$. This $2$-form defines
		\[\bigwedge^2_{\cO}\dT_{\stperf}[1]\simeq \lsym^2_{\cO}(\cE\otimes_\cO\cE^\vee)[2]\to \cO[2]\] by $(A,B)\mapsto \trace(AB)$, which is non-degenerate.
	\end{proof}

	\begin{corollary}\label{cor: perf_X}
		Let $X\to \spec R$ be a flat family of smooth proper $3$-dimensional Calabi--Yau varieties, where $R\in \scr_{\bF_p}$ and $p>2$. The derived moduli stack of perfect complexes over $X$
		\[\stperf_X:=\map_{\dSt_R}(X,\stperf)\]
		is locally Artin and carries a natural de Rham symplectic form $\omega_{\deRham}:=\int_{[X]}\ch^{\deRham}_2$ of degree $-1$. 
	\end{corollary}
	\begin{proof}
		Following \cite[Proposition 3.13, Corollary 3.29]{toen2007moduli}, the derived stack $\stperf_X$ is an increasing union of finitely presented Artin stacks $\stperf_X^{[a,b]}$. Therefore, we can apply the de Rham version of Theorem \ref{theorem: AKSZ} together with Theorem \ref{theorem: perf with symplectic structure}.
	\end{proof}
	{\bibliographystyle{alpha}
	\bibliography{bib.bib}}

\begin{thebibliography}{PTVV13}

\bibitem[Ant19]{antieau2019periodic}
Benjamin Antieau.
\newblock Periodic cyclic homology and derived de {R}ham cohomology.
\newblock {\em Annals of K-theory}, 4(3):505--519, 2019.

\bibitem[Ant25]{antieau2025filtrations}
Benjamin Antieau.
\newblock Filtrations and cohomology {I}: crystallization.
\newblock {\em arXiv preprint arXiv:2511.01567}, 2025.

\bibitem[BBJ19]{brav2019darboux}
Christopher Brav, Vittoria Bussi, and Dominic Joyce.
\newblock A {D}arboux theorem for derived schemes with shifted symplectic
  structure.
\newblock {\em Journal of the American Mathematical Society}, 32(2):399--443,
  2019.

\bibitem[BCN21]{brantner2021pd}
Lukas Brantner, Ricardo Campos, and Joost Nuiten.
\newblock P{D} operads and explicit partition {L}ie algebras.
\newblock {\em arXiv preprint arXiv:2104.03870}, 2021.

\bibitem[Beh09]{behrend2009donaldson}
Kai Behrend.
\newblock Donaldson-{T}homas type invariants via microlocal geometry.
\newblock {\em Annals of Mathematics}, pages 1307--1338, 2009.

\bibitem[BG13]{bouaziz2013d}
E~Bouaziz and I~Grojnowski.
\newblock A d-shifted {D}arboux theorem.
\newblock {\em arXiv preprint arXiv:1309.2197}, 2013.

\bibitem[BGT13]{blumberg2013universal}
Andrew~J Blumberg, David Gepner, and Gon{\c{c}}alo Tabuada.
\newblock A universal characterization of higher algebraic {K}-theory.
\newblock {\em Geometry \& Topology}, 17(2):733--838, 2013.

\bibitem[BM19]{brantner2019deformation}
Lukas Brantner and Akhil Mathew.
\newblock Deformation theory and partition {L}ie algebras.
\newblock {\em arXiv preprint arXiv:1904.07352}, 2019.

\bibitem[BMN25]{brantner2025formal}
Lukas Brantner, Kirill Magidson, and Joost Nuiten.
\newblock Formal integration of derived foliations.
\newblock {\em arXiv preprint arXiv:2502.05257}, 2025.

\bibitem[Bou68]{bousfield1968operations}
Aldridge~Knight Bousfield.
\newblock Operations on derived functors of non-additive functors.
\newblock {\em (No Title)}, 1968.

\bibitem[BS17]{bhatt2017projectivity}
Bhargav Bhatt and Peter Scholze.
\newblock Projectivity of the {W}itt vector affine {G}rassmannian.
\newblock {\em Inventiones mathematicae}, 209(2):329--423, 2017.

\bibitem[Cho24]{chough2024formal}
Chang-Yeon Chough.
\newblock Formal {D}erived {A}lgebraic {G}eometry.
\newblock {\em arXiv preprint arXiv:2412.13506}, 2024.

\bibitem[Fu24]{fu2024duality}
Jiaqi Fu.
\newblock A duality between {L}ie algebroids and infinitesimal foliations.
\newblock {\em arXiv preprint arXiv:2410.04950}, 2024.

\bibitem[Fu25]{Fu2025}
Jiaqi Fu.
\newblock {\em Homotopical algebra of foliations}.
\newblock Ph{D} thesis, Universit{\'e} de Toulouse, 2025.
\newblock NNT: 2025TLSES087.

\bibitem[GP18]{gwilliam2018enhancing}
Owen Gwilliam and Dmitri Pavlov.
\newblock Enhancing the filtered derived category.
\newblock {\em Journal of Pure and Applied Algebra}, 222(11):3621--3674, 2018.

\bibitem[Gro67]{EGAIV4}
Alexander Grothendieck.
\newblock \'el\'ements de g\'eom\'etrie alg\'ebrique : Iv. \'etude locale des
  sch\'emas et des morphismes de sch\'emas, quatri\`eme partie.
\newblock {\em Publications Math\'ematiques de l'IH\'ES}, 32:5--361, 1967.

\bibitem[HHR24]{hennion2024gluing}
Benjamin Hennion, Julian Holstein, and Marco Robalo.
\newblock Gluing invariants of {D}onaldson--{T}homas type--{P}art {I}: the
  {D}arboux stack.
\newblock {\em arXiv preprint arXiv:2407.08471}, 2024.

\bibitem[Lur04]{lurie2004derived}
Jacob Lurie.
\newblock {\em Derived algebraic geometry}.
\newblock PhD thesis, Massachusetts Institute of Technology, 2004.

\bibitem[Lur09]{HTT}
Jacob Lurie.
\newblock {\em Higher topos theory}.
\newblock Princeton University Press, 2009.

\bibitem[{Lur}11]{Lur}
Jacob {Lurie}.
\newblock {D}{A}{G} {X}: {F}ormal {M}oduli {P}roblems, 2011.

\bibitem[Lur17]{HA}
Jacob Lurie.
\newblock Higher algebra, 2017.

\bibitem[Lur18]{SAG}
Jacob Lurie.
\newblock Spectral algebraic geometry.
\newblock {\em preprint}, 2018.

\bibitem[Mat16]{mathew2016galois}
Akhil Mathew.
\newblock The {G}alois group of a stable homotopy theory.
\newblock {\em Advances in Mathematics}, 291:403--541, 2016.

\bibitem[MRT22]{moulinos2022universal}
Tasos Moulinos, Marco Robalo, and Bertrand To{\"e}n.
\newblock A universal {H}ochschild--{K}ostant--{R}osenberg theorem.
\newblock {\em Geometry \& Topology}, 26(2):777--874, 2022.

\bibitem[Ogu75]{ogus1975cohomology}
Arthur Ogus.
\newblock Cohomology of the infinitesimal site.
\newblock In {\em Annales scientifiques de l'{\'E}cole Normale Sup{\'e}rieure},
  volume~8, pages 295--318, 1975.

\bibitem[Pri10]{pridham2010unifying}
Jon~P Pridham.
\newblock Unifying derived deformation theories.
\newblock {\em Advances in Mathematics}, 224(3):772--826, 2010.

\bibitem[PT15]{PT2015youtube}
Tony {Pantev} and Bertrand {To\"en}.
\newblock I{H}{E}{S} presentation.
\newblock \url{https://www.youtube.com/watch?v=fXYO552lXf4&t=199s}, 2015.

\bibitem[PTVV13]{pantev2013shifted}
Tony Pantev, Bertrand To{\"e}n, Michel Vaqui{\'e}, and Gabriele Vezzosi.
\newblock Shifted symplectic structures.
\newblock {\em Publications math{\'e}matiques de l'IH{\'E}S}, 117:271--328,
  2013.

\bibitem[Rak20]{raksit2020hochschild}
Arpon Raksit.
\newblock Hochschild homology and the derived de {R}ham complex revisited.
\newblock {\em arXiv preprint arXiv:2007.02576}, 2020.

\bibitem[Smi06]{smith2006differential}
S~Paul Smith.
\newblock Differential operators on the affine and projective lines in
  characteristic $p> 0$.
\newblock In {\em S{\'e}minaire d'Alg{\`e}bre Paul Dubreil et Marie-Paule
  Malliavin: Proceedings, Paris 1985 (37{\`e}me Ann{\'e}e)}, pages 157--177.
  Springer, 2006.

\bibitem[sta]{stack}
The stack projects.
\newblock https://stacks.math.columbia.edu/.

\bibitem[Tho00]{thomas2000holomorphic}
Richard~P Thomas.
\newblock A holomorphic {C}asson invariant for {C}alabi-{Y}au 3-folds, and
  bundles on $ k3 $ fibrations.
\newblock {\em Journal of Differential Geometry}, 54(2):367--438, 2000.

\bibitem[To{\"e}06]{toen2006champs}
Bertrand To{\"e}n.
\newblock Champs affines.
\newblock {\em Selecta mathematica}, 12(1):39--134, 2006.

\bibitem[Tot18]{totaro2018hodge}
Burt Totaro.
\newblock Hodge theory of classifying stacks.
\newblock 2018.

\bibitem[TV07]{toen2007moduli}
Bertrand To{\"e}n and Michel Vaqui{\'e}.
\newblock Moduli of objects in dg-categories.
\newblock In {\em Annales scientifiques de l'Ecole normale sup{\'e}rieure},
  volume~40, pages 387--444, 2007.

\bibitem[TV08]{toen2008homotopical}
Bertrand To{\"e}n and Gabriele Vezzosi.
\newblock {\em Homotopical Algebraic Geometry II: Geometric Stacks and
  Applications: Geometric Stacks and Applications}, volume~2.
\newblock American Mathematical Soc., 2008.

\bibitem[TV11]{toen2011algebres}
Bertrand To{\"e}n and Gabriele Vezzosi.
\newblock Algebres simpliciales ${S}^1$-{\'e}quivariantes, th{\'e}orie de de
  {R}ham et th{\'e}oremes {H}{K}{R} multiplicatifs.
\newblock {\em Compositio Mathematica}, 147(6):1979--2000, 2011.

\bibitem[TV15]{toen2015caracteres}
Bertrand To{\"e}n and Gabriele Vezzosi.
\newblock Caracteres de {C}hern, traces {\'e}quivariantes et g{\'e}om{\'e}trie
  alg{\'e}brique d{\'e}riv{\'e}e.
\newblock {\em Selecta Mathematica}, 21(2):449--554, 2015.

\bibitem[TV23a]{toen2023algebraic}
Bertrand To{\"e}n and Gabriele Vezzosi.
\newblock Algebraic foliations and derived geometry: the {R}iemann--{H}ilbert
  correspondence.
\newblock {\em Selecta Mathematica}, 29(1):5, 2023.

\bibitem[TV23b]{toen2023infinitesimal}
Betrand To{\"e}n and Gabriele Vezzosi.
\newblock Infinitesimal derived foliations.
\newblock {\em arXiv preprint arXiv:2305.13010}, 2023.

\end{thebibliography}
		\begin{minipage}{.45\linewidth}
		\begin{flushleft}
			\textsc{Université de Toulouse\\
				Institut de Mathématiques de Toulouse\\
				118, route de Narbonne\\
				F-31062 Toulouse Cedex 9\vspace{1em}}
		\end{flushleft}
	\end{minipage}\hfill\begin{minipage}{.45\linewidth}
		\begin{flushright}
			\textit{E-mail}: jiaqi.fu@univ-tlse3.fr\\
			\textit{Homepage}: \url{https://jiaqifumath.github.io/DAG/}
		\end{flushright}
	\end{minipage}

\end{document}